\documentclass{amsart}
\usepackage{stmaryrd}
\usepackage[all]{xy}
\usepackage{mathrsfs}
\usepackage{amscd,amssymb,color}

\newcommand{\ssdot}{\bullet}

\newcommand{\subdot}{_\ssdot}

\newcommand{\Spdot}[1][\ssdot]{S'_{#1}}
\newcommand{\Sdot}[1][\ssdot]{S_{#1}}

\newcommand{\Spdotmac}[2]{S'^{(#1)}_{#2}}

\newcommand{\Spdotq}[1][\seqdot]{\Spdotmac{q}{#1}}
\newcommand{\iSdot}[1][\ssdot]{S^{\infty}_{#1}}
\newcommand{\iSdotq}[1][\seqdot]{(S^{\infty})^{(q)}_{#1}}

\mathchardef\varDelta="7101

\let\sma\wedge
\newcommand{\htp}{\simeq}
\renewcommand{\to}{\mathchoice{\longrightarrow}{\rightarrow}{\rightarrow}{\rightarrow}}

\newcommand{\cA}{{\mathcal A}}
\newcommand{\cB}{{\mathcal B}}
\newcommand{\cC}{{\mathcal C}}
\newcommand{\cD}{{\mathcal D}}
\newcommand{\cE}{{\mathcal E}}
\newcommand{\cF}{{\mathcal F}}
\newcommand{\cM}{{\mathcal M}}
\newcommand{\cS}{{\mathcal S}}

\newcommand{\cU}{{\mathcal U}}

\let\catsymbfont\mathcal
\newcommand{\aA}{{\catsymbfont{A}}}
\newcommand{\aB}{{\catsymbfont{B}}}
\newcommand{\aC}{{\catsymbfont{C}}}
\newcommand{\aD}{{\catsymbfont{D}}}
\newcommand{\aE}{{\catsymbfont{E}}}

\newcommand{\aM}{{\catsymbfont{M}}}

\newcommand{\aP}{{\catsymbfont{P}}}

\newcommand{\aS}{{\catsymbfont{S}}}
\newcommand{\aT}{{\catsymbfont{T}}}
\newcommand{\aU}{{\catsymbfont{U}}}
\newcommand{\aV}{{\catsymbfont{V}}}

\newcommand{\A}{\mathcal{A}}

\newcommand{\C}{\catsymbfont{C}}
\renewcommand{\L}{\mathrm{L}}

\newcommand{\T}{\mathcal{T}}
\newcommand{\bbK}{I\mspace{-6.mu}K}

\newcommand{\bbZ}{\mathbb{Z}}
\newcommand{\bbN}{\mathbb{N}}

\newcommand{\bL}{\mathbb{L}}

\newcommand{\bN}{{\mathbb{N}}}

\newcommand{\bS}{\mathbb{S}}

\newcommand{\bZ}{{\mathbb{Z}}}

\def\quickop#1{\expandafter\DeclareMathOperator\csname
#1\endcsname{#1}}
\quickop{id}\quickop{Id}\quickop{colim}\quickop{hocolim}\quickop{op}
\quickop{Ar}\quickop{sing}\quickop{Hom}\quickop{w}\quickop{Ho}
\quickop{ob}\quickop{diag}\quickop{Stab}\quickop{Cat}\quickop{Mot}\quickop{Mod}
\quickop{map}\quickop{Cone}\quickop{End}\quickop{Idem}\quickop{Perf}\quickop{Ind}
\quickop{Gap}\quickop{Shv}\quickop{Spaces}\quickop{cof}\quickop{ex}\quickop{perf}
\quickop{tri}\quickop{Set}\quickop{Fib}\quickop{Nat}\quickop{haut}\quickop{holim}\quickop{Tor}\quickop{Ext}\quickop{cf}\quickop{cyc}

\newcommand{\lkan}[1]{(#1)_{\mathrm{iso}}}
\newcommand{\N}{\mathrm{N}}
\newcommand{\modules}[1]{\mathrm{Mod}(#1)}

\newcommand{\flt}{\mathrm{flt}}

\newcommand{\Pre}{\mathrm{Pre}}
\newcommand{\PPre}{\mathrm{Pre}}
\newcommand{\icat}{\Cat_\i}

\newcommand{\stabcat}{\Cat_\i^{\ex}}
\newcommand{\kstabcat}{\Cat_\i^{\ex(\!\kappa)}}
\newcommand{\idemstabcat}{\Cat_\i^{\perf}}

\newcommand{\idemtimes}{\widehat{\otimes}}
\newcommand{\idemfun}{\Fun^{\ex}}
\newcommand{\prcat}{{\mathcal{P}\mathrm{r}}^\mathrm{L}}
\newcommand{\stabprcat}{{{\mathcal{P}\mathrm{r}}^{\mathrm{L}}_\mathrm{St}}}
\newcommand{\stabprcatcg}{{{\mathcal{P}\mathrm{r}}^{\mathrm{L}}_\mathrm{St,cg}}}

\newcommand{\Loc}{\mathrm{Loc}}
\newcommand{\Fun}{\mathrm{Fun}} 
\newcommand{\Map}{\mathrm{Map}}
\newcommand{\bbS}{\mathbb{S}}
\newcommand{\bbA}{\mathbb{A}}
\newcommand{\ispec}{\aS_{\infty}}
\newcommand{\Spec}{\mathrm{Sp}}
\newcommand{\Spcat}{\Cat_\mathcal{S}}
\newcommand{\Spt}{\mathcal{S}}
\newcommand{\Scat}{\Cat_{\mathcal{T}}}

\newcommand{\sset}{\Set_{\Delta}}

\newcommand{\ie}{{i.e.,}\ }

\newcommand{\too}{\longrightarrow}

\renewcommand{\i}{\infty}
\newcommand{\Sp}{\mathrm{Sp}}
\newcommand{\rep}{\mathrm{rep}} 

\newcommand{\Motadd}{\cM_{\mathrm{add}}}

\newcommand{\Motaddun}{\cM_{\mathrm{add}}^{\mathrm{un}}}
\newcommand{\uMotadd}{\underline{\cM_{\mathrm{add}}^\kappa}}
\newcommand{\uwMotloc}{\underline{\cM_{\mathrm{wloc}}^\kappa}}
\newcommand{\wMotloc}{\cM_{\mathrm{loc}}^\kappa}
\newcommand{\Motloc}{\cM_\mathrm{loc}}

\newcommand{\Umot}{\cU_{\mathrm{add}}}

\newcommand{\Umotun}{\cU_{\mathrm{add}}^{\mathrm{un}}}
\newcommand{\uUmot}{\underline{\cU_{\mathrm{add}}^\kappa}}
\newcommand{\uwUloc}{\underline{\cU_{\mathrm{wloc}}^\kappa}}
\newcommand{\wUloc}{\cU_{\mathrm{loc}}^\kappa}
\newcommand{\Uloc}{\cU_{\mathrm{loc}}}
\newcommand{\Uadd}{\cU_{\mathrm{add}}}

\newcommand{\spect}{\Upsilon}
\newcommand{\kcone}{\cF_{\kappa}}
\newcommand{\ksusp}{\Sigma_{\kappa}}

\numberwithin{equation}{section}

\newtheorem{theorem}[equation]{Theorem}
\newtheorem*{theorem*}{Theorem}
\newtheorem{corollary}[equation]{Corollary}
\newtheorem{lemma}[equation]{Lemma}
\newtheorem{proposition}[equation]{Proposition}
\theoremstyle{definition}
\newtheorem{definition}[equation]{Definition}
\newtheorem{remark}[equation]{Remark}
\newtheorem{example}[equation]{Example}
\newtheorem{notation}[equation]{Notation}

\xyoption{arrow}
\xyoption{matrix}
\xyoption{cmtip}
\SelectTips{cm}{}

\newdir{ >}{{}*!/-5pt/\dir{>}}

\begin{document}

\title[A universal characterization of higher algebraic $K$-theory]{A universal characterization of \\higher algebraic $K$-theory}

\author{Andrew J. Blumberg}
  \address{Department of Mathematics, University of Texas,
Austin, TX \ 78703, USA}
\email{blumberg@math.utexas.edu}
\author{David Gepner}
\address{Fakult\"at f\"ur Mathematik,
Universit\"at Regensburg, 93040 Regensburg, Germany}
\email{djgepner@gmail.com}
\author{Gon{\c c}alo~Tabuada}
\address{Gon{\c c}alo Tabuada, Department of Mathematics, MIT, Cambridge, MA 02139, USA and Departamento de Matem{\'a}tica e CMA, FCT-UNL, Quinta da Torre, 2829-516, Caparica, Portugal}
\email{tabuada@math.mit.edu}

\begin{abstract}
In this paper we establish a universal characterization of higher
algebraic $K$-theory in the setting of small stable $\i$-categories.
Specifically, we prove that connective algebraic $K$-theory is the
universal {\em additive} invariant, \ie the universal functor with
values in spectra which inverts Morita equivalences, preserves
filtered colimits, and satisfies Waldhausen's additivity theorem.
Similarly, we prove that non-connective algebraic $K$-theory is the
universal {\em localizing} invariant, \ie the universal functor that
moreover satisfies the ``Thomason-Trobaugh-Neeman'' localization
theorem.

To prove these results, we construct and study two stable
$\i$-categories of ``non-commutative motives''; one associated to
additivity and another to localization.  In these stable
$\i$-categories, Waldhausen's $\Sdot$ construction corresponds to the
suspension functor and connective and non-connective algebraic
$K$-theory spectra become corepresentable by
the non-commutative motive of the sphere spectrum. In
particular, the algebraic $K$-theory of every scheme, stack, and ring
spectrum can be recovered from these categories of non-commutative
motives. In the case of a connective ring spectrum $R$, we
prove moreover that its negative $K$-groups are isomorphic to the negative $K$-groups of the ordinary ring $\pi_0R$.

In order to work with these categories of non-commutative motives, we
establish comparison theorems between the category of spectral
categories localized at the Morita equivalences and the category
of small idempotent-complete stable $\i$-categories.  We also explain
in detail the comparison between the $\i$-categorical version of
Waldhausen $K$-theory and the classical definition.

As an application of our theory, we obtain a complete classification
of the natural transformations from higher algebraic $K$-theory to
topological Hochschild homology ($THH$) and topological cyclic
homology ($TC$).  Notably, we obtain an elegant conceptual description
of the cyclotomic trace map.
\end{abstract}

\maketitle
\setcounter{tocdepth}{2}
\tableofcontents

\section{Introduction}


Algebraic $K$-theory is a fundamental algebro-geometric invariant,
capturing information in arithmetic, algebraic geometry, and topology.
The algebraic $K$-theory of a ring encodes many of its classical
number-theoretic invariants, such as its Picard and Brauer groups.
More generally, the algebraic $K$-theory of a scheme encodes
arithmetic information as well as information about its singularities.
The extension of algebraic $K$-theory from classical algebraic objects
to ring spectra and to derived schemes provides a connection to
geometric topology; notably, Waldhausen's $A$-theory (i.e., the
$K$-theory of the sphere spectrum) is essentially equivalent to stable
pseudo-isotopy theory \cite{WaldA2}.

The subject originated with Grothendieck's definition of $K_0$ (the
``Grothendieck group'') in the course of his work on the Riemann-Roch
theorem.  By construction, $K_0$ is the universal receptacle for Euler
characteristics, i.e. functions $\chi$ from the set of isomorphism classes of objects of a category $\C$ equipped with a suitable notion of ``equivalence'' and ``exact sequence'' to abelian groups which satisfy the relation
\[
\chi(X) - \chi(Y) + \chi(Z) = 0
\]
whenever there is an exact sequence $X\to Y\to Z$ in $\C$.

During the 70's and early 80's, Quillen~\cite{Quillen} and then
Waldhausen~\cite{Wald} extended Grothendieck's work making use of
tools from algebraic topology.  They defined the {\em connective}
algebraic $K$-theory spectrum $K(\aC)$ of a suitable category $\aC$; the homotopy
groups of this spectrum are the higher algebraic $K$-groups $K_i,
i \geq 0$.  Subsequently, Thomason and Trobaugh~\cite{TT} generalized
Bass' work \cite{Bass} on negative algebraic $K$-groups and introduced
the {\em non-connective} algebraic $K$-theory spectrum $\bbK(\aC)$ of $\aC$ in order to 
properly capture Mayer-Vietoris phenomena for schemes.

In contrast with the Grothendieck group, however, the construction
of these algebraic $K$-theory spectra does not provide a universal
characterization.
The most basic theorem about connective algebraic $K$-theory, the
additivity theorem \cite{Wald}, essentially says that the
$\Sdot$ construction forces the chosen cofiber sequences to split.
McCarthy's simplicial proof of the additivity theorem \cite{McCarthy2}
for any theory satisfying a few simple axioms suggests that the
$\Sdot$ construction is a universal construction for imposing
additivity on a functor from categories to spectra.  Moreover, the
construction of the connective algebraic $K$-theory spectrum
in terms of iterating the $\Sdot$ construction (along with the
additivity theorem applied to a simplicial path fibration) implies
that the $\Sdot$ construction functions as a kind of delooping
functor.

However, although these shadows of a universal description have been
known for a long time, a precise universal characterization of
algebraic $K$-theory proved elusive.  A major technical impediment has
been the absence of a framework in which to systematically express
homotopical constructions in the category of categories.  This
impediment has been lifted by recent
developments in the foundations of higher category theory, such as the theory of
derivators or $\i$-categories.

Over the past few years the third author of this paper and Cisinski
have carried out a program~\cite{CisTab0,CisTab, Duke} of
providing universal characterizations of algebraic $K$-theory in the
setting of dg-categories via the formalism of derivators.  In this
paper we adapt this approach to the setting of Lurie's 
theory of stable $\i$-categories.  The $\i$-category of stable
$\i$-categories provides a natural home for many examples of
interest coming from algebraic geometry and algebraic topology: the
$\i$-category of perfect complexes associated to a scheme (or a
stack), the $\i$-category of compact module spectra for a
ring spectrum, and the $\i$-category of stable retractive spaces are
all examples of stable $\i$-categories.

\subsection{Universal characterization}
Let $\stabcat$ be the $\i$-category of small stable $\i$-categories
and $\idemstabcat$ the full subcategory of {\em idempotent-complete}
small stable $\i$-categories; see \S\ref{sub:categories}. An exact
functor $F: \cA \to \cB$ in $\stabcat$ is called a {\em Morita
equivalence} when its idempotent completion
$\Idem(F): \Idem(\cA) \to \Idem(\cB)$ is an equivalence of
$\i$-categories; see definition~\ref{def:imorita}.  This is equivalent
to the condition that the induced map $F_!:\Mod(\aA)\to\Mod(\aB)$ is an equivalence of $\i$-categories.

A sequence $\aA \to \aB \to \aC$ in $\idemstabcat$ is called
{\em exact} if the composite is zero, $\aA \to \aB$ is fully faithful, and the map
$\aB/\aA\to\aC$, from the cofiber of the inclusion of $\aA$ into $\aB$
to $\aC$, is an equivalence; see proposition~\ref{prop:stabquot}. The
sequence is called {\em split-exact} if there exist adjoint splitting
maps such that the relevant composites are the respective identities;
see definition~\ref{def:splitexact}.  More generally, a sequence
$\aA\to\aB\to\aC$ in $\stabcat$ is (split-) exact if
$\Idem(\aA)\to\Idem(\aB)\to\Idem(\aC)$ is (split-) exact.
Now, let
\[
E: \stabcat \too \cD
\]
be a functor with values in a stable presentable $\i$-category $\cD$.
We say that $E$ is an {\em additive invariant} if it inverts Morita
equivalences, preserves filtered colimits, and sends split-exact
sequences to (split) cofiber sequences; see
definition~\ref{def:additive}.  This last condition corresponds to the
stable $\i$-categorical analogue of Waldhausen's additivity theorem.  If $E$
in fact sends {\em all} exact sequences to cofiber sequences we say
that it is a {\em localizing invariant}; this property corresponds to
Neeman's generalization of the Thomason-Trobaugh localization
theorem \cite{NeemanLoc, TT}.  

Every localizing invariant is additive, but not every additive
invariant is necessarily localizing.  Examples of localizing
invariants include the non-connective version of
algebraic $K$-theory and topological Hochschild homology.  Connective
algebraic $K$-theory is an example of an additive invariant which is
not localizing.  Note that, as we discuss in
sections~\ref{sec:kthy} and \ref{sec:non-connective}, the
$\i$-categorical versions of these invariants agree with the
``classical'' definitions.  For instance, we give a precise comparison
in section~\ref{sec:kthy} between Waldhausen's algebraic $K$-theory of
a Waldhausen category $\aC$ and the $\i$-categorical algebraic
$K$-theory of the $\i$-category underlying $\aC$.

Our first main result is the following.

\begin{theorem}{(see theorems~\ref{thm:main1}
and \ref{thm:univ-loc})}\label{thm:Imain1}
There are stable presentable $\i$-categories $\Motadd$ and $\Motloc$
and {\em universal} additive and localizing invariants 
\begin{eqnarray}\label{eq:univ-add/loc}
\Umot \colon \stabcat \too \Motadd && \Uloc \colon \stabcat \too \Motloc\,.
\end{eqnarray}
That is, given any stable presentable $\i$-category
$\aD$, we have induced equivalences of $\i$-categories
\begin{eqnarray*}
(\Umot)^*\colon\Fun^{\L}(\Motadd,\aD)&\stackrel{\sim}{\too}& \Fun_{\mathrm{add}}(\stabcat,\aD)\\
(\Uloc)^*\colon\Fun^{\L}(\Motloc,\aD)&\stackrel{\sim}{\too}& \Fun_{\mathrm{loc}}(\stabcat,\aD)\,,
\end{eqnarray*}
where the left-hand sides denote the $\i$-categories of
colimit-preserving functors and the right-hand sides the
$\i$-categories of additive and localizing invariants.
\end{theorem}

From a motivic perspective, the $\i$-categories $\Motadd$ and
$\Motloc$ should be considered as candidate categories of
non-commutative motives.  In fact, theorem~\ref{thm:Imain1} shows us
that every additive (respectively localizing) invariant factors
uniquely through $\Motadd$ (resp., through $\Motloc$).  That is, all the
information concerning additive (resp., localizing) invariants is encoded in
$\Motadd$ (resp., in $\Motloc$).

Our second main result is the following characterization of the higher
algebraic $K$-theory of a stable $\i$-category.  Any
stable $\i$-category (and in particular $\Motadd$ and $\Motloc$)
admits natural mapping spectra; that is, a stable $\i$-category is
naturally enriched over the $\i$-category $\ispec$ of spectra (see
sections~\ref{sec:stabilization} and ~\ref{sec:specinf}).

\begin{theorem}{(see theorems~\ref{thm:main2} and \ref{thm:main3})}\label{thm:Imain2}
Let $\aA$ be an idempotent-complete small stable $\i$-category. Then,
there are natural equivalences of spectra   
\begin{eqnarray}
\Map(\Umot(\ispec^{\omega}), \Umot(\aA)) &\simeq& K(\aA) \label{eq:1} \\
\Map(\Uloc(\ispec^{\omega}), \Uloc(\aA)) &\simeq& \bbK(\aA)  \label{eq:2} \,,
\end{eqnarray}
where $\ispec^{\omega}$ is the small stable $\i$-category of
compact spectra. In particular, for all $n\in\bbZ$, we have isomorphisms of abelian groups
\begin{align}
\Hom(\Umot(\ispec^{\omega})), \Sigma^{-n}\Umot(\cA)) &\simeq K_n(\cA)\\
\Hom(\Uloc\,(\ispec^{\omega}))\,, \Sigma^{-n}\Uloc\,(\cA)) &\simeq \bbK_n(\cA)
\end{align}
in the triangulated categories $\Ho(\Motadd)$ and $\Ho(\Motloc)$. 
\end{theorem}

\begin{remark}\label{rk:stronger}
In fact, stronger results are true.  In equivalence \eqref{eq:1},
$\ispec^\omega$ can be replaced by any compact idempotent-complete 
small stable $\i$-category $\aB$ and the right-hand side by the
$K$-theory spectrum of $\Fun^{\ex}(\aB,\aA)$; see
theorem~\ref{thm:main2}. In equivalence \eqref{eq:2}, $\ispec^{\omega}$
can be replaced by any smooth and proper (i.e., dualizable) small
stable $\i$-category $\aB$ and the right-hand side by
$\bbK(\aB^{\op}  \idemtimes \aA)$; see theorem~\ref{thm:main34gen}.
\end{remark}

In particular, when $\cA$ is the $\i$-category of perfect complexes
over a suitable scheme (or stack), we recover the $K$-theory spectra of the
scheme (stack).  Taking $\cA$ to be the $\i$-category of compact
modules over a ring spectrum $R$, we recover the $K$-theory of $R$.
When $R$ is a connective ring spectrum, we show in
theorem~\ref{thm:nonconnKconn} that the negative homotopy groups of the non-connective $K$-theory of $R$ are the same as those of $\pi_0R$.
However, we expect the non-connective $K$-theory of a non-connective
ring spectrum to be an interesting new invariant.

Note that the left-hand sides of the natural equivalences and
isomorphisms of theorem~\ref{thm:Imain2} are defined solely in terms
of universal constructions on presheaf categories; algebraic
$K$-theory is not used in their construction.  Rather, a variant of
Waldhausen's path-fibration argument (see
proposition~\ref{prop:suspensions}) shows that Waldhausen's $\Sdot$ construction acts as the suspension functor
in $\Motadd$ and also $\Motloc$.

Therefore, theorem~\ref{thm:Imain2} (combined with
theorem~\ref{thm:Imain1}) provides an intrinsic characterization of
algebraic $K$-theory as a functor of stable $\i$-categories.
Furthermore, theorem~\ref{thm:inats} below (see also
theorem~\ref{thm:nats} in the text) shows that the co-representability
result coupled with the Yoneda lemma provides a complete
classification of all natural transformations from
algebraic $K$-theory to an arbitrary additive (or localizing) functor
from small stable categories to spectra.

\begin{remark}
Analogues of theorems~\ref{thm:Imain1}, \ref{thm:Imain2},
and \ref{thm:inats} in the setting of dg-categories were previously
known (due to the third author and
Cisinski~\cite{CisTab0,CisTab,Duke}).  Our arguments follow the
same general outline.  
\end{remark}

\subsection{Morita theory}

The main technical device in our proofs of theorems~\ref{thm:Imain1}
and \ref{thm:Imain2} is the Morita theory of stable categories and
spectral categories.  In particular, we prove theorem~\ref{thm:Imain2}
by using a comparison result between the theory of small spectral
categories (see \S\ref{sub:spectralR}) and the theory of small stable
$\i$-categories to rigidify questions about the algebraic $K$-theory
of $\i$-categories to corresponding questions in the (classical)
Waldhausen $K$-theory of Waldhausen categories.

The category $\Spcat$ of small spectral categories
carries a Quillen model category structure in which the weak
equivalences are the {\em DK-equivalences}, i.e., the functors that
are fully faithful and essentially surjective up to weak homotopy
equivalence; see \cite{Spectral} (reprised below in
theorem~\ref{thm:Quillenstable}).  As a consequence, we can form the
associated $\i$-category $(\Spcat)_{\i}$ of small spectral
categories.

A spectral functor $F\colon \aA \to \aB$ is called a {\em triangulated
equivalence} if it induces a weak equivalence on the triangulated
closures of $\aA$ and $\aB$, and it is called a {\em Morita
equivalence} if it induces a weak equivalence on the thick closures of
$\aA$ and $\aB$; see definition~\ref{def:Morita}.  Our comparison
result, which can be regarded as a generalization of the
Morita theory of \cite{SS}, is the following.

\begin{theorem}{(see theorems~\ref{thm:stable} and
    \ref{thm:morita})}\label{thm:models}
The accessible localization of $(\Spcat)_{\i}$ along the triangulated equivalences
is equivalent to $\stabcat$, and the (further) localization of
$(\Spcat)_{\i}$ along the Morita equivalences is equivalent to
$\idemstabcat$.
\end{theorem}

We use this comparison result to deduce structural properties of the
categories $\stabcat$ and $\idemstabcat$, notably that they are compactly
generated, complete, and cocomplete; see corollary~\ref{cor:cgen}.
Furthermore, we prove theorem~\ref{thm:Imain2} by using
theorem~\ref{thm:models} to lift split-exact sequences of small
$\i$-categories to split-exact sequences of small spectral categories
so we can take the $K$-theory in the setting of
Waldhausen categories.  More generally, this comparison result
explains the relationship between the classical versions of algebraic
$K$-theory (and topological Hochschild and cyclic homology) and the
$\i$-categorical versions.
We believe that theorem~\ref{thm:models} is of independent interest
and expect it will find applications in the future.  For
instance, this theorem provides clean and concise proofs of the
main theorems of To{\"e}n's work \cite{Toen} on internal $\hom$
objects in the category of dg-categories and its (previously
unknown) extension to the context of spectral categories.

\subsection{Symmetric monoidal structure and dualizable objects}

The category $\idemstabcat$ is a symmetric monoidal $\i$-category 
(in the sense of \cite[\S 2]{HAG}), in which the tensor product
$\idemtimes$ is characterized by the property that maps out of
$\aA \idemtimes \aB$ correspond to maps out of the product
$\aA \times \aB$ which preserve colimits in each variable; see
section~\ref{sec:monoidal} for a discussion of this structure,
following the work of Lurie \cite{HAG} and Ben-Zvi, Francis, and Nadler \cite{BFN}.
We will reserve a careful study of the structure of $\Motadd$ and
$\Motloc$ as symmetric monoidal $\i$-categories for the
forthcoming paper \cite{BGT2}. However, in order to carry out the
extension of the non-connective co-representability theorem described
in remark~\ref{rk:stronger}, we study the theory of dualizable objects
in $\idemstabcat$, using the theory of \cite[\S 4.2.5]{HAG}.

In analogy with the situation for dg-categories \cite[\S4]{CisTab},
we obtain a characterization of the dualizable objects in
$\idemstabcat$ as the {\em smooth} and {\em proper} objects.  We
define these notions as follows.  Implicit in the comparison between
small stable $\i$-categories and spectral categories of
theorem~\ref{thm:models} is the fact that for objects $a$ and $b$ in a
small stable $\i$-category $\aA$ there exists a natural mapping
spectrum $\aA(a,b)$ (see definition~\ref{defn:mappingspectrum}).
Using this fact, we say that a small stable $\infty$-category $\aA$ is 
{\em proper} if, for all pairs of objects $a$ and $b$ of $\aA$, the mapping spectrum $\aA(a,b)$ is compact.  We say that a small stable $\infty$-category $\aA$ is {\em smooth} if it is perfect as an $\aA^{\op}\idemtimes\aA$-module.  (Here we use the fact that any small
stable $\i$-category can be regarded as a bimodule over itself.)  We
then have the following theorem characterizing the dualizable objects
in these terms:

\begin{theorem}{(see theorem~\ref{thm:dualizable})}
An idempotent-complete small stable $\infty$-category $\aA$ is
dualizable (as an object of the symmetric monoidal $\infty$-category
$\idemstabcat$) if and only if $\aA$ is smooth and proper.
Moreover, the dual of a dualizable idempotent-complete small stable $\infty$-category $\aA$ is the opposite $\infty$-category $\aA^{\op}$.
\end{theorem}

\subsection{Trace methods}

One of the major revolutions in the calculational study of algebraic
$K$-theory of rings and schemes in the past two decades has been the
development of trace methods, following the ideas of Goodwillie and
Bokstedt-Hsiang-Madsen \cite{BokstedtHsiangMadsen}.  The cyclotomic
trace from $K$-theory to topological cyclic homology $TC$ and
topological Hochschild homology $THH$ (stable homotopy theory
generalizations of negative cyclic homology and Hochschild homology)
has allowed major calculational advances.  The fiber of this map is
well understood by work of Goodwillie, McCarthy, and
Dundas \cite{McCarthyFib, DundasFib}, and the target is relatively
computable using the methods of equivariant stable homotopy theory
(e.g., see the extensive body of work by Hesselholt and Madsen on the
Quillen-Lichtenbaum conjecture \cite{HesselholtMadsen}).  One
application of the co-representability of algebraic $K$-theorem
(theorem~\ref{thm:Imain2}) is the complete
classification of all natural transformations with source the
algebraic $K$-theory functor.

\begin{theorem}(see theorem~\ref{thm:nats})\label{thm:inats}
Given an additive invariant
\[
E: \stabcat \too \ispec
\]
with values in the stable $\i$-category of spectra, we have a natural
equivalence of spectra
\[
\Map(K,E) \simeq E(\ispec^{\omega}),
\]
where $\Map(K,E)$ denotes the spectrum of natural transformations of additive invariants.  The
analogous result for localizing invariants holds.  In the particular
case where $E$ is topological Hochschild homology, we obtain an
isomorphism 
\[
\pi_0\Map(K, THH) \simeq \pi_0 THH(\ispec^{\omega}) \simeq \pi_0
THH(\bbS) \simeq \bZ.
\]
\end{theorem}

A calculation then provides a canonical construction and conceptual
description of the topological Dennis trace map $K \to THH$.

\begin{corollary}{(see section~\ref{sec:cyclo})}
The set of homotopy classes of natural transformations of additive
invariants from connective algebraic $K$-theory to $THH$ is isomorphic
to $\bZ$; furthermore, the topological Dennis trace is
characterized up to homotopy as the natural transformation $K \to THH$
corresponding to the unit $1 \in \bZ$.
\end{corollary}

That is, up to scaling, the trace is the only natural
transformation of additive invariants between connective algebraic
$K$-theory and $THH$.  This provides a direct proof that all known
constructions of the topological Dennis trace map agree up to
homotopy.

Working directly with topological cyclic homology ($TC$) is somewhat
more complicated; $TC$ does not preserve filtered colimits in general,
and is therefore not an additive or localizing functor.  However, we
deduce an analogous identification of the cyclotomic trace as
determined by the unit map; see section~\ref{sec:cyclo}.

Finally, we note that in the localizing setting our results provide an
extension of the cyclotomic trace from non-connective algebraic
$K$-theory to the non-connective versions of $TC$ and $THH$.  This
generalizes and extends the non-connective traces constructed
in \cite{GeisserHesselholtRelative} and \cite{BM} for rings and
schemes.

\subsection{Related works}
The ``motivic'' idea of constructing universal invariants is not new
and appears in several different subjects: for example,
Corti{\~n}as-Thom's work~\cite{Cortinas} on bivariant algebraic
$K$-theory, Higson's work~\cite{Higson} on Kasparov's bivariant
$K$-theory, Meyer-Nest's work~\cite{Nest} on $C^{\ast}$-algebras,
Morel-Voevodsky's work~\cite{Morel} on $\bbA^1$-homotopy theory of
schemes, and Voevodsky's work~\cite{Voevodsky} on (mixed) motives.

In this vein, over the past few years the third author and Cisinski
have carried out a program~\cite{CisTab0,CisTab,Duke} of
providing universal characterizations of algebraic $K$-theory in the
setting of dg-categories, using the formalism of derivators.  One of the main goals of our
work in this paper completes this program by extending these results
to stable $\i$-categories and by solving a key open question left open
in the previous work of the third author, namely identifying the
cyclotomic trace in terms of the co-representability results.

Finally, we would like to mention that Barwick~\cite{Bar} has recent work on a
universal characterization of higher algebraic $K$-theory, in the
context of a detailed study of the algebraic $K$-theory of
$\i$-categories. 

\bigskip

\noindent\textbf{Acknowledgments:} The authors would like to thank
Mike Mandell for many helpful conversations and Haynes Miller for
asking motivating questions.  They are grateful to David Ben-Zvi,
Chris Brav, Bob Bruner, Jonathan Campbell, Denis-Charles Cisinski,
Bjorn Dundas, Tom Goodwillie, Jeremiah Heller, Kathryn Hess, John
Lind, Peter May, Jack Morava, Markus Spitzweck and Bertrand To\"en for
helpful comments on a previous draft.  The authors would like also to
thank the anonymous referees for very careful readings and detailed
comments and corrections which substantially improved this paper.
This project was initiated during a visit by the second and third
authors to Stanford's math department, and they would like to thank
the department for its hospitality.  Finally, the authors would like
to thank the Midwest Topology Network for funding various trips which
facilitated the conduct of this research.  A.~J.~Blumberg was
supported in part by NSF grant DMS-0906105.  G.~Tabuada was supported by the {\em Estimulo {\`a} Investiga{\c c}{\~a}o} Award 2008 - 
Calouste Gulbenkian Foundation, by the FCT-Portugal grants PTDC/MAT/098317/2008 and SFRH/BSAB/1116/2011 and by the NEC award-2742738.

\section{Spectral categories and stable $\i$-categories}

The purpose of this section is to collect and recall the results about
spectral categories and stable $\i$-categories we will require for our
constructions.

\subsection{Review of spectral categories}\label{sub:spectralR}
We write $\aT$ for the symmetric monoidal simplicial model category of simplicial sets and $\Spt$ for the symmetric monoidal simplicial model category
of symmetric spectra \cite{HSS}.  Recall that a {\em
spectral category $\cA$} is a category 
enriched in the category of symmetric spectra.  Specifically, a
spectral category is given by:
\begin{itemize}
\item A class of objects $\mbox{obj}(\cA)$,
\item for each pair of objects $(x,y)$ of $\cA$, a symmetric
spectrum $\cA(x,y)$,
\item for each triple of objects $(x,y,z)$ of $\cA$, a composition
morphism in $\Spt$ $$\cA(y,z)\wedge \cA(x,y) \too \cA(x,z)\,,$$
satisfying the usual associativity condition, and
\item for any object $x$ of $\cA$, a morphism $\bS \to \cA(x,x)$
in $\Spt$, satisfying the usual unit condition with respect to the
above composition.
\end{itemize}

A spectral category is said to be small if its class of objects forms
a set.  We write $\Spcat$ the category of small spectral categories
and spectral (enriched) functors.  References on spectral categories
are \cite[\S2]{BM}, \cite[Appendix\,A]{SS} and \cite[\S 2]{Spectral}.

We now briefly recall the Quillen model structure on spectral
categories we work with in this paper.  Given a spectral category
$\cA$, we can form a genuine category $[\cA]$ by keeping the same set
of objects and defining the set of morphisms between $x$ and $y$ in
$[\cA]$ to be the set of morphisms in the homotopy category $\Ho(\Spt)$
from the sphere spectrum $\bS$ to $\cA(x,y)$.  We obtain in this way a
functor $$ [-] \colon \Spcat \too \Cat \,,$$ with values in the
category of small categories.  Equivalently, we can think of $[-]$ as
computed by passing to $\pi_0$ on the morphism spectra, and so we will
also refer to $[\cA]$ as the homotopy category $\Ho(\cA)$.

\begin{definition}\label{def:stable}
A spectral functor $F:\cA \to \cB$ is a {\em DK-equivalence}, if:
\begin{itemize}
\item for all objects $x,y \in \cA$, the morphism in $\Spt$
$$ F(x,y):\cA(x,y) \too \cB(Fx,Fy) $$
is a stable equivalence and
\item the induced functor
$$ [F]:[\cA] \too [\cB]$$
is an equivalence of categories.
\end{itemize}
\end{definition}

\begin{theorem}{(\cite[5.10]{Spectral})}\label{thm:Quillenstable}
The category $\Spcat$ carries a right proper Quillen model structure
whose weak equivalences are the DK-equivalences.
\end{theorem}

Recall from \cite[\S\,2]{Spectral} the natural adjunction 
\begin{equation}\label{eqn:specsimp}
\xymatrix{
\Spcat \ar@<1ex>[d]^{\Omega^\i} \\
 \Scat \ar@<1ex>[u]^{\Sigma^{\infty}_+} \,,
}
\end{equation}
between spectral and simplicial categories, where $\Omega^\i$ (also
denoted $(-)_0$) is the space of maps from the unit (equivalently,
restriction to the $0$-th space of the spectrum).  We will use this
adjunction to pass between spectral categories and $\i$-categories.
Using the model structure on simplicial categories of \cite{Bergner},
the pair $(\Sigma^{\infty}_+, \Omega^\i)$ is a Quillen adjunction. 

For technical control, we require the following corollary which
sharpens the description of the model structure, providing a
combinatorial model category.  (For references for Jeff Smith's theory
of combinatorial model categories, see \cite{BarwickLoc}
or \cite{Dugger}.)

\begin{corollary}\label{cor:Quillencombi}
The category $\Spcat$ endowed with the model structure of
theorem~\ref{thm:Quillenstable} is a combinatorial model category and
is Quillen equivalent (via a zig-zag) to a simplicial 
category with a left proper combinatorial simplicial model structure.
There are simplicial cofibrant and fibrant replacement functors.  The
adjunction $(\Sigma^{\infty}_+, \Omega^\i)$ can be lifted to a
simplicial Quillen adjunction.
\end{corollary}

\begin{proof}
The proof of this theorem follows from a refinement of the proof of
Theorem~\ref{thm:Quillenstable}.  The model structure therein 
arises as the Bousfield localization of a cofibrantly-generated model
structure on $\Spcat$ in which the weak equivalences are the levelwise
equivalences~\cite[\S4]{Spectral}, i.e., the spectral functors
$F: \cA \to \cB$ such that for all objects $x, y \in \cA$, the
morphism $\cA(x,y) \to \cB(Fx,Fy)$ is a levelwise equivalence, and the
induced simplicial functor $\Omega^\i(\cA) \to \Omega^\i(\cB)$ is a
DK-equivalence.

First, we observe that the category $\Spcat$ is locally
presentable; a set of small generators is given by applying the
functor $U$ (see \cite[A.1]{Spectral}) to a set of small 
generators for the category of symmetric spectra. Since the leverwise model structure on $\Spcat$ is cofibrantly generated, it follows that it is combinatorial.  Next, the
arguments of \cite{Spectral} produce a generating set of
DK-equivalences at which to localize $\Spcat$.  The main
theorem about the existence of left Bousfield localization for
combinatorial model categories (e.g., see the treatment
in~\cite{BarwickLoc}) now implies that we can localize and obtain a
combinatorial model structure on $\Spcat$.

The machinery of Dugger's approach to universal homotopy
theories \cite{Dugger} now permits us to replace $\Spcat$ with a
Quillen equivalent simplicial model category (the simplicial objects
over $\Spcat$) which is combinatorial and left proper.  By applying
the techniques of \cite{Dugger, RezkSchwedeShipley}, we can promote
this adjunction to a simplicial Quillen adjunction.  Specifically, the
simplicial prolongation of the adjunction forms a Quillen pair on the
categories of simplicial objects \cite[6.1]{RezkSchwedeShipley}.
\end{proof}

Let $\cA$ be a (fixed) small spectral category and let $\cA^{\op}$ denote the
opposite spectral category, defined by $\cA^{\op}(x,y)=\cA(y,x)$.

\begin{definition}
A {\em $\cA$-module} is a spectral functor from $\cA^{\op}$ to the
spectral category $\Spt$ of symmetric spectra.  We denote by
$\widehat{\aA}$ the spectral category of $\cA$-modules.
\end{definition}

By \cite[A.1.1]{SS}, $\widehat{\aA}$ can be given a
combinatorial spectral model structure in which the weak equivalences
are the pointwise stable equivalences and the
fibrations are pointwise fibrations (referred to as the projective
model structure).  We will denote by $\widehat{\aA}^{\cf}$ the full
spectral subcategory of $\widehat{\cA}$ on the cofibrant and fibrant
$\cA$-modules, and by $\cD(\cA)$ the {\em derived category} of $\cA$,
i.e., the homotopy category $\Ho(\widehat{\cA})$ associated to the
model structure.  As usual, there is an equivalence
$[\widehat{\aA}^{\cf}]\simeq\cD(\cA)$.

Notice that we have a (fully faithful) spectral Yoneda embedding
$\cA \to \widehat{\cA}$ which sends the object $z$ to the functor
$\cA(-,z)\colon\cA^{\op}\to\Spt$ represented by $z$.  Note that when $\cA$ is
fibrant, the Yoneda embedding lands in $\widehat{\aA}^{\cf}$.
By \cite[\S A.1]{SS}, a spectral functor $F\colon\cA \to \cB$ gives rise to
a restriction/extension Quillen adjunction
\[
\xymatrix{
\widehat{\cB} \ar@<1ex>[d]^{F^{\ast}} \\
 \widehat{\cA} \ar@<1ex>[u]^{F_!}
}
\]
and therefore a total left-derived functor $\bL
F_! \colon \cD(\cA)\to\cD(\cB)$.

We will be most interested in spectral categories $\cA$ for which the
homotopy category $\Ho(\cA)$ has a triangulated structure compatible
with the mapping spectra; we refer to \cite[4.4]{BM} for the
definition of a pretriangulated spectral category, and highlight the
essential consequence \cite[4.6]{BM} that the homotopy category of a
pretriangulated spectral category is triangulated (with distinguished
triangles given by the Puppe sequences).  This is the stable homotopy
theory analogue of the notion of a pretriangulated dg-category.  A
spectral functor between pretriangulated spectral categories is a
DK-equivalence if and only if it induces an equivalent on homotopy
categories \cite[5.7]{BM2}.  Using the Yoneda embedding, we can
construct minimal pretriangulated categories containing the
spectral category $\aA$.

Given a spectral category $\aA$, the proof of \cite[4.5]{BM}
constructs a functorial ``triangulated closure''
$\widehat{\aA}_{\tri}$ which is a pretriangulated spectral category.
Briefly, $\widehat{\aA}_{\tri}$ consists of the subcategory of
cofibrant-fibrant objects in $\widehat{\aA}$ which have the homotopy
type of finite cell objects (in the projective model structure).
Using retracts of finite cell objects instead \cite[4.5]{BM} produces
a functorial ``thick closure'' $\widehat{\aA}_{\perf}$, which is an
idempotent-complete pretriangulated spectral category.  

\begin{remark}\label{rem:spectralcard}
In order for the preceding definitions to produce small spectral
categories, we need to restrict the sizes of the sets in the spaces of
the mapping spectra.  A careful discussion of this issue appears
in \cite[\S 4]{BM}; see also \cite{BM3}.  We return to the issue of
set-theoretic considerations in Section~\ref{sec:compact}.
\end{remark}

We start with a spectral functor $F \colon \aA \to \aB$, and tacitly
assume we have performed a functorial fibrant replacement.  We denote by 
$F_!^{\cf} \colon \widehat{\aA}^{\cf}\to\widehat{\aB}^{\cf}$ the
composite of $F_!$ with a fibrant replacement (note that we do not
need a cofibrant replacement here since $F_!$ preserves cofibrant
objects) to obtain 
\[ 
F_!^{\cf}\colon\widehat{\aA}^{\cf}\to\widehat{\aB}^{\cf}.
\]
Since this is a model of the derived functor of $F_!$ as a left
Quillen functor, it preserves homotopy colimits and thus 
sends modules of the homotopy type of finite cell $\aA$-modules to
modules of the homotopy type of finite cell $\aB$-modules and hence 
perfect $\aA$-modules to perfect $\aB$-modules.  Therefore, the
following definitions make sense.

\begin{definition}\label{def:Morita}
A spectral functor $F\colon \aA\to\aB$ is called

\begin{itemize}
\item a {\em triangulated equivalence} if
the induced functor
\[
F_!^{\cf}\colon \widehat{\aA}_{\tri} \to \widehat{\aB}_{\tri}
\]
is a DK-equivalence of spectral categories.

\item a {\em Morita equivalence} if the induced functor
\[
F_!^{\cf}\colon \widehat{\aA}_{\perf} \too \widehat{\aB}_{\perf}
\]
is a DK-equivalence of spectral categories.
\end{itemize}
\end{definition}

\begin{remark}
Suppose we are given a spectral functor $F\colon \aA\to\aB$.  Since
$\widehat{\aA}^{\cf}$ is generated by $\widehat{\aA}_{\perf}$ under
filtered homotopy colimits and Morita equivalences are stable under
filtered homotopy colimits, it follows that $F$ is a Morita
equivalence if and only if 
$F_!^{\cf}\colon\widehat{\aA}^{\cf}\to\widehat{\aB}^{\cf}$ is a
DK-equivalence.
\end{remark}

We can relate these notions to definitions purely on the level of
triangulated categories (the relationship between triangulated
constructions and enriched constructions is discussed further in
Section~\ref{sec:exact}).  For a spectral category $\aA$, let
$\cD_{\tri}(\aA)$ denote the smallest triangulated subcategory of
$\cD(\cA)$ containing the image of the $\aA$ under the Yoneda
embedding, and $\cD_{\perf}(\aA)$ denote the smallest thick
subcategory of $\cD(\cA)$ containing the image of $\aA$ under the
Yoneda embedding.  Observe that
$\cD_{\tri}(\aA)\simeq\Ho(\widehat{\aA}_{\tri})$ and
$\cD_{\perf}(\aA)\simeq\Ho(\widehat{\aA}_{\perf})$.  As a consequence,
we obtain the following proposition.

\begin{proposition}\label{prop:Morita}
A spectral functor $F\colon \aA \to \aB$ is

\begin{itemize}
\item a triangulated equivalence if and only if the induced derived functor
\[
\bL F_!\colon \cD_{\tri}(\aA) \to \cD_{\tri}(\aB)
\]
is an equivalence of (triangulated) categories,

\item a Morita equivalence if and only if the induced derived functor
\[
\bL F_!\colon \cD_{\perf}(\aA) \too \cD_{\perf}(\aB)
\]
is an equivalence of (triangulated) categories.
\end{itemize}
\end{proposition}

\begin{proof}
This follows immediately from \cite[5.7]{BM2}.
\end{proof} 

Finally, note that we can use $\Omega^\i$ to obtain simplicial models
of the triangulated and thick closures.  Define $\Mod(\aA)$ to be the
simplicial category $\Omega^\i \widehat{\aA}^{\cf}$, $\Mod(\aA)_{\tri}$ to
be the simplicial category $\Omega^\i \widehat{\aA}_{\tri}$, and
$\Mod(\aA)_{\perf}$ to be the simplicial category
$\Omega^\i \widehat{\aA}_{\perf}$.  Of course, it is also possible to
give intrinsic definitions of the latter two categories in terms of
$\Mod(\aA)$.

Summarizing the relationships between the various categories, we have
the following commutative diagram (with horizontal arrows induced by
the Yoneda embedding and subsequent inclusions):

\[
\xymatrix{
\aA \ar[r] \ar[d]^{\Omega^\i} & \widehat{\aA}_{\tri} \ar[r] \ar[d]^{\Omega^\i} &
\widehat{\aA}_{\perf} \ar[r]\ar[d]^{\Omega^\i} & \widehat{\aA}^{\cf} \ar[d]^{\Omega^\i} \\
\Omega^\i\aA \ar[r] \ar[d]^{\Ho} & \Mod(\aA)_{\tri} \ar[r] \ar[d]^{\Ho} &
\Mod(\aA)_{\perf} \ar[r]\ar[d]^{\Ho} & \Mod(\aA) \ar[d]^{\Ho} \\
\Ho(\aA) \ar[r] & \aD_{\tri}(\aA) \ar[r] & \aD_{\perf}(\aA) \ar[r] & \aD(\aA) \\
}
\]

\subsection{The $\i$-categories $\stabcat$ and
$\idemstabcat$}\label{sub:categories}

The basic setting for our work is the theory of $\i$-categories (and
particularly stable $\i$-categories), which provide a tractable way to
handle a ``homotopical category of homotopical categories'' as well as
homotopically meaningful categories of homotopical functors.  There
are now many competing models of $\i$-categories, including Rezk's
complete Segal spaces~\cite{Rezk}, the Segal
categories~\cite{HirschowitzSimpson, Tamsamani} of Simpson and
Tamsamani, the quasicategories (weak Kan complexes) of Boardman
and Vogt, the homotopy theory of simplicial categories as studied by
Dwyer-Kan and Bergner ~\cite{DwyerKan, Bergner}, and others, all of
which are known to be equivalent (see \cite{BergnerCompare} for a nice
discussion of the situation).  In a sense the situation is analogous
to the situation with the varied modern categories of spectra (e.g.,
symmetric spectra, orthogonal spectra, EKMM $S$-modules).  None of the
work of this paper depends in any way on particular properties of the
model of $\i$-categories chosen; given certain basic structural
properties, one could carry out our arguments in any of them.

We have chosen to work in this paper with the theory of
quasicategories.  These first appeared in the work of Boardman and
Vogt, where they were referred to as weak Kan
complexes~\cite{BoardmanVogt}.  The theory was subsequently developed
by Joyal~\cite{Joyal} and then extensively studied by Lurie.  In this
section we give a rapid review of the relevant background on the
theory of quasicategories as a model of $\i$-categories.  Our basic
references for this material are Lurie's books \cite{HTT, HAG}. 

We will write $\icat$ to denote the $\i$-category of small
$\i$-categories and functors, which we explicitly model as the
category of simplicial sets with the Joyal model
structure~\cite{Joyal}.  There is a simplicial nerve functor $\N$ from
simplicial categories to simplicial sets which is the right Quillen
functor of a Quillen equivalence~\cite[\S 1.1.5.5, 1.1.5.13,
2.2.5.1]{HTT}
\[
\xymatrix{
\Scat \ar@<1ex>[d]^N \\
 \Set_\Delta \ar@<1ex>[u]^{\mathfrak{C}}.
}
\]
Here the model structure on the top is Bergner-Dwyer-Kan's model
structure on simplicial categories \cite{Bergner} and the model
structure on the bottom is Joyal's model structure on simplicial sets.

There are a number of options for producing the ``underlying''
$\i$-category of a category equipped with a notion of ``weak
equivalence''. The most structured setting is that of a simplicial
model category $\aC$, where the $\i$-category can be obtained by
restricting to the full simplicial subcategory $\aC^{\cf}$ of
cofibrant-fibrant objects and then applying the simplicial nerve functor $\N$.
More generally, if $\aC$ is a category equipped with a subcategory of
weak equivalences $w\aC$, the Dwyer-Kan simplicial localization
$L \aC$ \cite{DwyerKan} provides a corresponding simplicial category,
and then $\N((L\aC)^{\textrm{fib}})$, where $(-)^{\textrm{fib}}$ denotes fibrant 
replacement in simplicial categories, yields an associated
$\i$-category.  Barwick and 
Kan \cite{BarwickKan} have studied this procedure in the context of
Segal spaces and Lurie has given a version of this approach
in \cite[\S 1.3.3]{HAG}: we associate to a (not necessarily
simplicial) category $\aC$ with weak equivalences $W$ an $\i$-category
$\N(\aC)[W^{-1}]$; when $\aC$ is a model category, for functoriality
reasons it is usually convenient to restrict to the cofibrant objects
$\aC^{\mathrm{c}}$ and consider $\N(\aC^{\mathrm{c}})[W^{-1}]$.

All of these constructions produce equivalent
$\i$-categories \cite[1.3.7]{HAG}.  Furthermore, all of them are
functorial.  Although a simplicial left Quillen functor $\aC \to \aD$
does not typically induce a functor between the subcategories of
cofibrant-fibrant objects in $\aC$ and $\aD$ respectively, composing
with a fibrant replacement functor as in definition~\ref{def:Morita},
does yield an induced functor on simplicial nerves.  Furthermore,
given a simplicial Quillen adjunction $(F,G)$, there is an induced
adjunction of functors on the level of $\i$-categories
by \cite[5.2.4.6]{HTT}.  (Alternatively, it can be seen directly that
a functor $\aC \to \aD$ which preserves weak equivalences between
cofibrant objects induces a functor
$\N(\aC)[W^{-1}] \to \N(\aD)[W^{-1}]$.)

Given an $\i$-category $\aC$, we can form its homotopy category
  $\Ho(\aC)$, which is an ordinary category \cite[\S 1.2.3]{HTT}.  In
  addition, given an $\i$-category $\aC$, there is a maximal
  $\i$-groupoid (Kan complex) $\aC_{\mathrm{iso}}$ inside of $\aC$,
  obtained by restricting to the subcategory of $\aC$ consisting of
  those arrows which become isomorphisms in the homotopy category
  $\Ho(\aC)$.  The functor which associates to the $\i$-category $\aC$
  its maximal subgroupoid $\aC_{\mathrm{iso}}$ is right adjoint to the
  inclusion of $\i$-groupoids into $\i$-categories.  We have the
  following proposition relating this to other, possibly more
  familiar, notions (see also \cite[2.3]{ToenVezzosi}).

\begin{proposition}\label{prop:nerve}
Let $\aC$ be a small category with a subcategory $w\aC$ of weak
equivalences which satisfies a homotopy calculus of two-sided
fractions (in the sense of Dwyer and Kan \cite[6.1]{DwyerKan2}).  Then
there is a weak equivalence of simplicial sets
\[
\N(w \aC)\simeq\lkan{\N((L^H \aC)^{\mathrm{fib}})},
\]
where here $L^H \aC$ denotes the hammock version of the simplicial
localization \cite{DwyerKan2} and $(L^H \aC)^{\mathrm{fib}}$ is a fibrant
replacement of $L^H\aC$ as a simplicial category.
\end{proposition}

\begin{proof}
There is an ``inclusion'' functor $\aC \to L^{H}\aC$.  Restricting to
the weak equivalences and passing to nerves via $\N$, we
obtain a map of simplicial sets
\begin{equation}\label{eqn:comp}
\N(w\aC) \to \N(L^{H} w\aC) \to \N((L^H w\aC)^{\mathrm{fib}});
\end{equation}
note that the nerve of $w\aC$ is the same whether we regard it as a category or as a category (trivially) enriched in simplicial sets \cite[1.1.5.8]{HTT}.  
Since $\lkan{\N((L^{H} w\aC)^{\mathrm{fib}})}$ is isomorphic to $\N((L^{H}
w\aC)^{\mathrm{fib}})$, the inclusion $L^H w\aC \to L^H \aC$
induces a natural map $\N((L^{H} w\aC)^{\mathrm{fib}}) \to \lkan{\N((L^H \aC)^{\mathrm{fib}})}$;
under the hypothesis that $\aC$ satisfies a homotopy calculus of
fractions, this map is a weak equivalence \cite[6.4]{DwyerKan2}.
Therefore, it suffices to show that the map of equation~\ref{eqn:comp}
is a weak equivalence.  We consider the map on components; for each
homotopy equivalence class $[x]$, both sides are equivalent to
$B\haut(x)$ and it is straightforward to see that the map induces the
equivalence.
\end{proof}

The $\i$-category of functors between two $\i$-categories $\aC$ and
$\aD$ is denoted $\Fun(\aC, \aD)$.  As a point-set object, in this
setting $\Fun(\aC, \aD)$ is the simplicial set of maps between the
quasicategories $\aC$ and $\aD$, which is itself a quasicategory.
Note that the {\em space} of functors from $\aC$ to $\aD$ is precisely
the maximal subgroupoid $\Fun(\aC,\aD)_{\mathrm{iso}}$ of
$\Fun(\aC,\aD)$ \cite[1.2.5.3, 3.0.0.1]{HTT}.

\begin{definition}
An $\i$-category is {\em stable} \cite[1.1.1.9]{HAG} if it has finite limits
and colimits and pushout and pullback squares
coincide \cite[1.1.3.4]{HAG}.  Let $\stabcat$ denote the (pointed)
$\infty$-category of small stable $\i$-categories and exact functors (i.e.,
functors which preserve finite limits and colimits) \cite[\S
1.1.4]{HAG}.  The $\i$-category of exact functors between $\aA$ and
$\aB$ is denoted by $\Fun^{\ex}(\aA, \aB)$; this is the full
$\i$-subcategory of $\Fun(\aA,\aB)$ spanned by the exact functors.
\end{definition}

For a small stable $\i$-category $\aC$, the homotopy
category $\Ho(\aC)$ is triangulated, with the exact triangles
determined by the cofiber sequences in $\aC$ \cite[1.1.2.13]{HAG}.

\begin{remark}\label{rem:pretriangstab}
A small stable $\i$-category corresponds to the notion of a
pretriangulated spectral category, and the weak equivalences are
given by exact functors which induce triangulated equivalences on
passage to the homotopy category.  We will make this correspondence
precise in section~\ref{sec:specinf}, but for now observe that given a
pretriangulated spectral category $\aC$, the $\i$-category
$\N((\modules{\aC})^{\cf})$ is stable.  Recall that a stable model
category is a pointed model category $\aC$ for which the functors
$\Sigma$ and $\Omega$ on $\Ho(\aC)$ are inverse equivalences.  Given a
stable simplicial model category $\aC$, the $\i$-category
$\N(\aC^{\cf})$ is stable.  More generally, if $\aC$ is a stable
model category, $\N(\aC^{\mathrm{c}})[W^{-1}]$ is a stable $\i$-category.
\end{remark}

Recall that an $\i$-category $\aC$ is idempotent-complete if the image
of $\aC$ under the Yoneda embedding $\aC \to \PPre(\aC)$ is closed
under retracts (see also \cite[\S 4.4.5]{HTT}); here $\PPre(\aC)$
denotes the $\i$-category $\Fun(\aC^{\op}, \N(\aT^{\cf}))$ of
presheaves of spaces on $\aC$.  Let $\idemstabcat$ denote the
$\i$-category of small idempotent-complete stable $\i$-categories.
There is an idempotent completion functor given as the left adjoint to
the inclusion $\idemstabcat \to \stabcat$ \cite[5.1.4.2]{HTT}, which
we denote by $\Idem$.

\begin{definition}\label{defn:iMor}
Let $\aA$ and $\aB$ be small stable $\i$-categories.  Then we will say
that $\aA$ and $\aB$ are {\em Morita equivalent} if $\Idem(\aA)$ and
$\Idem(\aB)$ are equivalent.
\end{definition}

We will verify shortly that this notion of Morita equivalence is
compatible with the definition given in terms of spectral categories
in definition~\ref{def:Morita}.

\subsection{Stabilization of $\i$-categories}\label{sec:stabilization}

Given any $\i$-category $\aC$ with finite limits, we can form the
stabilization $\Stab(\aC)$ \cite[\S 1.4]{HAG}. 
The $\i$-category $\Stab(\aC)$ is stable
and comes equipped with a limit-preserving functor  
\[
\Omega^{\infty} \colon \Stab(\aC) \to \aC.
\]
If in addition $\aC$ is presentable, then $\Omega^{\infty}$ admits a left adjoint
\[
\Sigma^\infty_+ \colon \aC \to \Stab(\aC)
\]
by \cite[1.4.4.4]{HAG}.

We now recall an explicit model of the stabilization of an
$\i$-category in terms of {\em spectrum objects} \cite[\S
1.4.2]{HAG}.  Recall that a spectrum object of a pointed $\i$-category
$\C$ consists of a functor $\N(\bZ \times \bZ) \to \C$.  In
particular, there are families of objects $A(i,j)$ of $\C$ and maps
$A(i,j)\to A(i+1,j)$, $A(i,j)\to A(i,j+1)$ such that $A(i,j)$ is zero
object whenever $i\neq j$ and the square 
\[
\xymatrix{
A(i,i)\ar[r]\ar[d] & A(i,i+1)\ar[d]\\
A(i+1,i)\ar[r] & A(i+1,i+1)}
\]
is cartesian for all $i$; consult \cite[1.4.2.4]{HAG} for further
details. Since the restriction of $A$ to the diagonal carries the
nontrivial objects in $A$, we set $A_i=A(i,i)$ and often refer to $A$
simply by the collection of pointed objects $\{A_i\}$.  We write
$\Spec(\C)$ for the $\i$-category of spectrum objects in $\C$;
$\Spec(\C)$ comes equipped with a functor
$\Omega^\i \colon \Spec(\C)\to\C$ which associates to the spectrum
object $A$ its zero space $A_0=A(0,0)$.  This is an explicit model for
the stabilization $\Stab(\aC)$ discussed previously.  To ease
notation, we will usually just write $\aT_\i\simeq\N(\T^{\cf})$ for the $\i$-category of spaces and $\ispec\simeq\Spec(\aT_\i)\simeq\N(\cS^{\cf})$ for the $\i$-category of spectra.

Now suppose that $\C$ is an arbitrary $\infty$-category. The Yoneda
embedding $\C\to\Fun(\C^{\op},\aT_\i)$ preserves finite limits (when they
exist), so it induces a functor 
\[
\Spec(\C_*)\to\Spec(\Fun(\C^{\op},\aT_\i)_*)\simeq\Fun(\C^{\op},\ispec)
\]
on the level of spectrum objects, where $\C_*$ denotes the category of
pointed objects in $\C$.  Here the last equivalence follows
from the fact that limits in functor categories are computed
pointwise, and observe also that $\C_*$ will be empty unless $\C$ has
a final object. 
On the other hand, if $\C$ is a stable $\i$-category, then
$\C\simeq\C_*$ and $\Omega^\i \colon \Spec(\C)\to\C$ is an equivalence
with inverse $\Sigma^\i \colon \C\to\Spec(\C)$ given by $(\Sigma^\i
a)_i=\Sigma^i a$ \cite[1.4.2.20]{HAG}. 
This motivates the following definition:

\begin{definition}\label{defn:mappingspectrum}
Let $\aC$ be stable $\i$-category.
The {\em spectral Yoneda embedding} is the composite
\[
\C\simeq\Spec(\C_*)\to\Spec(\Fun(\C^{\op},\N(\aT)^{\cf})_*)\simeq\Fun(\C^{\op},\ispec).
\]
The {\em mapping spectrum} functor
\[
\Map\colon \C^{\op}\times\C\to \ispec
\]
is the adjoint of the spectral Yoneda embedding.
\end{definition}

Informally, the mapping spectrum is described by the formula
\[
\Map(b,a)_i\simeq\map(b,\Sigma^i a).
\]
Note that this is a functor to the $\i$-category of spectra;
this is in contrast to the (point-set) mapping space functors from the
category of quasicategories to the category of simplicial sets described
in \cite[1.2.2]{HTT} or \cite{DuggerSpivak}.

We wish to characterize the image of $\aC$ under the spectral Yoneda
embedding:

\begin{definition}
Let $\C$ be an $\i$-category.  Then we will say that a functor
$X\colon\C^{\op}\to\ispec$ is {\em stably representable} if there exists a
spectrum object $A\in\Spec(\C_*)$ and an equivalence $\Map(-,A)\simeq
X$, where $\Map(-,A)$ denotes the functor $\C^{\op}\to\ispec$
represented by $A$ via the spectral Yoneda embedding
$\Spec(\C_*)\to\Fun(\C^{\op},\ispec)$. 
\end{definition}

When $\C$ is stable already, the following proposition gives an easy
characterization of stably representable functors.

\begin{proposition}
Let $\C$ be a stable $\i$-category.  Then a functor
$X\colon \C^{\op}\to\ispec$ is stably representable if and only if it is
represented by the suspension spectrum $\Sigma^\i z$ of a unique (up
to equivalence) object $z$ of $\C$. 
\end{proposition}

\begin{proof}
It suffices to show that any spectrum object $A$ of $\C$ is of the
form $\Sigma^\i z$ for a uniquely determined object $z$ of $\C$.
This follows from the fact that since $\aC$ is stable,
$\Omega^\i \colon \Spec(\C)\to\C$ is an equivalence with inverse
$\Sigma^\infty \colon \C\to\Spec(\C)$.
\end{proof}

\subsection{Compact objects and compactly-generated
$\i$-categories}\label{sec:compact} 

The categorical data which serves as the input to algebraic $K$-theory
is typically obtained as the objects in a larger ambient category
(with weak equivalences and extension sequences) that satisfy some
sort of ``smallness'' condition; e.g., the perfect complexes as a
subcategory of all complexes.  A key insight initially codified by
Thomason-Trobaugh \cite{TT} and subsequently elaborated upon by Neeman
\cite{Neeman} is that this example is generic in algebraic $K$-theory, and the typical
situation involves working with the {\em compact objects} in some
model of a triangulated category, which is generated under homotopy
colimits by those compact objects.  Thus, we will systematically
regard the small stable idempotent-complete $\i$-categories that are
the domain of the algebraic $K$-theory functor as arising as the
compact objects in a larger category.

This notion of looking at large categories which are in some sense
determined by the compact objects is axiomatized in category theory
with the formalism of {\em accessible} and {\em locally presentable}
categories, introduced by Makkai and Par\'{e}~\cite{MakkaiPare} and
further developed by Ad\'{a}mek and Rosick\'{y}~\cite{AdamekRosicky}.
This theory was integrated into homotopy theory in Jeff Smith's theory
of combinatorial model categories and developed further in this
context in the seminal work of Dugger~\cite{Dugger}.

A version of this theory forms the basis for Lurie's theory
of {\em presentable $\i$-categories}, which is the analogue in the
$\i$-category setting of the homotopy theories encoded by presentable
combinatorial model category structures (see also Simpson's related
work in the context of Segal spaces \cite{Simpson}).  We use this
approach to handle the set-theoretic issues that arise in our work,
along the lines described in \cite[1.2.15, 5.4.1]{HTT}.  As indicated
in remark~\ref{rem:spectralcard}, it is also possible to handle some
of the set-theoretic technicalities that arise (i.e., in the context
of the Yoneda lemma) by explicit size bounds.  

This framework is related to Grothendieck's universe formalization,
allowing us to handle small and large $\i$-categories on similar
grounds.  In particular, \cite[\S 5]{HTT} has extensive discussion of
the interaction of the Yoneda embeddings (which arise pervasively in
this context) with set-theoretic concerns.  In addition to Lurie's
work, the paper of Ben-Zvi, Francis, and Nadler \cite{BFN} provides a
nice expos\'{e} of this theory in the context of the study of
geometric function theory from a perspective with its origin in
Thomason-Trobaugh, and we refer the interested reader to sections 2 and
4.1 of that paper.

Roughly speaking, presentable $\i$-categories are large
$\i$-categories that are generated under sufficiently large filtered
colimits by some small $\i$-category.  To make this precise, we need
to discuss the notion of the $\Ind$-category.
Given any small $\i$-category $\aC$, we can form the $\i$-category
$\PPre(\aC)$ of presheaves of simplicial sets on $\aC$, which is the
formal closure of $\aC$ under colimits; that is, there is a
fully faithful Yoneda embedding $\aC
\to \PPre(\aC)$, and $\PPre(\aC)$ is generated by the image of $\aC$
under small colimits \cite[5.1.5.8]{HTT}.  For any $\i$-category $\aC$
and infinite regular cardinal $\kappa$, we can form the $\Ind$-category
$\Ind_\kappa(\aC)$, which is the formal closure under $\kappa$-filtered
colimits of $\aC$ \cite[\S 5.3.5]{HTT}.  The $\i$-category
$\Ind_\kappa(\aC)$ is a full subcategory of $\PPre(\aC)$, and the
Yoneda embedding $\aC \to \PPre(\aC)$ factors as
$\aC \to \Ind_{\kappa}(\aC)
\to \PPre(\aC)$.  We record here the following useful properties of
the construction of the $\Ind$-category.

\begin{proposition}\label{prop:Ind}
Let $\aC$ be a small $\i$-category and $\kappa$ an infinite regular
cardinal.
\begin{itemize}
\item The $\i$-category $\Ind_{\kappa}(\aC)$ admits all $\kappa$-small
colimits that exist in $\aC$ \cite[5.3.5.14, 5.5.1.1]{HTT}. 
\item The functor $\aC \to \Ind_{\kappa}(\aC)$ preserves
$\kappa$-filtered colimits \cite[5.3.5.2,
5.3.5.3]{HTT}. 
\item $\Ind_{\kappa}(\aC)$ is a stable $\i$-category \cite[1.1.3.6]{HAG}.
\item The image of $\aC$ in $\Ind_{\kappa}(\aC)$ provides a set of
compact objects which generates $\Ind(\aC)$ under $\kappa$-filtered
colimits \cite[5.3.5.5,5.3.5.11]{HTT}.
\item The category $\Ind_{\kappa}$ is characterized by the property
that it has $\kappa$-small filtered colimits, admits a functor
$\aC \to \Ind_{\kappa}(\aC)$, and this functor induces an equivalence
\[
\Fun_{\kappa}(\Ind(\aC), \aD) \to \Fun(\aC,\aD), 
\]
for any $\aD$ which admits $\kappa$-filtered colimits (here
$\Fun_{\kappa}(-,-)$ denotes the $\i$-category of functors that preserve
$\kappa$-small filtered colimits) \cite[5.3.5.10]{HTT}.
\end{itemize}
\end{proposition}

We now recall the following definitions \cite[5.4.2.1,5.5.1.1]{HTT}.

\begin{definition}
An $\i$-category $\aC$ is accessible if there exists a regular
cardinal $\kappa$ and a small $\i$-category $\aC^{0}$ such that there
is an equivalence 
\[
\Ind_{\kappa}(\aC^{0}) \htp \aC.
\] 
An $\i$-category $\aC$ is presentable if it arises as
$\Ind_{\kappa}(\aD)$ for a small $\i$-category $\aD$ which admits
$\kappa$-small colimits \cite[5.5.1.1]{HTT}.  
\end{definition}

A morphism of presentable $\i$-categories is a left adjoint functor; by the adjoint functor theorem \cite[5.5.2.9]{HTT}, a functor between presentable $\i$-categories is a left adjoint if and only if it preserves colimits.  We let $\prcat$ denote the
$\i$-category of presentable $\i$-categories and colimit-preserving
functors; the $\i$-category of colimit-preserving functors is denoted
by $\Fun^{\L}(-,-)$.  In fact, $\Fun^{\L}(-,-)$ is in fact itself a
presentable $\i$-category \cite[5.5.3.8]{HTT}, yielding an
internal $\hom$ object for $\prcat$.

We now restrict attention to the situation in which $\kappa = \omega$.
Recall that an object $x$ of an $\i$-category $\aC$ is compact if the
functor $\aC^{\op}\to\aT$ represented by $x$ commutes with filtered
colimits \cite[\S 5.3.4]{HTT}.  Given an an
$\i$-category $\aC$, let $\aC^{\omega}$ denote the full subcategory of
$\aC$ consisting of the compact objects
of $\aC$.  A presentable $\i$-category $\aC$ is {\em compactly
  generated} if the natural functor
\[
\Ind(\aC^{\omega})\to\aC,
\]
which sends a filtered diagram in $\aC^\omega$ to its colimit in $\aC$,
is an equivalence.  There is a correspondence between small
idempotent-complete $\i$-categories and 
compactly generated $\i$-categories given by the construction of the
$\Ind$-category \cite[\S 5.5.7]{HTT}.
More generally, the construction of the
$\Ind$-category sets up a correspondence between the $\i$-category of
compactly-generated presentable $\i$-categories with morphisms
colimit-preserving functors that preserve compact objects and
$\icat$; the other direction is given by passage to
compact objects \cite[5.5.7.10]{HTT}.

The preceding discussion carries over when we restrict attention to
stable categories.  In this setting, the stabilization $\Stab(\aC)$ is
initial amongst presentable stable $\i$-categories admitting a
functor from $\aC$ \cite[1.4.5.5]{HAG}, in the sense that if $\aD$ is
a presentable stable $\i$-category then $\Sigma^\infty_+$ induces an
equivalence   
\[
\Fun^{\L}(\Stab(\aC), \aD) \to \Fun^{\L}(\aC, \aD). 
\]

The $\i$-category of stable presentable
$\i$-categories $\stabprcat$ is a full subcategory of $\prcat$, and
the $\Ind$-category sets up a correspondence between $\idemstabcat$
and compactly generated stable $\i$-categories.  We may also apply
$\Ind$ to non-idempotent-complete stable $\i$-categories to obtain a
correspondence between $\stabcat$ and compactly generated stable
$\i$-categories; however, these two $\i$-categories are rather less
closely related, as the full subcategory of compact objects is always
idempotent-complete.

\begin{lemma}
$\idemstabcat$ is a reflective subcategory of $\stabcat$, and the
localization functor $\Idem\colon\stabcat\to\idemstabcat$ is given by
the formula $\Idem(\aC)\simeq\Ind(\aC)^\omega$.
\end{lemma}

\begin{proof}
The subcategory of compact objects $\Ind(\aC)^\omega$ of
$\Ind(\aC)$ is an idempotent-complete stable $\i$-category, so that
$\Idem$ is indeed a functor $\stabcat\to\idemstabcat$.  Now for small
stable $\i$-categories $\aC$ and $\aD$ with $\aD$ idempotent-complete,
we have a commuting square
\[
\xymatrix{
\Fun^{\ex}(\Idem(\aC),\aD)\ar[r]\ar[d]
& \Fun^{\L}(\Ind(\Idem(\aC)),\Ind(\aD))\ar[d] \\
\Fun^{\ex}(\aC,\aD)\ar[r] & \Fun^{\L}(\Ind(\aC),\Ind(\aD))}
\]
in which the horizontal maps are the inclusions of the full
subcategories of functors which preserve compact objects, and the
right vertical map is an equivalence as the natural map
$\Ind(\aC)\to\Ind(\Idem(\aC))$ is an equivalence.  Hence
$\Ind(\aC)^\omega\to\Ind(\Idem(\aC))^\omega$ is an equivalence, and
thus the left vertical map is as well.
\end{proof}

\subsection{Localization of $\i$-categories}\label{sec:bousloc}

Given an $\i$-category $\aC$ and a suitable collection of morphisms 
$S$, one might hope to form the localization $\aC [S^{-1}]$.
This is by definition an $\i$-category equipped with a functor
$f\colon \aC\to\aC[S^{-1}]$ which satisfies the following universal
property: for any other $\i$-category $\aD$, restriction along $f$
identities 
\[
\Fun(\aC[S^{-1}],\aD)\to\Fun(\aC,\aD)
\]
as the full subcategory of $\Fun(\aC,\aD)$ spanned by those functors
which send the morphisms in $S$ to equivalences in $\aD$.  If $S$ is a
proper set, then $\aC[S^{-1}]$ exists in the same universe as $\aC$;
indeed, without loss of generality we may assume that $S$ contains all
degenerate edges of the simplicial set $\aC$, in which case
$\aC[S^{-1}]$ may be constructed as a fibrant replacement of $(\aC,S)$
in the model category of marked simplicial sets.

Often in practice, however, $S$ is not small, and the existence of
$\aC[S^{-1}]$ (without passing to a higher universe) requires more
delicate analysis. One standard method is to show that $\aC$ is
presentable and $S$ is (generated by) a small set of arrows in a certain sense: this is
the theory of Bousfield localization, following Bousfield's seminal
work on the subject~\cite{Bousfield}.  In this case we may identify
the localization as the full subcategory of $\aC$ spanned by the
$S$-local objects.

In model categories, there is a well-developed theory of
Bousfield localization (e.g., Hirschhorn's comprehensive discussion in
~\cite{Hirschhorn}, Goerss and Jardine's treatment in the simplicial
setting~\cite{GoerssJardine}, or the exposition of Smith's theory for
combinatorial model categories in ~\cite{BarwickLoc}).  Because
localization is a central technical device in our work, in this
section we provide a brief review of Lurie's version of localization
in the setting of presentable $\i$-categories from \cite[\S 5.2.7]{HTT}
and \cite[\S 5.5.4]{HTT}.

Specifically, we say that a colimit preserving functor $f \colon \aC \to \aD$
of presentable $\infty$-categories $\aC$ and $\aD$ is a {\em Bousfield
localization} if the right adjoint of $f$ (which exists by the adjoint
functor theorem) is fully faithful \cite[5.2.7.2]{HTT}.  When the
context is clear, we tend to abuse notation and simply refer to this
as a localization. A useful observation is that this data induces an
equivalence between $\aD$ and a full subcategory of $\aC$, called the
subcategory of {\em local objects}.  In fact, \cite[5.2.7.4]{HTT}
gives a useful criterion for an endofunctor $L \colon \aC \to \aC$ to
be a localization.  Specifically, the following are equivalent: 

\begin{enumerate}
\item There exists a functor $f \colon \aC \to \aD$ with a fully
faithful right adjoint $g$ and an equivalence $g \circ f \htp L$.
\item When regarded as a functor $\aC \to L\aC$, $L$ is the left
adjoint of the inclusion $L\aC \to \aC$.
\item There exists a natural transformation
$\alpha \colon \id_{\aC} \to L$ such that for objects $X$ in $\aC$,
the morphisms $L(\alpha(X))$ and $\alpha(LX)$ are both equivalences.
\end{enumerate}

Recall that a functor is accessible if it is $\kappa$-continuous
(preserves $\kappa$-filtered colimits) for some sufficiently large
regular cardinal $\kappa$ \cite[5.4.2.5]{HTT}.  A localization is
accessible if $g$ or $L$ are accessible functors (equivalently,
see \cite[5.5.1.2]{HTT}) or the the essential image $L\aC$ is an
accessible subcategory.

Accessible localizations of presentable categories can be completely
classified as follows.  Recall from \cite[5.5.4]{HTT} that associated
to any {\em set} of arrows $S$ in a presentable $\i$-category $\aC$,
the Bousfield localization $S^{-1}\aC$ is equivalent to the ordinary
localization $\aC[T^{-1}]$ of $\aC$ at the {\em strongly saturated
class} $T$ generated by $S$ \cite[5.5.4.5]{HTT}.  In particular, many
different sets $S$ can generate the same strongly saturated class $T$;
they all define the same full subcategory $S^{-1}\aC$ of $\aC$ of
$S$-local objects \cite[5.5.4.15]{HTT}, where $S$-local is defined in
the standard fashion \cite[5.5.4.1]{HTT}.  An accessible localization
of a presentable $\i$-category is presentable.  As the notation
suggests, Bousfield localization is characterized by the following
universal property \cite[5.5.4.20]{HTT}: for any other presentable
$\i$-category $\aD$, composition with $L$ induces a functor
\[
\Fun^{\L}(S^{-1}\aC, \aD) \to \Fun^{\L}(\aC, \aD)
\]
which is fully faithful and whose essential image consists of those
colimit-preserving functors which take elements of $S$ to equivalences.

\section{Symmetric monoidal structure and dualizable objects}\label{sec:monoidal}

In this section we study the theory of dualizable objects in
$\idemstabcat$.  To this end, we need to give a very brief review of
the construction of the symmetric monoidal structure on
$\idemstabcat$.  We do not give a full review of the theory of
monoidal $\i$-categories in this section, since we need only a small
piece of the theory.

\subsection{Tensor products of stable $\i$-categories}
The $\i$-category $\stabprcat$ of presentable stable $\i$-categories is a
closed symmetric monoidal $\i$-category with product $\otimes$ and internal
mapping object given by the presentable stable $\i$-category
$\Fun^{\L}(\aA,\aB)$ of colimit-preserving functors \cite[6.3.1.14,
6.3.1.17]{HAG}.  Following \cite[\S 4.1.2]{BFN}, we can then define
the tensor product on small idempotent-complete stable $\i$-categories
as
\[
\aC \idemtimes \aD = (\Ind(\aC) \otimes \Ind(\aD))^{\omega}.
\]
The tensor product of idempotent-complete small stable $\i$-categories
is characterized by the universal property that maps out of $\aA
\otimes \aB$ correspond to maps out of the product $\aA \times \aB$
which preserve finite colimits in each variable \cite[4.4]{BFN}.
If $\aA$ and $\aB$ are arbitrary small stable $\i$-categories, then we set $\aA\idemtimes\aB:=\Idem(\aA)\idemtimes\Idem(\aB)$.

More precisely, we can define $\idemstabcat$ as a symmetric monoidal
$\i$-category as follows.  Let $\stabprcatcg$ denote the full
subcategory of $\stabprcat$ on the compactly-generated stable
$\infty$-categories.  The criterion of \cite[2.2.1.2]{HAG} implies
that $\stabprcatcg$ is a symmetric monoidal subcategory of
$\stabprcat$; the tensor product of compactly-generated stable
$\infty$-categories is itself compactly-generated, as is the unit
$\ispec\simeq\Ind(\ispec^\omega)$.

For a small stable idempotent-complete $\i$-category $\aA$ and a
presentable $\i$-category $\aB$, $\Fun^{\ex}$ and $\Fun^{\L}$ are
related by the formula
\[
\Fun^{\ex}(\aA, \aB) \htp \Fun^{\L}(\Ind(\aA), \aB),
\]
which follows from \cite[5.3.5.10]{HTT} and the fact that functors
which preserve filtered colimits and finite colimits preserve all
colimits.
Note that 
\[
\Ind \colon \idemstabcat\to\stabprcat
\]
factors through the full subcategory $\stabprcatcg$ by
definition.  This gives an equivalence of $\infty$-categories
between $\idemstabcat$ and the subcategory $\stabprcat_\omega$ of
$\stabprcat$ whose objects are the compactly-generated stable
$\infty$-categories and whose maps
\[
\Fun^{\ex}(\aA,\aB)\simeq\Fun^{\L}_\omega(\Ind(\aA),\Ind(\aB))\subset\Fun^{\L}(\Ind(\aA),\Ind(\aB)), 
\]
are the full subcategory of the colimit-preserving functors
$\Ind(\aA)\to\Ind(\aB)$ which preserve compact
objects \cite[5.5.7.10]{HTT}.  We regard $\idemstabcat$ as a symmetric 
monoidal $\i$-category via this equivalence.  The observation
of \cite[6.3.1.17]{HAG} implies that $\idemstabcat$ is closed.  Hence
we have the following result.

\begin{theorem}\label{thm:idemsymm}
The $\i$-category of small idempotent-complete stable $\i$-categories is a
closed symmetric monoidal category with respect to $\idemtimes$.
The unit is the $\i$-category $\ispec^{\omega}$ of compact spectra and the
internal mapping object is given for small idempotent-complete stable
$\i$-categories $\aA$ and $\aB$ by $\Fun^{\ex}(\aA, \aB)$.
\end{theorem}

Given a small stable idempotent-complete $\i$-category $\aA$, we have
the $\i$-category of $\aA$-modules, given by the compactly-generated stable
$\i$-category $\Fun^{\ex}(\aA^{\op},\ispec)$.  
The stable Yoneda embedding provides an exact functor \cite[5.3.5.2]{HTT}
\[
\aA \to \Fun^{\ex}(\aA^{\op},\ispec).
\]

\begin{proposition}\label{prop:idempotentchar}
For any small stable $\i$-category $\aA$, the stable Yoneda embedding
\[
\aA\longrightarrow\Fun^{\ex}(\aA^{\op},\ispec)
\]
induces an equivalence $\Ind(\aA)\simeq\Fun^{\ex}(\aA^{\op},\ispec)$.
\end{proposition}

\begin{proof}
Clearly $\Fun^{\ex}(\aA^{\op},\ispec)$ admits filtered colimits, as the
filtered colimit of finite colimit preserving functors itself
preserves finite colimits.  This gives a map
$\Ind(\aA)\to\Fun^{\ex}(\aA^{\op},\ispec)$ which is evidently fully
faithful since, using the fact that the usual Yoneda embedding is
fully faithful and that mapping spaces between representables in
$\Fun^{\ex}(\aA^{\op},\ispec)$ are computed as the limit 
\[
\lim_n\Omega^n\map(a,\Sigma^n b)\simeq\lim_n\map(a,\Omega^n\Sigma^n b)\simeq\map(a,b).
\]
To show that this map is also essentially surjective, we must show
that any exact functor $f\colon\aA^{\op}\to\ispec$ is ind-representable. 
Consider the pullback
\[
\xymatrix{\aA_{/f}\ar[r]\ar[d] & \aA\ar[d]\\
\Fun^{\ex}(\aA^{\op},\ispec)_{/f}\ar[r] & \Fun^{\ex}(\aA^{\op},\ispec),}
\]
where the right vertical map is the stable Yoneda embedding.
We claim that the $\i$-category $\aA_{/f}$ is filtered: to see this, let $K$ be
a finite simplicial set and $K\to\aA_{/f}$ a functor.  Since both
$\aA$ and $\Fun^{\ex}(\aA^{\op},\ispec)_{/f}$ admit finite colimits and
both functors to $\Fun^{\ex}(\aA^{\op},\ispec)$ preserve finite colimits,
we may extend $K\to\aA_{/f}$ to a colimit diagram
$K^\triangleright\to\aA_{/f}$.  In particular, this gives a cone on
$K\to\aA_{/f}$, which shows that $\aA_{/f}$ is a filtered
$\i$-category.  Finally, since filtered colimits in
$\Fun^{\ex}(\aA^{\op},\ispec)$ are computed pointwise, it follows that
$f$ is a colimit of the diagram
$\aA_{/f}\longrightarrow\Fun^{\ex}(\aA^{\op},\ispec)$, which is to say
that it is ind-representable.
\end{proof}

Provided $\aA$ is idempotent-complete, ~\cite[5.4.2.4]{HTT} tells us that the essential image of
the Yoneda embedding is precisely the $\infty$-category of compact $\aA$-modules  
\[
\aA \htp \Fun^{\ex}(\aA^{\op},\ispec)^{\omega}.
\]
Moreover, we know that if $\aA$ is an arbitrary small stable $\i$-category, then the Yoneda map $\aA \to \Fun^{\ex}(\aA^{\op},\ispec)^{\omega}$ models the idempotent-completion 
of $\aA$.

We use the preceding results to characterize $\Fun^{\ex}(\aA,\aB)$ in
terms of a certain subcategory of $\Fun^{\L}(\aA \idemtimes
\aB^{\op}, \ispec)$, the $\i$-category of $\aA$-$\aB$-bimodules.  Specifically,
the Yoneda embedding $\aB \to \Fun^{\ex}(\aB^{\op},\ispec)$ provides the
following composite
\begin{align*}
\Fun^{\ex}(\aA, \aB)
&\to \Fun^{\ex}(\aA, \Fun^{\ex}(\aB^{\op},\ispec)) \\
&\to \Fun^{\L}(\Ind(\aA), \Fun^{\L}(\Ind(\aB^{\op}), \ispec)) \\
&\to \Fun^{\L}(\Ind(\aA) \otimes \Ind(\aB^{\op}), \ispec) \\
&\to \Fun^{\ex}(\aA \idemtimes \aB^{\op}, \ispec),
\end{align*}
which exhibits $\Fun^{\ex}(\aA, \aB)$ as a full subcategory of
$\aA^{\op} \idemtimes \aB$-modules.

We have the following useful corollary, which is the analogue of a
characterization originally written down by To\"en \cite{Toen}.  For
each object $a \in \aA$, we have a map of small idempotent-complete
stable $\infty$-categories $\ispec^\omega\to\aA$ given by sending
$\bS\in\ispec^\omega$ to $a\in\aA$.  Since $\ispec^\omega$ is the unit of
the tensor $\idemtimes$, we obtain a restriction map
\[
\nu_{a} \colon \Fun^{\ex}(\aA \idemtimes \aB^{\op},\ispec) \to
\Fun^{\ex}(\aB^{\op}, \ispec).
\]
If the image of an element of $\Fun^{\ex}(\aA \idemtimes
\aB^{\op},\ispec)$ under $\nu_a$ is compact for every $a \in \aA$, we will
say that the element is {\em right-compact}.

\begin{corollary}\label{cor:toencor}
Let $\aA$ and $\aB$ be small stable idempotent-complete
$\i$-categories.  There is an equivalence of small stable
idempotent-complete $\i$-categories between $\idemfun(\aA, \aB)$
and $\i$-category of right-compact $\aA^{\op} \idemtimes
\aB$-modules.
\end{corollary}

\begin{proof}
Since the Yoneda embedding is fully faithful, it suffices to look at
the essential image of the composite.  The image of
$\Fun^{\ex}(\aA, \aB)$ in $\Fun^{\ex}(\aA \idemtimes
\aB^{\op}, \ispec)$ under $\nu_a$ is identified with the image of
$\Fun^{\ex}(\aA, \aB)$ in $\Fun^{\ex}(\aA,
\Fun^{\ex}(\aB^{\op},\ispec))$ under the corresponding map $\ispec^\omega\to\aA$.
Since this lies inside the image of $B$ inside
$\Fun^{\ex}(\aB^{\op},\ispec)$ under the Yoneda embedding, the result
follows.
\end{proof}

\subsection{Smooth and proper stable $\i$-categories}
Our final goal in this section is to characterize the dualizable
objects of $\idemstabcat$.  To do so, we need to introduce certain
smallness conditions on small stable $\infty$-categories. 

\begin{definition}\label{def:proper}
A small stable $\infty$-category $\aA$ is {\em proper} if, for all pairs of objects $a$ and $b$ of $\aA$, the mapping spectrum $\aA(a,b)$ (recall
definition~\ref{defn:mappingspectrum}) is compact. 
\end{definition}

Note that a small stable $\infty$-category $\aA$ is proper if and only if its idempotent-completion $\Idem(\aA)$ is proper, as retracts of compact objects are compact.

\begin{definition}\label{def:smooth}
A small stable $\infty$-category $\aA$ is {\em
smooth} if it is perfect as an $\aA^{\op}\idemtimes\aA$-module (i.e.,
in the smallest subcategory of $\idemfun(\aA\idemtimes\aA^{\op}, \ispec)$ generated by
the representables under finite colimits and retracts).  If
$\aA$ is idempotent-complete, we may equivalently require that $\aA$ is a representable
$\aA^{\op}\idemtimes\aA$-module: since $\aA^{\op}\idemtimes\aA$ is an
idempotent-complete stable $\infty$-category, it is closed under
finite colimits and retracts, and so any perfect
$\aA^{\op}\idemtimes\aA$-module is representable.
\end{definition}

We will typically only be interested in smoothness and properness of small stable $\i$-categories which are also idempotent-complete.
This is because these are the situations which arise in algebra and geometry, e.g. when $\aA$ is the stable $\i$-category of perfect modules for a ring spectrum or perfect complexes for a scheme, and such categories
are always idempotent complete. Conversely (as we will show in
section~\ref{sec:specinf}) any idempotent-complete small stable
$\infty$-category $\aA$ is equivalent to the stable $\infty$-category
of perfect modules for some spectral category.

\subsection{Dualizability}
We now recall the definitions of dualizability in symmetric monoidal
$\i$-categories from \cite[\S 4.2.5]{HAG}.  The salient fact here is
that dualizability can be detected in the (symmetric monoidal)
homotopy category:

\begin{definition}
Let $\aC^\otimes$ be a symmetric monoidal $\infty$-category.
An object of the underlying $\infty$-category $\aC$ of $\aC^\otimes$ is said to be {\em dualizable} if it is
dualizable as an object of the symmetric monoidal homotopy category of $\aC^\otimes$.
\end{definition}

In other words, an object $A$ of $\aC$ is dualizable if there 
exists an object $DA$ together with an evaluation map
$\epsilon:A\otimes DA\to 1$ and a coevaluation map $\delta:1\to
DA\otimes A$ such that the composites 
\[
A\simeq A\otimes 1\overset{A\otimes\delta}{\to} A\otimes DA\otimes
A\overset{\epsilon\otimes A}{\to} 1\otimes A\simeq A 
\]
and
\[
DA\simeq 1\otimes DA\overset{\delta\otimes DA}{\to} DA\otimes A\otimes
DA\overset{DA\otimes\epsilon}{\to} DA\otimes 1\simeq DA 
\]
are the respective identities in $\Ho(\aC)$.  The object $DA$ is called the {\em dual} of $A$, and is unique up to equivalence in $\aC$. 

Recall that the discussion preceding theorem~\ref{thm:idemsymm} above
identifies $\idemstabcat$ as a symmetric monoidal subcategory of
$\stabprcatcg$; in particular, the functor $\Ind$ preserves
dualizable objects.  Next, observe that $\stabprcatcg$
is a rigid symmetric monoidal category; that is, all objects in
$\stabprcatcg$ are dualizable.  This is because, for $\aA\in\idemstabcat$,
\[
\Ind(\aA^{\op})\simeq\Fun^{\ex}(\aA,\ispec)\simeq\Fun^\L(\Ind(\aA),\ispec)
\]
is the dual of $\aA$ and the coevaluation map
\[
\ispec\to\Ind(\aA^{\op})\otimes\Ind(\aA)\simeq\Ind(\aA^{\op}\idemtimes\aA)\simeq\Fun^{\ex}(\aA\idemtimes\aA^{\op},\ispec)
\]
is given by formation of mapping spectra in $\aA^{\op}$.

In analogy with the situation for dg-categories \cite[\S4]{CisTab},
this allow us to obtain the following characterization of the
dualizable objects.

\begin{theorem}\label{thm:dualizable}
An idempotent-complete small stable $\infty$-category $\aA$ is
dualizable (as an object of the symmetric monoidal $\infty$-category
$\idemstabcat$ of idempotent-complete small stable
$\infty$-categories) if and only if $\aA$ is smooth and proper.  Moreover, the dual of a dualizable object
$\aA$ is its opposite $\infty$-category $\aA^{\op}$.
\end{theorem}

\begin{proof}
By the proceeding discussion, $\Ind(\aA)$ is a dualizable object of $\stabprcatcg$ with dual $\Ind(\aA^{\op})$. Thus the dual of $\aA$ in $\idemstabcat$ is $\aA^{\op}$, and $\aA$ is dualizable in $\idemstabcat\simeq\stabprcat_\omega$ if and only if the evaluation and coevaluation maps lie in the subcategory $\stabprcat_\omega\subset\stabprcatcg$. But the evaluation map
\[
\Ind(\aA\idemtimes\aA^{\op})\simeq\Ind(\aA)\otimes\Ind(\aA)^*\longrightarrow\Ind(\ispec^\omega)\simeq\ispec
\]
is induced by the mapping spectrum functor $\Map_\aA \colon \aA^{\op}\idemtimes\aA\longrightarrow\ispec$ in $\aA$; dually, the coevaluation map 
\[
\ispec\simeq\Ind(\ispec^\omega)\longrightarrow\Ind(\aA^{\op})\otimes\Ind(\aA)\simeq\Ind(\aA^{\op}\idemtimes\aA)
\]
given by the map $\ispec\too\Ind(\aA^{\op}\idemtimes\aA)$ which classifies $\aA$ as an $\aA^{\op}\idemtimes\aA$-module.  Hence the evaluation map lies in this subcategory if and only if the mapping spectra $\aA(a,b)$
in $\A$ are compact, and the coevaulation map lies in this subcategory if and only if $\aA$ is a compact $\aA^{\op}\idemtimes\aA$-module.
Therefore, by definition, $\aA$ is a dualizable object of $\idemstabcat$ if and only if $\aA$ is smooth and proper.
\end{proof}

\section{Morita theory}\label{sec:specinf}

There is a close connection between stable $\i$-categories and
spectral categories.  On the one hand, for every pair of objects in a
stable $\i$-category we can extract a mapping spectrum, as we
discussed in definition~\ref{defn:mappingspectrum}.  On the other
hand, given a category $\aA$ enriched in spectra, the category of
(right) $\aA$-modules has a standard projective model structure and
the associated $\i$-category is stable.

The purpose of this section is to provide a precise account of the
relationship between small spectral categories and small stable
$\i$-categories.  The moral of the story is that the homotopy theory
of small spectral categories, localized at the Morita equivalences, is
the same as the homotopy theory of small stable idempotent-complete
$\i$-categories.  Specifically, we prove theorem~\ref{thm:models} from
the introduction, which can be thought of as a generalization of the Morita theory of \cite{SS}; that is, small stable idempotent-complete $\i$-categories are $\i$-categories of
modules, and the $\i$-category of exact functors between two such is a
stable subcategory of the $\i$-category of bimodules.

Establishing this correspondence serves several purposes for us.  For
one thing, having models of $\stabcat$ and $\idemstabcat$ as
accessible localizations of an $\i$-category which arises as the nerve
of a model category provides technical control on $\stabcat$ and
$\idemstabcat$; we use this to show that $\stabcat$ and $\idemstabcat$
are compactly generated in corollary~\ref{cor:cgen}.  For another, it
permits us to rectify diagrams of small stable $\i$-categories to
strict diagrams in $\Spcat$.  We exploit this to pass to rigid models
for the purposes of using Waldhausen's $K$-theory machinery in
Section~\ref{sec:kthy}.

As described in Section~\ref{sub:categories}, we have several
equivalent options for producing a model of the $\i$-category of
spectral categories (with respect to the DK-equivalences): we can use
the combinatorial simplicial model structure of
Corollary~\ref{cor:Quillencombi} and take $\N((\Spcat)^{\cf})$, we can
use the Dwyer-Kan simplicial localization followed by fibrant
replacement to obtain $\N((L^H \Spcat)^\mathrm{fib})$, or we can invert the weak
equivalences and form $\N((\Spcat)^{\mathrm{c}})[W^{-1}]$.  We will
refer interchangably to the underlying $\i$-category as ``the''
$\i$-category of small spectral categories.

\subsection{Stable envelopes of spectral categories}
Given any spectral category $\aC$, we can produce an $\i$-category by
passing to the associated simplicial category, fibrantly replacing, and applying the simplicial nerve to obtain $\N(\Omega^\i(\aC)^{\textrm{fib}})$.
This process yields a functor
\[
\Spcat \too \sset,
\]
from the category $\Spcat$ of small spectral categories to the
category of simplicial sets.  Precomposing with the functors
$\widehat{(-)}_{\perf}$ and $\widehat{(-)}_{\tri}$, we obtain functors 
\[
\psi_{\tri},\psi_{\perf} \colon \Spcat\too \sset
\]
and a natural transformation $\psi_{\tri} \to \psi_{\perf}$.  First,
we observe that these functors are compatible with the weak
equivalences of Theorem~\ref{thm:Quillenstable}.

\begin{lemma}\label{lem:preservesequiv}
Let $\aA$ and $\aB$ be small spectral categories, and let $f \colon
\aA \to \aB$ be a DK-equivalence.  Then the induced maps
$\psi_{\tri}(f)$ and $\psi_{\perf}(f)$ are categorical equivalences of
simplicial sets.
\end{lemma}

\begin{proof}
If $f\colon \aA\to\aB$ is a DK-equivalence, then one can check that
$(f_!,f^*)$ gives a Quillen equivalence between the spectral model
categories $\widehat{\aA}$ of $\aA$-modules and the spectral model
category $\widehat{\aB}$ of $\aB$-modules.  Passing to underlying
simplicial categories of cofibrant and fibrant objects, we see that
$\Omega^\i(\widehat{\aA})^{\cf}=\Mod(\aA)^{\cf}$ and
$\Omega^\i(\widehat{\aB})=\Mod(\aB)^{\cf}$ are DK-equivalent
simplicial categories.  Finally, applying the simplicial nerve yields
categorically equivalent simplicial sets.  Restricting to various full
subcategories yields the result for $\psi_{\tri}(f)$ and
$\psi_{\perf}(f)$.
\end{proof}

Therefore, we have induced functors $\Psi_{\tri}$ and $\Psi_{\perf}$
connecting $\N((\Spcat)^{\mathrm{c}})[W^{-1}]$ and
$\N((\sset)^{\mathrm{c}})[W^{-1}]$, equipped with a natural
transformation connecting them:
\[
\Psi_{\tri}\to \Psi_{\perf} \colon \N((\Spcat)^{\mathrm{c}})[W^{-1}] \too
\N((\sset)^{\mathrm{c}})[W^{-1}] \simeq \icat.
\]
In fact, by construction these functors preserve triangulated and
Morita equivalences respectively.  Furthermore, $\Psi_{\tri}$ lands in
small stable $\i$-categories and $\Psi_{\perf}$ lands in
idempotent-complete stable $\i$-categories.

\begin{lemma}\label{lem:sptostab}
$\Psi_{\tri} \aC$ factors through the subcategory $\stabcat\subset\icat$, and $\Psi_{\perf} \aC$ factors through the subcategory $\idemstabcat\subset\icat$.
\end{lemma}

\begin{proof}
As noted in remark~\ref{rem:pretriangstab}, since $\psi_{\tri}$ is the
underlying simplicial category associated to a pretriangulated
spectral category, $\Psi_{\perf}$ is characterized as the idempotent-completion of $\psi_{\tri}$ given by proposition~\ref{prop:idempotentchar} coupled with
corollary~\ref{cor:toencor}. Finally, note that maps of spectral categories induce, by left Kan extension, finite colimit-preserving on the level of stable $\i$-categories.
\end{proof}

Consequently, we may regard $\Psi_{\tri}$ as a functor $\N((\Spcat)^{\mathrm{c}})[W^{-1}]\to\stabcat$ and $\Psi_{\perf}$ as a functor $\N((\Spcat)^{\mathrm{c}})[W^{-1}]\to\idemstabcat$.

\begin{remark}\label{rem:modelloc}
Using the machinery of combinatorial simplicial model categories, we
can also localize the combinatorial model structure of
corollary~\ref{cor:Quillencombi} on $\Spcat$ at the triangulated or
Morita equivalences directly to obtain ``triangulated'' or ``Morita''
simplicial model categories on small spectral categories and then pass
to simplicial nerves; this is equivalent to localizing the
$\i$-category $\N((\Spcat)^{\cf})$.
\end{remark}

The content of theorem~\ref{thm:models} is that these functors are
equivalences.  We prove this theorem by producing an ``inverse'' to
$\Psi_{\tri}$ and $\Psi_{\perf}$ such that the composite is a
localization functor on $\N((\Spcat)^{\mathrm{c}})[W^{-1}]$.  We begin
with the following definition:

\begin{definition}
A simplicial category $\aA$ is {\em stable} if the simplicial nerve of
a fibrant replacement of $\aA$ is a stable $\i$-category.  A spectral
category $\aA$ is {\em stable} if its underlying simplicial category
$\Omega^\i \aA$ is stable. 
\end{definition}

We have the following characterization of equivalences between stable
spectral categories (see also \cite[5.7]{BM2}).

\begin{proposition}\label{prop:stabspecDK}
Let $\aA$ and $\aB$ be stable spectral categories.
Then a spectral functor $f\colon \aA\to\aB$ is a DK-equivalence if and only
if 
\[
\Ho(\Omega^\i f) \colon \Ho(\Omega^\i\aA)\to\Ho(\Omega^\i\aB)
\] 
is an equivalence.
\end{proposition}

\begin{proof}
Certainly essential surjectivity is determined on the level of the
homotopy category, so it suffices to show that, for all pairs of
objects $a$ and $b$ of $\aA$, $\pi_n\Map(a,b)\to\pi_n\Map(fa,fb)$ for
all integers $n$ whenever this is the case for $n=0$. 
Since $\aA$ and $\aB$ are stable,
\[
\pi_0\Map(\Sigma^n
a,b)\cong\pi_n\Map(a,b)\to\pi_n\Map(fa,fb)\cong\pi_0\Map(\Sigma^n fa,fb), 
\]
so this is immediate.
\end{proof}

We write $\Cat_\T^{\ex}$ for the simplicial category of small
stable simplicial categories.  We can model this as the subcategory of
the simplicial category $L^H(\Cat_\T)$ of small simplicial categories
where the objects are the stable simplicial categories and the mapping
spaces are computed by restriction of vertices to those simplicial functors which represent exact functors upon passage to the simplicial nerve.

\begin{proposition}
The $\i$-category obtained by applying the simplicial nerve to (a
fibrant replacement of) $\Cat_\T^{\ex}$ is equivalent to the
$\i$-category $\stabcat$.  That is, the equivalence (induced by the
simplicial nerve \cite[2.2.0.1]{HTT}) 
\[
\N((L^H(\Cat_\T)^{\mathrm{fib}}) \to\Cat_\i
\]
restricts to an equivalence
\[
\N((\Cat_\T^{\ex})^{\mathrm{fib}}) \to\Cat_\i^{\ex}.
\]

\end{proposition}

\begin{proof}
It suffices to show that the mapping spaces in $\Cat_\T^{\ex}$
have the correct homotopy type, and this follows from the comparison
between the mapping spaces of $\Cat_\T$ and
$\icat$ \cite[2.2.0.1]{HTT} and the fact that on both sides we define
the mapping spaces by the same restriction of vertices.
\end{proof}

\subsection{Spectral enrichment of stable $\i$-categories}
The mapping spaces of a stable $\i$-category $\aC$ are naturally the underlying spaces of mapping {\em spectra}, as discussed in section~\ref{sec:stabilization}.
We now use this to construct a cofibrant and fibrant
spectral category $\Upsilon(\C)$ whose underlying $\infty$-category
$\N(\Omega^\infty\Upsilon(\C))$ is equivalent to $\C$.  Indeed, the
simplicial category of presheaves of spectra
\[
\Fun_\Delta(\mathfrak{C}[\C]^{\op},\cS)
\]
on the associated (cofibrant) simplicial category $\mathfrak{C}[\C]$
is simultaneously a simplicial model category as well as a spectral
category, where the spectral enrichment is inherited from the
spectral structure on $\cS$ itself.
Moreover, we have an equivalence
\[
\N(\Fun_\Delta(\mathfrak{C}[\C]^{\op},\cS)^{\cf})[W^{-1}]\simeq\Fun(\C^{\op},\ispec),
\]
so it makes sense to ask whether or not a given presheaf of spectra is
stably representable (in the underlying $\infty$-category
$\Fun(\C^{\op},\ispec)$). 

\begin{definition}
Let
\[
\spect(\C) \subset \Fun_\Delta(\mathfrak{C}[\C]^{\op},\cS)
\]
denote the full spectral subcategory spanned by those (projectively)
cofibrant and fibrant functors which are stably representable.
\end{definition}

\begin{proposition}\label{prop:Nomega}
For a small stable $\i$-category $\C$, there is a natural equivalence
of $\infty$-categories $\C\to \N(\Omega^\infty \spect(\C))$.
\end{proposition}

\begin{proof}
The spectral Yoneda embedding
\[
\C\to\Fun(\C^{\op},\ispec)\simeq\N(\Fun_\Delta(\mathfrak{C}[\C]^{\op},\cS)^{\cf})
\]
is adjoint to a simplicial functor
\[
\mathfrak{C}[\C]\to\Fun_\Delta(\mathfrak{C}[\C]^{\op},\cS)^{\cf})
\]
which evidently factors through the full simplicial subcategory
$\Omega^\infty \spect(\C)$ spanned by the stably representable functors. 
The map $\C\to\N(\Omega^\infty \spect(\C))$ is the adjoint of the resulting
map $\mathfrak{C}[\C]\to\Omega^\infty \spect(\C)$. 

To see that this map is an equivalence, we observe first that it is
essentially surjective: indeed, a stably representable cofibrant and
fibrant functor $X\colon\mathfrak{C}[\C]^{\op}\to\cS$ is necessarily of the form
$X\simeq\Map(-,A)$ for some spectrum object $A=\{a_i\}$ of
$\C\simeq\N(\mathfrak{C}[\C]^{\mathrm{fib}})$.  Since $\C$ is stable,
$a_i\simeq\Sigma^i a$ for some object $a$ of $\C$, so $X$ is in the
image of $\C$ (which sends $a$ to the presheaf represented by
$\Sigma^\infty a$).  This map is also fully faithful, because if $a$
and $b$ are any pair of objects of $\C$, then
\[
\map(\Sigma^\i b,\Sigma^\i a)\simeq\map(b,\Omega^\i\Sigma^\i a)\simeq\map(b,a)
\]
since $a\simeq\Omega^\i\Sigma^\i a$.
\end{proof}

We have the following description of $\spect$ in terms of the stable
Yoneda embedding.

\begin{proposition}\label{prop:stabrepisexact}
Let $\C$ be a small stable $\i$-category.
The fully-faithful inclusion
\[
\spect(\C)\to\Fun_\Delta(\mathfrak{C}[\C]^{\op},\cS)
\]
factors, on the level of underlying $\i$-categories, as the composite
\begin{align*}
\N(\Omega^\i\spect(\C))\simeq\C\rightarrow
 \Fun^{\ex}(\C^{\op},\ispec)&^\omega\\
 \subseteq\Fun(\C^{\op},\ispec)&\simeq\N\Fun_\Delta(\mathfrak{C}[\C]^{\op},\cS)^{\cf}[W^{-1}].
\end{align*}
\end{proposition}

\begin{proof}
Any stably representable functor $\C^{\op}\to\ispec$ is exact, giving
the factorization  
\[
\C\to\Fun^{\ex}(\C^{\op},\ispec)\subseteq\Fun(\C^{\op},\ispec).
\]
By proposition \ref{prop:idempotentchar}, we may rewrite this as
$\C\to\Ind(\C)\simeq\Fun^{\ex}(\C^{\op},\ispec)$ to see that, as an
exact functor $\C^{\op}\to\ispec$, any stably representable functor is
also compact. 
\end{proof}

Our model of $\Cat^{\ex}_\T$ allows us to check that the construction
of $\spect$ induces a simplicial functor:  

\begin{proposition}\label{prop:specfunc}
The assignment which associates to the stable simplicial category $\C$
the spectral category $\spect(\C)$ defines a simplicial functor
\[
\spect \colon \Cat^{\ex}_\T \too L^H(\Spcat)
\]
and hence a functor of $\i$-categories
\[
\N(\spect) \colon \stabcat \too \N((L^H(\Spcat))^{\mathrm{fib}}).
\]
\end{proposition}

\begin{proof}
We first check that the construction of $\spect$ induces a
functor $\Cat^{\ex}_\aT \to \Spcat$.  Let $f\colon\aC\to\aD$ be a map of
stable simplicial categories and write
\[
f_!^\mathrm{cf}\colon\Fun_\Delta(\C^{\op},\cS)^\mathrm{cf}\to\Fun_\Delta(\aD^{\op},\cS)^\mathrm{cf}\]
for
the induced spectral functor.  Suppose that $X\colon\C^{\op}\to \cS$ is
projectively cofibrant and fibrant and that $\N(X)\colon\N(\C)^{\op}\to\ispec$ is
stably representable via the spectrum object $A=\{a_i\}$ in $\N\C$.
Since the diagram
\[
\xymatrix{
\N\C\ar[r]\ar[d] & \N\aD\ar[d]\\
\Fun(\N\C^{\op},\ispec)\ar[r] & \Fun(\N\aD^{\op},\ispec)}
\]
commutes (where the vertical maps are the {\em stable} Yoneda
embeddings), we see that $f_!$ restricts to a spectral functor
$\spect(\aC)\to\spect(\aD)$.

To verify that $\spect$ induces a simplicial functor, we must check
that it preserves equivalences of stable simplicial categories.  So
suppose that $f\colon\C\to\aD$ is an equivalence of stable simplicial
categories.  Then it follows that $f_!^\mathrm{cf}$ is a
DK-equivalence of spectral categories, as is its restriction to the
stably representable objects.
\end{proof}

\begin{proposition}\label{prop:ex/spec}
Let $\cA$ be a spectral category.
Then there are natural equivalences of compactly-generated stable $\i$-categories
\[
\N(\Fun_\cS(\cA^{\op},\cS)^\mathrm{c})[W^{-1}]\simeq\Ind(\Psi_{\tri}\cA)\simeq\Fun^{\ex}(\Psi_{\tri}\cA^{\op},\ispec).
\]
\end{proposition}

\begin{proof}
The first equivalence follows from the definition of $\Psi_{\tri}\aA$ as the smallest stable subcategory of the stable $\i$-category $\N(\Fun_\cS(\cA^{\op},\cS)^\mathrm{c})[W^{-1}]$ containing the representables, together with the observations that $\N(\Fun_\cS(\cA^{\op},\cS)^\mathrm{c})[W^{-1}]$ is compactly generated with compact objects
\[
\Psi_{\perf}\cA\simeq\N(\Fun_\cS(\cA^{\op},\cS)^\mathrm{c})[W^{-1}]^\omega
\]
and $\Ind(\Psi_{\tri}\cA)\simeq\Ind(\Psi_{\perf}\cA)$ (as $\Ind$-categories are automatically idempotent-complete).
The second equivalence follows immediately from proposition \ref{prop:idempotentchar}.
\end{proof}

\subsection{The triangulated and Morita localizations}
Let
\[
M\colon \N((\Spcat)^\mathrm{c})[W^{-1}] \to \N((\Spcat)^\mathrm{c})[W^{-1}]
\]
denote the composite functor 
\begin{equation}\label{eq:mcomp}
\xymatrix{
\N((\Spcat)^\mathrm{c})[W^{-1}] \ar[r]^-{\Psi_{\tri}}
& \stabcat \ar[r]^-{\N(\spect)}
& \N((L^H \Spcat)^{\textrm{fib}}) \htp \N((\Spcat)^\mathrm{c})[W^{-1}].\\ 
}
\end{equation}
As the previous proposition suggests, $M\cA$ is essentially the same as the pretriangulated spectral closure $\widehat{\aA}_{\tri}$ of $\cA$.

\begin{proposition}\label{prop:spec=tri}
There is an equivalence
\[
\widehat{\cA}_{\tri}\simeq M\cA
\]
in $\N((\Spcat)^\mathrm{c})[W^{-1}]$, natural in spectral categories $\cA$.
\end{proposition}

\begin{proof}
By proposition \ref{prop:ex/spec}, we have natural equivalences
\[
\N(\Fun_\cS(\cA^{\op},\cS)^\mathrm{c})[W^{-1}]^\omega\simeq\Psi_{\perf}\cA\simeq\Fun^{\ex}(\Psi_{\tri}\cA^{\op},\ispec)^\omega.
\]
These allows us to identify the smallest stable subcategory $\Psi_{\tri}\cA\subseteq\Psi_{\perf}\cA$ spanned by the representable functors $\widehat{a}\colon\cA^{\op}\to\cS$ with the stably representable functors $\Sigma^\i\widehat{a}\colon\Psi_{\tri}\aA^{\op}\to\ispec$.
\end{proof}

There is a natural transformation $\eta\colon\id\to M$ which can be
described as follows:  
On the level of spectral categories, the Yoneda embedding
\[
\cA\to\widehat{\cA}
\]
factors through the inclusion of the essentially small full spectral subcategory
\[
\widehat{\cA}_{\tri}\to\widehat{\cA}.
\]
The result is a natural transformation $\id\to\widehat{(-)}_{\tri}$ of endofunctors of $\Spcat$.
By proposition \ref{prop:spec=tri}, there is a natural equivalence
\[
\widehat{\cA}_{\tri}\overset{\sim}{\longrightarrow}\spect(\Psi_{\tri}\cA)=M\cA
\]
in $\N((\Spcat)^{\mathrm{c}})[W^{-1}]$; composing with this natural equivalence
gives the desired natural transformation $\eta\colon\id\to M$. 

\begin{proposition}\label{ff}
For any spectral category $\aA$,
$\eta_{\aA}\colon\aA\to M\aA$ is fully faithful.
\end{proposition}

\begin{proof}
By Yoneda's lemma, mapping spectra in $M\aA$ between stably representable objects are given by mapping spectra between the representing spectrum objects, giving an equivalence
\[
\Map_{M\aA}(\eta_{\aA}(a),\eta_{\aA}(b))_n\simeq\map_{\Psi_{\tri}\aA}(\widehat{a},{\Sigma^n}{\widehat{b}}).
\]
Since $\Psi_{\tri}\aA$ is a stable $\i$-category of spectral functors, Yoneda's lemma also gives an equivalence
\[
\Map_{\aA}(a,b)_n\simeq\map_{\Psi_{\tri}\aA}(\widehat{a},{\Sigma^n}{\widehat{b}}).
\]
Hence $\Map_{\aA}(a,b)\simeq\Map_{M\aA}(\eta_{\aA}(a),\eta_{\aA}(b))$.
\end{proof}

\begin{proposition}\label{es}
The functor $\eta_{\aA}\colon{\aA}\to M\aA$ is essentially surjective if and
only if ${\aA}$ is stable.
\end{proposition}

\begin{proof}
Indeed, ${\aA}\to M{\aA}$ is essentially surjective if and only if $\Omega^\i
\aA\to\Omega^\i M\aA$ is essentially surjective, which is the case if and
only if $\aA$ is already stable.
\end{proof}

Combining propositions~\ref{ff} and~\ref{es}, we obtain the following
corollary.

\begin{corollary}\label{Mequiv}
The functor $\eta_{\aA}\colon\aA\to M\aA$ is an equivalence of spectral
categories if and only if $\aA$ is a stable spectral category.
\end{corollary}

Next, we want to verify that $M$ is a localization.

\begin{proposition}\label{commutes}
The pair of natural transformations $\eta_{M\aA},M\eta_{\aA}\colon
M\aA\to M^2\aA$ induce a homotopy commutative square
\[
\xymatrix{
\aA\ar[r]^{\eta_{\aA}}\ar[d]_{\eta_{\aA}} & M\aA\ar[d]^{M_{\eta_{\aA}}}\\
M\aA\ar[r]_{\eta_{M\aA}}            & M^2 \aA.}
\]
\end{proposition}

\begin{proof}
First note that $M\eta_{\aA}\colon M\aA\to M^2 \aA$ sends
$x\colon\widehat{\aA}\to\ispec$ to 
the functor $\widehat{\eta_{\aA}}_! x\colon\widehat{M\aA}\to\ispec$ induced by
homotopy left Kan extension along
$\widehat{\eta_{\aA}}:\widehat{\aA}\to\widehat{M\aA}$.
If $x=\Map_A(-,a)$ is represented by the object $a$ of $\aA$, then the
universal properties of representable functors and homotopy left Kan
extensions force an equivalence $\widehat{\eta_{\aA}}_!
x\cong\Map(-,\eta_{\aA}(a))$, so that $\widehat{\eta_{\aA}}_! x$ is
represented by $\eta_{\aA}(a)$.  It follows that the restrictions of
$\eta_{M{\aA}}$ and $M\eta_{\aA}$ to $\aA$ are equivalent.
\end{proof}

\begin{corollary}\label{etaM=Meta}
The spectral functors $\eta_{M\aA}$ and $M\eta_{\aA}$ are equivalent.
In particular, both $\eta_{M\aA}$ and $M\eta_{\aA}$ are equivalences.
\end{corollary}

\begin{proof}
Since $\aA$ generates $M\aA$ under finite homotopy colimits and
desuspensions, it suffices to show that $M\eta_{\aA}$ preserves finite
homotopy colimits and desuspensions.  The fact that $M\eta_{\aA}$
preserves finite homotopy colimits follows from the fact that
$M\eta_{\aA}$ is a homotopy left Kan extension along $\widehat{\eta_{\aA}}$.
But suspension is an example of a finite homotopy colimit, so we have
that $M\eta_{\aA}(\Sigma x)\simeq\Sigma M\eta_{\aA}(x)$.  Hence
$M\eta_{\aA}(x)\simeq\Sigma M\eta_{\aA}(\Sigma^{-1} x)$, and as $M^2 \aA$ is
stable we see that $M\eta_{\aA}(\Sigma^{-1} x)\simeq\Sigma^{-1}
M\eta_{\aA}(x)$.  The final statement is a consequence of corollary
\ref{Mequiv} and the fact that $M\aA$ is a stable spectral category.
\end{proof}

\begin{corollary}\label{local}
The functor $M$ defines a localization of $\N((\Spcat)^{\mathrm{c}})[W^{-1}]$
with essential image the stable spectral categories.
\end{corollary}

\begin{proof}
This follows from the previous proposition and corollary by \cite[5.2.7.4]{HTT}.
\end{proof}

To see that we have an accessible localization, we need the following
proposition:

\begin{proposition}\label{prop:specfiltered}
Let $\aC\simeq\colim_i\aC_i$ be a filtered colimit of stable
$\i$-categories.  Then there is an equivalence of spectral categories
\[
\spect(\aC)\simeq\colim_i\spect(\aC_i).
\] 
\end{proposition}

\begin{proof}
We must show that the natural map
\[
\colim_i\spect(\aC_i)\to\spect(\aC)
\]
is a DK-equivalence of spectral categories.
Since $\aC$ and the $\aC_i$ are all stable spectral categories and
$\Omega^\infty$ and $\N$ commute with filtered colimits, this follows
from propositions \ref{prop:stabspecDK} and \ref{prop:Nomega}.
\end{proof}

Following Definition~\ref{def:Morita}, we make the following definitions.

\begin{definition}\label{def:imorita}
A map of small spectral $\i$-categories $f:\aA\to\aB$ is:
\begin{itemize}
\item A {\em triangulated equivalence} if
  $\Psi_{\tri}{f}:\Psi_{\tri}{\aA}\to\Psi_{\tri}{\aB}$ is an
  equivalence of (stable) $\i$-categories, and
\item A {\em Morita equivalence} if $\Psi_{\perf}f \colon \Psi_{\perf}
  \aA \to \Psi_{\perf} \aB$ is an equivalence of (idempotent-complete)
  stable $\i$-categories.
\end{itemize}
\end{definition}

Assembling the work of this section we obtain the following two results:

\begin{theorem}\label{thm:stable}
The functor
\[
\Psi_{\tri}(-) \colon \N((\Spcat)^{\mathrm{c}})[W^{-1}] \to\stabcat
\]
admits a fully faithful and accessible right adjoint
\[
\spect \colon \stabcat \to \N((\Spcat)^{\mathrm{c}})[W^{-1}].
\]
That is, the $\i$-category of stable $\i$-categories is an accessible 
localization of the $\i$-category of spectral categories obtained by
inverting the triangulated equivalences.
\end{theorem}

\begin{proof}
The follows from the factorization of $M$ given in
equation~\ref{eq:mcomp} and proposition~\ref{prop:specfiltered}.
\end{proof}

Recall that we have a stable idempotent completion functor
$\Idem\colon \Cat_\i^{\ex}\to\idemstabcat$.  Since $\Idem$ is left
adjoint to the (fully faithful) inclusion
$\idemstabcat\to\Cat_\i^{\ex}$, $\idemstabcat$ is the
localization of $\Cat_\i^{\ex}$ obtain by inverting idempotent
completion maps.  Further, recall that there is an equivalence $\Idem
\circ \Psi_{\tri} \htp \Psi_{\perf}$.

\begin{theorem}\label{thm:morita}
The functor
\[
\Psi_{\perf} \colon \N((\Spcat)^{\mathrm{c}})[W^{-1}] \to\idemstabcat
\]
admits a fully faithful and accessible right adjoint
\[
\spect \colon \idemstabcat\to \N((\Spcat)^{\mathrm{c}})[W^{-1}].
\]
That is, the $\i$-category of idempotent-complete stable
$\i$-categories is an accessible localization of the $\i$-category of
spectral categories obtained by inverting the Morita equivalences.
\end{theorem}

\begin{proof}
The $\i$-category of idempotent-complete stable $\i$-categories is a
localizing subcategory of the $\i$-category of stable $\i$-categories,
and idempotent-completion is an accessible functor as the
inclusion $\Psi_{\tri} \to \Psi_{\perf}$ preserves filtered colimits.
\end{proof}

\begin{remark}
Theorems~\ref{thm:stable} and \ref{thm:morita} imply that computing
the localizations of the model category structure on spectral
categories from Corollary~\ref{cor:Quillencombi} at the triangulated
and Morita equivalences (as discussed in Remark~\ref{rem:modelloc})
and passing to the simplicial nerve also yields the $\i$-categories
$\stabcat$ and $\idemstabcat$ respectively.
\end{remark}

We conclude the section with the promised applications of the theory.
First, the fact that we have accessible localizations provides the
following corollary about the structure of $\stabcat$ and
$\idemstabcat$.

\begin{corollary}\label{cor:cgen}
The $\i$-categories $\stabcat$ and $\idemstabcat$ are
compactly generated, complete, and cocomplete.
\end{corollary}

We will use the comparison above to lift small stable $\i$-categories
to spectral categories.  To this end, we make the following definition.

\begin{definition}\label{def:rep}
Let $\aA$ and $\aB$ be small idempotent-complete stable
$\i$-categories.  We write $\rep(\aB,\aA)=\spect(\Fun^{\ex}(\aB,\aA))$
for the small pretriangulated spectral category associated to the
small stable $\i$-category of exact functors from $\aB$ to $\aA$.
\end{definition}

\begin{corollary}\label{cor:rep}
Let $\aA$ and $\aB$ be idempotent-complete small stable
$\i$-categories and let $\spect(A)$ and $\spect(B)$
be spectral categories lifting $\aA$ and $\aB$.
Then $\N(\rep(\aA,\aB))\simeq\idemfun\!(\aA,\!\aB)$ is equivalent to the
$\i$-category of right-compact
$\spect(A)^{\op}\sma\spect(B)$-modules.
\end{corollary}

\begin{proof}
This follows from theorem~\ref{thm:morita} and
corollary~\ref{cor:toencor}.
\end{proof}

We also record the following result concerning lifting diagrams of
small stable $\infty$-categories to diagrams of spectral categories. 

\begin{proposition}\label{prop:genlift}
Let $I$ be a small category.  Given a diagram $\cD$ of small stable
$\i$-categories indexed by $\N(I)$, there exists an $I$-diagram of
pretriangulated spectral categories $\widetilde{\cD}$ lifting $\cD$.
\end{proposition}

\begin{proof}
This is a consequence of \cite[4.2.4.4]{HTT}.  Given a diagram of
small stable $\i$-categories, the equivalence in
theorem~\ref{thm:stable} gives rise to a diagram in the localization
of $\N((\Spcat)^{\mathrm{c}})[W^{-1}]$.  Including the localization
into $\N((\Spcat)^{\mathrm{c}})[W^{-1}]$, we now obtain a diagram in
$\N((\Spcat)^{\mathrm{c}})[W^{-1}] \htp \N((\Spcat)^{\cf})$ and we can
use \cite[4.2.4.4]{HTT} to lift this to a rigid diagram in $\Spcat$.
\end{proof}

\section{Exact sequences}\label{sec:exact}

In this section we discuss the various definitions of exact sequence,
relating notions for triangulated categories, spectral categories, and
stable $\i$-categories.  Arguably the most fundamental definition is
that of an exact sequence of triangulated categories, as it turns out
that exact sequences in both spectral and stable $\i$-categories can
be detected on the level of the homotopy category.
Recall that a sequence of triangulated categories
\[
\aA \too \aB \too \aC
\]
is called {\em exact} if the composite is zero, the functor $\aA \to
\aB$ is fully faithful, and the induced functor from the Verdier
quotient $\aB/\aA$ to $\aC$ is {\em cofinal}, i.e., it becomes an
equivalence after idempotent completion.  Said differently, a
triangulated functor $\aC'\to\aC$ is {\em cofinal} if every object of
$\aC$ is a summand of an object of $\aC'$.
The purpose of this section is to develop analogues of these notions
for stable $\i$-categories.  

\subsection{The Verdier quotient as the cofiber in $\idemstabcat$}
Let $\kappa$ denote an infinite regular cardinal.
We recall the
following terminology from \cite[\S 5.3.4]{HTT}.

\begin{definition}
Let $\aA$ be an $\i$-category.
We say that $\aA$ is {\em $\kappa$-cocomplete} if $\aA$ admits all
$\kappa$-small colimits. 
\end{definition}

Most of the small $\i$-categories which arise in this paper can be
realized as the full subcategory $\aC^\kappa\subset\aC$ of
$\kappa$-compact objects in a stable presentable $\i$-category $\aC$.
In this case, we can reconstruct $\aC$ itself as
$\Ind_\kappa(\aC^\kappa)$, which formally adjoins $\kappa$-filtered
colimits.  To make this precise, we recall the notions of
$\kappa$-filtered $\infty$-category, $\kappa$-filtered colimit, and
$\kappa$-continuous functor.

\begin{definition}
An $\i$-category $\aC$ is {\em $\kappa$-filtered} if every map
$K\to\aC$ from a $\kappa$-small simplicial set $K$ extends to a
functor $K^{\triangleright}\to\aC$ (see \cite[2.1.4.2]{HTT}
for the cone notation).  A simplicial set $K$ is $\kappa$-filtered if
there exists a categorical equivalence $K\to\aC$ for some
$\kappa$-filtered $\i$-category $\aC$.  Lastly, a $\kappa$-filtered
colimit is a colimit indexed by a $\kappa$-filtered simplicial set.
\end{definition}

\begin{definition}
Let $\aA$ and $\aB$ be $\i$-categories and let $f:\aA\to\aB$ be a
functor.  We say that $f$ is {\em $\kappa$-continuous} if $f$
preserves $\kappa$-filtered colimits. 
\end{definition}

We write $\kstabcat$ for the $\i$-category of small
$\kappa$-cocomplete stable $\i$-categories and $\kappa$-small
colimit-preserving functors 
thereof; note that if $\kappa>\omega$, any small $\kappa$-cocomplete
stable $\i$-category $\aA$ is necessarily idempotent
complete \cite[5.4.2.4]{HTT}.  Given a small $\kappa$-cocomplete
stable $\i$-category $\aA$, the $\i$-category $\Ind_\kappa(\aA)$ is a
$\kappa$-compactly generated stable $\i$-category such that
$\Idem(\aA) \simeq\Ind_\kappa(\aA)^\kappa$ \cite[5.5.7.8,
5.5.7.10]{HTT}.

In fact, provided $\kappa>\omega$, restriction to subcategories of
$\kappa$-compact objects determines an equivalence between the
$\i$-category of $\kappa$-compactly generated stable $\i$-categories
$\stabprcat_\kappa$ and the $\i$-category $\kstabcat$ of small
$\kappa$-cocomplete stable $\i$-categories, with inverse
$\Ind_\kappa$ \cite[5.5.7.10]{HTT}.  As a consequence,
corollary~\ref{cor:cgen} implies that the $\i$-category of
$\kappa$-compactly generated stable $\i$-categories is cocomplete.

We now define an analogue of the Verdier quotient of triangulated
categories on the level of stable $\i$-categories.  Specifically, if
$\aA\to\aB$ is a fully faithful functor of stable $\infty$-categories,
then $\Ho(\aA)\to\Ho(\aB)$ is a fully faithful functor of triangulated
categories, and we may form the usual Verdier quotient
$\Ho(\aB)/\Ho(\aA)$.  This is defined as the initial triangulated
category $\Ho(\aB)/\Ho(\aA)$ equipped with a triangulated functor
$\Ho(\aB)\to\Ho(\aB)/\Ho(\aA)$ such that the composite
$\Ho(\aA)\to\Ho(\aB)\to\Ho(\aB/\aA)$ is trivial \cite[2.1.8]{Neeman}.

\begin{definition}\label{defn:verdquot}
Let $f:\aA\to\aB$ be a fully faithful functor of presentable stable
$\infty$-categories (this means that $f$ preserves colimits).  The
{\em Verdier quotient} $\aB/\aA$ of $\aB$ by $\aA$ is the cofiber of
$f$ in the $\infty$-category $\stabprcat$ of presentable stable
$\i$-categories.
\end{definition}

It is useful to identify the Verdier quotient in terms of a Bousfield
localization; specifically, we will see that the Verdier quotient is
the Bousfield localization of $\aB$ at the arrows with cofiber in
$\aA$.

\begin{lemma}\label{lem:upobl}
Let $\aC$ be a presentable $\i$-category and $S$ be a strongly
saturated class of arrows of $\aC$.  Then $S$ is of small generation
if and only if the full subfunctor 
\[
\Fun^{\L}_S(\aC,-)\subseteq\Fun^{\L}(\aC,-) \colon \prcat\to\widehat{\mathrm{Cat}}_\i 
\]
of $\Fun^{\L}(\aC,-)$, spanned by those colimit-preserving functors
$\aC\to\aD$ which carry the arrows in $S$ to equivalences in $\aD$, is
corepresentable by a presentable $\i$-category $\aC'$.  Moreover, in
this case, $\aC'\simeq S^{-1}\aC$. 
\end{lemma}

\begin{proof}
If $S$ is of small generation then $S^{-1}\aC$ is presentable and
corepresents the functor $\Fun^{\L}_S(\aC,-)$
by \cite[5.5.4.14, 5.5.4.20]{HTT}.  Conversely, if
this functor is corepresentable by $\aC'$ then the identity
$\aC'\to\aC'$ determines a colimit-preserving functor $\aC\to\aC'$.
Let $T$ be the class of arrows in $\aC$ which become invertible in
$\aC'$, and note that $S\subseteq T$, $T$ is strongly
saturated \cite[5.5.4.10]{HTT}, and $T$ is of small
generation \cite[5.5.4.16]{HTT} (the last claim uses the
fact that the equivalences in $\aC'$ is the strongly saturated class
generated by the identity of the initial object of $\aC'$, which
follows from \cite[5.5.4.5, 5.5.4.6]{HTT}).  Thus
$T^{-1}\aC\simeq\aC'$, so $\aC'$ also corepresents the functor
$\Fun^{\L}_T(\aC,-)$, showing that a colimit-preserving functor
$\aC\to\aD$ inverts the arrows of $S$ if and only if it inverts the
arrows of $T$.  Since $S$ is strongly saturated, we conclude that $S=T$.
\end{proof}

The preceding lemma now allows us to characterize the cofiber as a
localization.

\begin{proposition}\label{prop:loc=cof}
Let $\aA\to\aB$ be a fully faithful functor of presentable stable
$\infty$-categories and let $S$ denote the collection of arrows in
$\aB$ whose cones lie in the essential image of $\aA$.  Then $S$ is a
strongly saturated class of maps in $\aB$ of small generation, and the
Verdier quotient $\aB/\aA$ is equivalent to the Bousfield localization
$S^{-1}\aB$. 
\end{proposition}

\begin{proof}
Let $\aC$ be a presentable stable $\i$-category, and note that a
colimit-preserving functor $\aB\to\aC$ sends the arrows in $S$ to
equivalences in $\aC$ if and only if its restriction to $\aA$ is
trivial. 
We therefore may identify
\[
\Fun^{\L}(\aB/\aA,\aC)\subseteq\Fun^{\L}(\aB,\aC),
\]
with the full subcategory spanned by those colimit-preserving functors
$\aB\to\aC$ which send the arrows in $S$ to equivalences in $\aC$.  It
follows from lemma~\ref{lem:upobl} that $\aB/\aA\simeq
T^{-1}\aB$, where $T$ is the strongly saturated class of arrows of
$\aB$ which become equivalences in $\aB/\aA$.

We now show that $S$ is strongly saturated, so that $S=T$.  First,
suppose given a cofiber sequence $X\to Y\to Z$ in $\aB$ such that $Z$
lies in the essential image of $\aA$, and let $X\to X'$ be any map.
Then the cofiber of $X'\to X'\coprod_X Y$ is equivalence to $Z$, so is
also in the essential image of $\aA$.  Second, given a diagram
$f_\alpha\colon X_\alpha\to Y_\alpha$ in $\Fun(\Delta^1,\aB)$ with
colimit $f\colon X\to Y$, and suppose that the cofibers $Z_\alpha$ of
each $f_\alpha$ lies in the essential image of $\aA$.  Commuting
colimits implies that the cofiber $Z$ of $f$ is computed as the
colimit of the $Z_\alpha$, and this lies in the essential image of
$\aA$ since $\aA$ is closed under colimits and the functor $\aA\to\aB$
preserves colimits.  Lastly, suppose $h=g\circ f$ is a composite of
$f\colon X\to Y$ followed by $g\colon Y\to Z$, and write $Y/X$, $Z/Y$,
and $Z/X$ for the cofibers of $f$, $g$, and $h$, respectively.  Then
we have a cofiber sequence $Y/X\to Z/X\to Z/Y$, so if any two lie in
the essential image of $\aA$ then so does the third.
\end{proof}

In fact, we can be more precise about a generating set for the local
equivalences:

\begin{proposition}\label{prop:initialquotient}
Let $i\colon \aA\to\aB$ be a fully faithful inclusion of
$\kappa$-compactly generated stable $\i$-categories which preserves
$\kappa$-compact objects, let $S$ be the (small) collection of arrows
of $\aB^\kappa$ whose cofibers lie in the image of $\aA^\kappa$, and
let $T$ be the (large) collection of arrows of $\aB$ whose cofibers
lie in the image of $\aA$.  Then the natural map
\[
S^{-1}\aB\to T^{-1}\aB \htp \aB[T^{-1}] \htp \aB/\aA
\]
is an equivalence of $\i$-categories, where here $S^{-1}\aB$ and
$T^{-1} \aB$ denote the subcategories of local objects.
\end{proposition}

\begin{proof}
Without loss of generality we may identify $\aA$ with its essential
image in $\aB$, so that an arrow $f\colon X\to Y$ is in $T$ if and only
if any cofiber $Z$ of $f$ lies in $\aA$.  By \cite[5.5.4.15]{HTT} it
suffices to show that $T\subseteq\overline{S}$, the 
strongly saturated class of arrows of $\aB$ generated by $S$
(see \cite[5.5.4.5]{HTT}).
To see this, let $X\overset{f}{\to} Y\overset{g}{\to} Z$ be a cofiber 
sequence in $\aB$ such that $Z$ is in $\aA$.  Then $Z=\colim_\alpha
Z_\alpha$ is a $\kappa$-filtered colimit of objects
$Z_\alpha\in\aA^\kappa\subset\aB^\kappa$ and $Y=\colim_\alpha
Y_\alpha$ is a $\kappa$-filtered colimit of objects
$Y_\alpha=Y\times_Z Z_\alpha$.
Now $Y_\alpha$ may not be $\kappa$-compact, so write
$Y_\alpha=\colim_\beta Y_{\alpha\beta}$ for some
$Y_{\alpha\beta}\in\aB^\kappa$ and consider the resulting diagram of
cofiber sequences 
\[
\xymatrix{
X_{\alpha\beta}\ar[r]^{f_{\alpha\beta}}\ar[d] & Y_{\alpha\beta}\ar[r]^{g_{\alpha\beta}}\ar[d] & Z_{\alpha\beta}\ar[d]\\
X_\alpha\ar[r]^{f_\alpha}\ar[d] & Y_\alpha\ar[r]^{g_\alpha}\ar[d] & Z_\alpha\ar[d]\\
X\ar[r]^f & Y\ar[r]^g & Z}
\]
in which the lower right and upper left squares are cartesian, which
implies that these two squares are also cocartesian and that the maps
$X_\alpha\to X$ and $Z_{\alpha\beta}\to Z_\alpha$ are equivalences.
Hence $Z_{\alpha\beta}\in\aA^\kappa\subseteq\aB^\kappa$ and we
conclude that $g_{\alpha\beta}$ and therefore $f_{\alpha\beta}$ as
well are maps in $\aB^\kappa$; in particular, $f_{\alpha\beta}$ is an
arrow in $S$. 
It follows from \cite[5.5.4.5]{HTT} that the pushout $f_\alpha$ of
$f_{\alpha\beta}$ along $X_{\alpha\beta}\to X_\alpha\simeq X$ is in
$\overline{S}$, and we see from \cite[5.5.4.12]{HTT} that
$f\simeq\colim_\alpha f_\alpha:X\simeq\colim_\alpha X_\alpha\to\colim
Y_\alpha\simeq Y$ is then also in $\overline{S}$.  
\end{proof}

Definition~\ref{defn:verdquot} leads to the following definition of an
exact sequence.

\begin{definition}\label{defn:exact}
A sequence of presentable stable $\i$-categories
$\aA\to\aB\to\aC$ is {\em exact} if the composite is trivial,
$\aA\to\aB$ is fully faithful, and the map $\aB/\aA\to\aC$ is an
equivalence.
\end{definition}

Somewhat surprisingly, as a consequence of the hypothesis of stability
we can detect exact sequences on the level of homotopy categories,
despite the fact that functors which are fully faithful on homotopy
categories are not typically fully faithful as functors of
$\i$-categories.  The following proposition connects the
$\i$-categorical Verdier quotient of definition~\ref{defn:verdquot} to
the Verdier quotient of the triangulated homotopy categories.

\begin{proposition}\label{prop:quotient}
Let $\aA\to\aB$ be a fully faithful inclusion of presentable stable
$\i$-categories.  Then the natural map
$\Ho(\aB)/\Ho(\aA)\to\Ho(\aB/\aA)$ is an equivalence. 
\end{proposition}

\begin{proof}
By construction, $\aB/\aA\subseteq\aB$ is the full subcategory on
those objects $b$ such that $\map(a,b)\simeq *$ for all objects $a$ in
the image of $\aA$. This shows that, as full subcategories of
$\Ho(\aB)$, $\Ho(\aB/\aA)\subseteq\Ho(\aB)/\Ho(\aA)$. Conversely, if
$b$ is in $\Ho(\aB)/\Ho(\aA)$, then $\pi_0\map(a,b)\simeq\ast$ for
each object $a$ in the image of $\aA$, and we claim that in fact
$\map(a,b)\simeq *$. Indeed, $\aA$ is a stable subcategory of $\aB$,
so that $\pi_n\map(a,b)\simeq\pi_0\map(\Sigma^n a,b)\simeq\ast$. Hence
$\Ho(\aB)/\Ho(\aA)\subseteq\Ho(\aB/\aA)$ as well. 
\end{proof}

The argument for the previous proposition also implies the following
characterization of fully faithful maps; note that here we do not need
the hypothesis that the stable $\i$-categories are presentable, as we
are not working with localizations.

\begin{proposition}\label{prop:full}
A map of stable $\i$-categories $\aA\to\aB$ is fully faithful if and
only if $\Ho(\aA)\to\Ho(\aB)$ is fully faithful.
\end{proposition}

\begin{corollary}\label{cor:equiv}
A map of stable $\i$-categories $\aA\to\aB$ is an equivalence if and
only if $\Ho(\aA)\to\Ho(\aB)$ is an equivalence. 
\end{corollary}

As we are predominantly interested in sequences of small
$\i$-categories, we will now extend definition~\ref{defn:exact} to the
$\i$-category $\kstabcat$.

\begin{definition}\label{def:exact}
A sequence of $\kappa$-cocomplete small stable $\i$-categories and
$\kappa$-small colimit preserving functors $\aA \to \aB \to \aC$ is
{\em exact} if the sequence
\[
\Ind_\kappa(\aA)\to\Ind_\kappa(\aB)\to\Ind_\kappa(\aC)
\]
is an exact sequence of presentable stable $\i$-categories.
\end{definition}

Although we've defined exact sequences in $\kstabcat$ to be those
sequences which are exact in $\stabprcat$, we can give an intrinsic 
description.  Just as in the presentable case, the quotient $\aB/\aA$
will denote the cofiber of the fully faithful inclusion $\aA\to\aB$
of $\kappa$-cocomplete small stable $\i$-categories.

\begin{proposition}\label{prop:exact}
A sequence of $\kappa$-cocomplete small stable $\i$-categories and
$\kappa$-small colimit preserving functors $\aA \to \aB \to \aC$ is
exact if and only if the composite is trivial, $\aA\to\aB$ is fully
faithful, and the resulting map $\aB/\aA\to\aC$ is an equivalence
(after idempotent completion if $\kappa=\omega$).
\end{proposition}

\begin{proof}
The fully faithful inclusions $\aA\subset\Ind_\kappa(\aA)$ and
$\aB\subset\Ind_\kappa(\aB)$ show that $\aA\to\aB$ is fully faithful
if and only if $\Ind_\kappa(\aA)\to\Ind_\kappa(\aB)$ is fully faithful
(for the reverse direction, this follows from the definition of the
mapping spaces in $\Ind_{\kappa}(-)$).  Thus it remains to check that
$\aB/\aA\to\aC$ is an equivalence upon idempotent completion if and
only if
$\Ind_\kappa(\aB)/\Ind_\kappa(\aA)\simeq\Ind_\kappa(\aB/\aA)$. Since
$\Ind_\kappa$ preserves cofibers, it is enough to check that the
equivalence
$\Ind_\kappa(\aB)/\Ind_\kappa(\aA)\simeq\Ind_\kappa(\aB/\aA)$ implies
the equivalence $\aB/\aA\simeq\aC$ whenever the latter are idempotent
complete.  Thus, given a $\kappa$-cocomplete small stable
$\i$-category $\aD$ (which we assume is idempotent complete if
$\kappa=\omega$), we must show that 
\[
\Fun^{\ex(\kappa)}(\aC,\aD)\to\Fun^{\ex(\kappa)}(\aB,\aD)\to\Fun^{\ex(\kappa)}(\aA,\aD)
\]
is a fiber sequence of $\i$-categories.
Since $\aD\simeq(\Ind_\kappa(\aD))^\kappa$, by adjunction this is
equivalent to the sequence 
\[
\Fun^{\L}(\Ind_\kappa(\aC),\Ind_\kappa(\aD))\to\Fun^{\L}(\Ind_\kappa(\aB),\Ind_\kappa(\aD))\to\Fun^{\L}(\Ind_\kappa(\aA),\Ind_\kappa(\aD)),
\]
which is a fiber sequence by assumption.
\end{proof}

\subsection{The Thomason-Neeman localization theorem}
In fact, we can further reduce to a criterion on the level of homotopy
categories.  For this, we need the following proposition.  In the
proof, we take advantage of the detailed study of localization in the
context of (well-generated) triangulated categories by
Neeman~\cite{Neeman, NeemanLoc} and Krause~\cite{Krausenotes} and the
fact that the homotopy category of a presentable stable $\i$-category
is a well-generated triangulated category (see~\cite[1.4.5.2]{HAG}
and~\cite{Krausewell}).

\begin{proposition}\label{prop:kappaquotient}
Let $\aA\to\aB$ be a fully faithful and $\kappa$-small colimit
preserving functor of $\kappa$-cocomplete small stable
$\i$-categories.  Then the natural map
\[
\Ho(\aB)/\Ho(\aA)\to\Ho(\aB/\aA)
\]
is an equivalence.  In other words, the functor $\Ho(-)$ preserves quotients of fully faithful functors.
\end{proposition}

\begin{proof}
We have equivalences
\[
\Ho(\Ind_\kappa(\aB))/\Ho(\Ind_\kappa(\aA))\simeq\Ho(\Ind_\kappa(\aB)/\Ind_\kappa(\aA))\simeq\Ho(\Ind_\kappa(\aB/\aA)),
\]
where the first equivalence follows from
proposition~\ref{prop:quotient} and the last equivalence follows from
the fact that $\Ind_\kappa$ preserves cofibers.  We therefore obtain a
commutative (up to natural isomorphism) square 
\[
\xymatrix{
\Ho(\aB)/\Ho(\aA)\ar[r]\ar[d] & \Ho(\aB/\aA)\ar[d]\\
\Ho(\Ind_\kappa(\aB))/\Ho(\Ind_\kappa(\aA))\ar[r] & \Ho(\Ind_\kappa(\aB/\aA))}
\]
where the right vertical map is fully faithful and the bottom map
is an equivalence.  

To see that the top vertical map is fully faithful, we use 
Neeman's generalization of Thomason's localization theorem  
(see \cite[4.4.9]{Neeman} or ~\cite{NeemanLoc}) to show that the left
vertical map is fully faithful.  First, since $\Ind_{\kappa}(\aA)$ and
$\Ind_{\kappa}(\aB)$ are presentable, the criterion
of~\cite{Krausewell} (characterizing well-generated triangulated
categories) and~\cite[1.4.5.2]{HAG} imply that
$\Ho(\Ind_{\kappa}(\aB))$ is well-generated and (since the map 
$\Ind_{\kappa}(\aA) \to \Ind_{\kappa}(\aB)$ is fully-faithful) the
image of $\Ho(\Ind_{\kappa}(\aA))$ is a localizing subcategory
generated by a small set of objects.  Applying the form of Neeman's
theorem proved by Krause in~\cite[7.2.1]{Krausenotes} now implies that  
\[
\Ho\Ind_{\kappa}(\aB)^\kappa/\Ho\Ind_{\kappa}(\aA)^\kappa\to\Ho\Ind_{\kappa}(\aB)/\Ho\Ind_{\kappa}(\aA)
\]
is a fully faithful map.
Since~\cite[1.4.5.1]{HAG} implies that there is an equivalence
$\Ho(\Ind_{\kappa}(\aA)^{\kappa}) \htp \Ho(\Ind_{\kappa}(\aA))^{\kappa}$
and $\Ind_{\kappa}(\aA)^{\kappa} \htp \aA$ up to idempotent completion
(and similarly for $\aB$), we conclude that the left vertical map is
fully faithful.

Finally, this map is essentially surjective because there is a
commutative (up to natural isomorphism) triangle 
\[
\xymatrix{
& \Ho(\aB)\ar[ld]\ar[rd] & \\
\Ho(\aB)/\Ho(\aA)\ar[rr] & & \Ho(\aB/\aA)}
\]
such that both maps from $\Ho(\aB)$ are essentially surjective.
\end{proof}

Now we can obtain the following correspondence between exact sequences
of small stable $\i$-categories and exact sequences of triangulated
categories.

\begin{proposition}\label{prop:stabquot}
A sequence of $\kappa$-cocomplete small stable $\i$-categories and
$\kappa$-small colimit preserving functors $\aA\to\aB\to\aC$ is exact
if and only if the associated sequence
$\Ho(\aA)\to\Ho(\aB)\to\Ho(\aC)$ of triangulated categories is exact, in the sense that the composite is trivial, $\Ho(\aA)\to\Ho(\aB)$ is fully faithful, and the map $\Ho(\aB)/\Ho(\aA)\to\Ho(\aC)$ is an equivalence after idempotent completion. 
\end{proposition}

\begin{proof}
Suppose $\aA\to\aB\to\aC$ is exact.
Then the composite is trivial, $\aA\to\aB$ is fully faithful, and
$\aB/\aA\to\aC$ is an equivalence up to idempotent completion, and so
the same must be true on the level of triangulated homotopy
categories.  Thus it is enough to show that
$\Ho(\aB)/\Ho(\aA)\to\Ho(\aC)$ is an equivalence up to idempotent
completion, which follows from proposition \ref{prop:kappaquotient}.
Conversely, suppose that  
\[
\Ho(\aA)\to\Ho(\aB)\to\Ho(\aC)
\]
is exact.  Then $\aA\to\aB$ is fully faithful by
proposition \ref{prop:full}, and the equivalences
$\Ho(\aB/\aA)\simeq\Ho(\aB)/\Ho(\aA)\simeq\Ho(\aC)$ (the last up to
idempotent completion) implies that $\aB/\aA\simeq\aC$ by
corollary \ref{cor:equiv}.
\end{proof}

Finally, we record a technical proposition that is used in the context
of our construction of non-connective $K$-theory.  First, we need a
technical lemma about the behavior of the $\Ind$ functor.

\begin{lemma}\label{lem:indcpt}
Let $\aA\to\aB$ be an exact functor of small stable $\infty$-categories.
Then the induced map $\Ind(\aA)\to\Ind(\aB)$ of presentable stable
$\infty$-categories preserves $\kappa$-compact objects for all
infinite regular cardinals $\kappa$. 
\end{lemma}

\begin{proof}
Recall that a right adjoint preserves $\kappa$-filtered colimits if
and only if its left adjoint preserves $\kappa$-compact
objects~\cite[5.5.1.4]{HTT}.  Since functors which preserve filtered 
colimits also preserve $\kappa$-filtered colimits and $\Ind(-)$ takes exact
functors to functors which preserve compact objects, the result
follows.
\end{proof}

Note that in the statement of the following proposition, we implicitly
use the facts that stable $\i$-categories have all finite colimits and
exact functors preserve finite colimits.

\begin{proposition}\label{prop:notyet}
Let $\aA\to\aB\to\aC$ be an exact sequence of small stable
$\i$-categories.  Then for any infinite regular cardinal $\kappa$,
\[
\Ind(\aA)^\kappa\to\Ind(\aB)^\kappa\to\Ind(\aC)^\kappa
\]
is an exact sequence of idempotent-complete small stable
$\i$-categories.
\end{proposition}

\begin{proof}
First, by proposition~\ref{prop:stabquot}, it suffices to check that
\[ 
\Ho(\Ind(\aA)^\kappa) \to\Ho(\Ind(\aB)^\kappa)\to\Ho(\Ind(\aC)^\kappa)
\]
is an exact sequence of triangulated categories.  Again, we will
deduce this from Neeman's generalization of Thomason's localization
theorem  (see \cite[4.4.9]{Neeman} or ~\cite{NeemanLoc}), as follows.
First, observe that~\cite[1.4.5.1]{HAG} implies that there is an 
equivalence $\Ho(\Ind(\aA)^{\kappa}) \htp \Ho(\Ind(\aA))^{\kappa}$
(and analogous equivalences for the other terms in the sequence).
Next, since $\Ind(\aA)$ and $\Ind(\aB)$ are
presentable, the criterion of~\cite{Krausewell} (characterizing
well-generated triangulated categories) and~\cite[1.4.5.2]{HAG} imply
that $\Ho(\Ind(\aB))$ is well-generated and (since the map
$\Ind(\aA) \to \Ind(\aB)$ is fully-faithful) the image of
$\Ho(\Ind(\aA))$ is a localizing subcategory generated by a small set  
of objects.  Once again, the localization
theorem~\cite[7.2.1]{Krausenotes} implies that  
\[
\Ho(\Ind(\aB))^\kappa/\Ho(\Ind(\aA))^\kappa\to\Ho(\Ind(\aB/\aA))^\kappa
\]
is an equivalence up to idempotent completion.  The hypothesis that
$\aB/\aA \to \aC$ is an equivalence up to idempotent completion now
implies the result.
\end{proof}

\subsection{Split-exact sequences}
We will be particularly interested in exact sequences which are split in
the following sense.

\begin{definition}\label{def:splitexact}
An exact sequence of small $\kappa$-cocomplete stable $\i$-categories
and $\kappa$-small colimit preserving functors 
\[
\xymatrix{
\aA \ar[r]^{f} & \aB \ar[r]^{g} & \aC \\
}
\]
is called {\it split-exact} if there exist exact functors $i \colon \aB \to \aA$
and $j \colon \aC \to \aB$, right adjoint to $f$ and $g$, respectively,
such that $i \circ f \htp \Id$ and $g \circ j \htp \Id$ via the
adjunction morphisms.
\end{definition}

We will also be interested in (split-) exact sequences of spectral categories.

\begin{definition}
A sequence $\aA\to\aB\to\aC$ of spectral categories is exact if the
induced sequence of stable presentable $\i$-categories 
\[
\N(\Omega^\infty \Mod(\aA)^{\cf}) \to\N(\Omega^\infty \Mod(\aB)^{\cf}) \to\N(\Omega^\infty \Mod(\aC)^{\cf})
\]
is exact.
\end{definition}

The following characterization is an immediate corollary of proposition~\ref{prop:stabquot}.

\begin{proposition}
A sequence $\aA \to \aB \to \aC$ of spectral categories is exact if
and only if the induced sequence of triangulated categories
\[
\aD(\aA) \to \aD(\aB) \to \aD(\aC)
\]
is exact.
\end{proposition}

Next, observe that we can relate these notions as follows (the proof
of which is immediate):

\begin{proposition}
Let $\aA \to \aB \to \aC$ be an (split-) exact sequence of small spectral
categories.  Then
$\Psi_{\perf}(\aA) \to \Psi_{\perf}(\aB) \to \Psi_{\perf}(\aC)$ is a
(split-) exact sequence of small stable $\i$-categories.
\end{proposition}

We also have an essential converse statement.

\begin{proposition}\label{prop:speclift}
Let $\aA \to \aB \to \aC$ be a (split-) exact sequence of small stable
$\i$-categories.  Then there exists a (split-) exact sequence of small
stable spectral categories
\[
\widetilde{\aA} \to \widetilde{\aB} \to \widetilde{\aC}
\]
such that
$\Psi_{\perf}(\widetilde{\aA} \to \widetilde{\aB} \to \widetilde{\aC})$
is naturally equivalent to $\aA \to \aB \to \aC$.
\end{proposition}

\begin{proof}
This follows from proposition~\ref{prop:genlift}.
\end{proof}

\subsection{Approximating split-exact sequences}

In order to localize with respect to the (split-) exact sequences, we
need to be able to choose a set of representatives which generate them
under filtered colimits.

\begin{lemma}\label{lem:small}
The full subcategory
$(\idemstabcat)^{\omega} \subset \idemstabcat$ of compact small stable
idempotent-complete $\i$-categories is essentially small.
\end{lemma}

\begin{proof}
The result follows from the fact that $\idemstabcat$ is an accessible
localization of $\stabcat$, and $\stabcat$ itself is an
accessible localization of the finitely presentable $\i$-category of
small spectral categories via theorem~\ref{thm:models}.
\end{proof}

This has the following immediate and essential corollary:

\begin{corollary}\label{cor:setrep}
For any regular cardinal $\kappa$, there exists a set $\aE$ of
representatives of split-exact sequences of $\kappa$-compact small
idempotent-complete stable $\i$-categories.
\end{corollary}

It is straightforward to see that filtered colimits of exact sequences of such
$\i$-categories are exact. 

\begin{lemma}
Given a filtered diagram of exact sequences $\aA_\alpha\to \aB_\alpha\to
\aC_\alpha$ of compact idempotent-complete small stable $\i$-categories,
the colimit $\aA\to \aB\to \aC$ is an exact sequence of idempotent-complete
small stable $\i$-categories; that is, $\aA\to \aB$ is fully faithful with
cofiber $\aC$.
\end{lemma}

\begin{proof}
This follows from the fact that the filtered colimit of fully faithful
functors is a fully faithful functor and that the cofiber of a
filtered colimit of fully faithful functors is equivalent to the
filtered colimit of the cofibers.
\end{proof}

The $\i$-category of $\i$-categories equipped with a localization,
\[
\mathrm{Loc}(\Cat_\i)\subseteq\Fun(\Delta^1,\Cat_\i),
\]
is the subcategory of those functors $g:\aB\to\aC$ such that $g$ admits a right adjoint $j$ with $g\circ j\simeq\Id_\aC$, and maps those transformations which also commute with the adjoint.
We have obvious analogues $\mathrm{Loc}(\stabcat)$ and $\mathrm{Loc}(\idemstabcat)$, and in the stable setting a localization is part of the data of a split-exact sequence.
We write
\[
\mathrm{Split}(\stabcat)\subseteq\Fun(\Delta^2,\stabcat)
\]
for the subcategory consisting of those diagrams $\aA\overset{f}{\to}\aB\overset{g}{\to}\aC$ of small stable $\i$-categories such that $g\circ f\simeq 0$, $f$ is fully faithful with cofiber $g$, $f$ admits a right adjoint $i$ with $\Id_\aA\simeq i\circ f$, and $g$ admits a right adjoint $j$ with $g\circ j\simeq\Id_\aC$; maps are those transformations
\[
\xymatrix{
\aA\ar[r]^f\ar[d]^\alpha & \aB\ar[r]^g\ar[d]^\beta & \aC\ar[d]^\gamma\\
\aA'\ar[r]^{f'} & \aB'\ar[r]^{g'} & \aC'}
\]
which also commute with the adjoints, i.e., $\alpha\circ i\simeq i'\circ\beta$ and $\beta\circ j\simeq j'\circ\gamma$.

\begin{proposition}\label{prop:split=loc}
The functors $\mathrm{Split}(\stabcat)\to\mathrm{Loc}(\stabcat)$ and $\mathrm{Split}(\idemstabcat)\to\mathrm{Loc}(\idemstabcat)$, induced by the inclusion $\Delta^1\cong\Delta^{\{1,2\}}\to\Delta^2$, are equivalences.
\end{proposition}

\begin{proof}
First observe that a split-exact sequence $\aA\overset{f}{\to}\aB\overset{g}{\to}\aC$ is completely determined by the projection $g:\aB\to\aC$
together with its section $j:\aC\to\aB$.
This is because $f:\aA\to\aB$ is the fiber of $g$, which we may identify with the full subcategory of $\aB$ spanned by the $b\in\aB$ such that $g(b)\simeq 0$, and, since $f$ is fully faithful, $i:\aB\to\aA$ is determined by the composite $f\circ i:\aB\to\aA\to\aB$, the fiber
\[
f\circ i\to\id_{\aB}\to j\circ g
\]
of the unit map of the adjunction $(g,j)$.
Hence $\mathrm{Split}(\idemstabcat)\to\mathrm{Loc}(\idemstabcat)$
has contractible (homotopy) fibers and is therefore and equivalence.
\end{proof}

\begin{proposition}\label{prop:splitapprox}
The $\i$-category $\mathrm{Split}(\idemstabcat)$ of split-exact
sequences of small stable $\i$-categories is accessible.  In
particular, there exists a cardinal $\kappa$ such that any split-exact
sequence in $\idemstabcat$ is a $\kappa$-filtered (and hence filtered)
colimit of $\kappa$-compact split-exact sequences in $\idemstabcat$.
\end{proposition}

\begin{proof}
By proposition \ref{prop:split=loc}, we may equivalently show that
$\mathrm{Loc}(\idemstabcat)$ is accessible.  Recall that an adjunction
of $\i$-categories can be described as a map $\aM \to \Delta^1$ which
is both a cocartesian fibration and a cartesian
fibration~\cite[5.2.2.1]{HTT}.  This leads us to consider the
commutative diagram of pullback squares
\[
\xymatrix{
\mathrm{Loc}(\idemstabcat)\ar[r]\ar[d] & \mathrm{Loc}(\Cat_\i)\ar[r]\ar[d] & \Cat_{\i/\Delta^1}^\mathrm{cart,ff}\ar[d]\\
\Fun(\Delta^1,\idemstabcat)\ar[r] & \Cat_{\i/\Delta^1}^\mathrm{cocart}\ar[r] & \Cat_{\i/\Delta^1}}
\]
in which $\Cat_{\i/\Delta^1}^\mathrm{cocart}\subset\Cat_{\i/\Delta^1}$
(respectively,
$\Cat_{\i/\Delta^1}^\mathrm{cart,ff}\subset\Cat_{\i/\Delta^1}$) denote
the subcategories of cocartesian fibrations (respectively, cartesian
fibrations whose straightenings are fully faithful) and functors which
preserve cocartesian (respectively, cartesian) edges.

Since $\idemstabcat\subseteq\Cat_\i$ is an accessible functor between
accessible $\i$-categories, it suffices, using the~\cite[5.4.4.3,
5.4.5.16, 5.4.6.6]{HTT} and the duality between cartesian and
cocartesian fibrations, to show that
$\Cat_{\i/\Delta^1}^\mathrm{cart,ff}$ is accessible, and that the
inclusions 
$\Cat_{\i/\Delta^1}^\mathrm{cart,ff}\subseteq\Cat_{\i/\Delta^1}^\mathrm{cart}\subseteq\Cat_{\i/\Delta^1}$
are accessible functors. 
The straightening functor gives an equivalence
$\Cat_{\i/\Delta^1}^\mathrm{cart}\simeq\Pre_{\Cat_\i}(\Delta^1)$
between cartesian fibrations over $\Delta^1$ and presheaves of
$\i$-categories on $\Delta^1$~\cite[3.2.0.1]{HTT}.

In order to understand the condition of being fully faithful, we write
$\Cat_\i$ as an accessible localization
$\Cat_\i\subseteq\Pre(\mathrm{N}(\Delta))$ of simplicial
spaces~\cite{JT}.  A functor is fully faithful when the corresponding
map of (local) simplicial spaces is fully faithful, and recall that a
map of simplicial spaces $j\colon X\to Y$ is fully faithful if
and only if $\map(\Delta^1,X) \to \map(\partial\Delta^1,
X)\times_{\map(\partial\Delta^1,Y)} \map(\Delta^1, Y)$ is an 
equivalence.  It follows that $\Cat_{\i/\Delta^1}^\mathrm{cart,ff}$ is
the accessible localization of
$\Pre(\Delta^1\times\mathrm{N}(\Delta))$ obtained by also inverting
the pushout product of $\Id_{\Delta^1}$ and
$\partial\Delta^1\to\Delta^1$.  Thus
$\Cat_{\i/\Delta^1}^\mathrm{cart,ff}$ and
$\Cat_{\i/\Delta^1}^\mathrm{cart,ff}\subseteq\Cat_{\i/\Delta^1}^\mathrm{cart}$
are accessible.

Finally, it remains to show that the inclusion
$\Cat_{\i/\Delta^1}^\mathrm{cart}\subseteq\Cat_{\i/\Delta^1}$ is
accessible.  First, observe that finite limits commute with filtered colimits in $\Cat_\i$, as $\Cat_\i\simeq\Ind(\Cat_\i^\omega)$ is compactly generated, the inclusion $\Ind(\Cat_\i^\omega)\subseteq\Pre(\Cat_\i^\omega)$ preserves limits and filtered colimits \cite[5.3.5.3]{HTT}, and finite limits commute with filtered colimits in presheaf $\i$-categories (this uses \cite[5.3.3.3]{HTT} and the fact that (co)limits in presheaf $\i$-categories are computed objectwise). It follows that the filtered colimit $\C\simeq\colim_i\C_i$ of cartesian fibrations
$p_i\colon\C_i\to\Delta^1$, computed in $\Cat_\i$, is itself a
cartesian fibration $p\colon\C\to\Delta^1$; indeed, the inclusions
$\C_i\to\C$ preserve cartesian edges over $\Id_{\Delta^1}$, and
inspection of the fibers 
\[
\C\times_{\Delta^1}\Delta^0\simeq(\colim\C_i)\times_{\Delta^1}\Delta^0\simeq\colim(\C_i\times_{\Delta^1}\Delta^0)
\]
over each vertex $\Delta^0\to\Delta^1$ shows that $p\colon \C\to\Delta^1$ is also the colimit in $\Cat_{\i/\Delta^1}^\mathrm{cart}$. 
\end{proof}

\subsection{Strict-exact sequences}\label{sec:strict}

\begin{definition}\label{def:strictexact}
An exact sequence of small stable $\i$-categories of the form
\begin{equation}\label{eq:strict}
\cA \too \cB \too \cB/\cA
\end{equation}
is called {\em strict-exact} if $\aA\to\aB$ is the inclusion of a full subcategory and any object of $\cB$ which is a summand of an object of $\cA$ is also in $\cA$. In particular, every split-exact sequence (see definition~\ref{def:splitexact}) is equivalent to a strict-exact exact sequence.
\end{definition}
We denote by $\underline{\cE_\mathrm{wL}^\kappa}$ a set of representatives of
strict-exact sequences $\aA\to\cB\to\cB/\cA$ with $\cB$ in
$(\stabcat)^\kappa$. 
\begin{proposition}\label{prop:approx}
Any strict-exact sequence $\cA \to \cB \to \cB/\cA$ is a $\kappa$-filtered colimit of strict-exact sequences $\cA_{\alpha} \to \cB_{\alpha} \to \cB_{\alpha}/\cA_{\alpha}$ in $\underline{\cE_\mathrm{wL}^\kappa}$.
\end{proposition}
\begin{proof}
Write $\cB\simeq\colim_\alpha\cB_\alpha$ as a $\kappa$-filtered colimit of $\kappa$-compact stable $\i$-categories $\cB_\alpha$, and define $\aA_\alpha=\aA\times_\aB\aB_\alpha$ to be the full subcategory of $\aA$ consisting of those objects of $\aA$ which lie in the image of $\aB_\alpha$.
Evidently, $\aA\to\aB\to\aB/\aA$ is the $\kappa$-filtered colimit of the exact sequences $\aA_\alpha\to\aB_\alpha\to\aB_\alpha/\aA_\alpha$, and $\aA_\alpha\to\aB_\alpha\to\aB_\alpha/\aA_\alpha$ is strict-exact because if $Y\in\aB_\alpha$ is a summand of $X\in\aA_\alpha$ then $Y\in\aA_\alpha$ because the image of $Y$ in $\aB$ lies in $\aA$.
\end{proof}

We denote by $\cE^\kappa_\mathrm{L}$ a set of representatives of maps of the form $\cA \to \Idem(\cA)$ with $\cA$ in $(\stabcat)^\kappa$.
\begin{proposition}\label{prop:approx1}
Any map of the form $\cA \to \Idem(\cA)$ is a $\kappa$-filtered colimit of elements of $\cE^\kappa_\mathrm{L}$.
\end{proposition}
\begin{proof}
Write $\cA\simeq\colim_\alpha\aA_\alpha$ as a $\kappa$-filtered
colimit of $\kappa$-compact small stable $\i$-categories
$\cA_{\alpha}$.  Then
$\Idem(\cA)\simeq\colim_\alpha\Idem(\aA_\alpha)$, since $\Idem$ 
(viewed as an endofunctor of $\stabcat$) commutes with
$\kappa$-filtered colimits --- this follows from the characterization
of $\Idem$ in terms of a subcategory of the $\Ind$
category \cite[5.4.2.4]{HTT} and the fact that filtered colimits in
$\stabcat$ can be computed in $\icat$ \cite[1.1.4.6]{HAG}.
\end{proof}

\section{Additivity}\label{sec:additive}

In this section we construct the universal additive invariant of small
stable $\i$-categories; see theorem~\ref{thm:main1}.  Its construction
is divided in two steps\,: first, using
proposition~\ref{prop:splitapprox}, we construct the unstable version;
see theorem~\ref{thm:weakadd}. Then, by stabilizing, we obtain the
universal additive invariant.  Our arguments follow the pattern of the
analogous result for dg-categories given in \cite{Duke}.

\begin{definition}\label{def:additive}
Let $\aD$ be a stable presentable $\i$-category. A functor
$$ E: \stabcat \too \aD$$
is called an {\em additive invariant} of small stable $\i$-categories if it inverts Morita equivalences (see definition~\ref{defn:iMor}), preserves filtered colimits, and satisfies {\em additivity}, i.e., given a split-exact sequence
\begin{equation}\label{number3}
\xymatrix{
\cA \ar[r]_{i} & \cC \ar@<-1ex>[l]_f
\ar[r]_j & \cB \ar@<-1ex>[l]_{g} \,,
}
\end{equation}
the functors $i$ and $g$ induce an equivalence in $\cD$
$$ E(\cA) \vee E(\cB) \stackrel{\sim}{\too} E(\cC)\,.$$
We denote by $\Fun_{\mathrm{add}}(\stabcat, \cD)$ the $\i$-category of additive invariants with values in $\cD$.
\end{definition}

\begin{example}
As we discuss in Sections~\ref{sec:kthy} and ~\ref{sec:cyclo},
appropriate versions of algebraic $K$-theory and topological
Hochschild homology provide additive invariants of small stable
$\i$-categories.
\end{example}

\subsection{Unstable version}
Let us denote by $\Pre((\idemstabcat)^{\omega})_*$ the $\i$-category $$\Fun(((\idemstabcat)^{\omega})^{\op}, \aT_{\i})_*$$ of presheaves
of pointed spaces on the essentially small $\i$-category
$(\idemstabcat)^{\omega}$ of compact idempotent-complete small stable
$\i$-categories.

\begin{lemma}\label{lem:filtrant}
Let $\cD$ be a pointed presentable $\i$-category. Then, we have an
equivalence of $\i$-categories
$$ \Fun^{\L}(\Pre((\idemstabcat)^{\omega})_*, \aD)
\simeq \Fun_{\flt}(\idemstabcat, \aD)\,,$$
where the right-hand side denotes the $\i$-category of morphisms of
$\i$-categories which preserve filtered colimits.
\end{lemma}
\begin{proof}
The proof is a consequence of the equivalences
\begin{align*}
\Fun^{\L}(\Pre((\idemstabcat)^{\omega})_*, \aD)
&\htp \Fun((\idemstabcat)^{\omega}, \aD) \\
&\htp \Fun_{\flt}(\Ind((\idemstabcat)^{\omega}), \aD) \\
&\htp \Fun_{\flt}(\idemstabcat, \aD)\,,
\end{align*}
where the first follows from \cite[5.1.5.6]{HTT} and the fact that $\aD$ is pointed, and the last follows from corollary~\ref{cor:cgen}. 
\end{proof}

Let
\[
\phi \colon \idemstabcat \to \Pre((\idemstabcat)^{\omega})_*
\]
be the functor obtained by first taking the Yoneda embedding and then
restricting the presheaves to the category
$(\idemstabcat)^{\omega}$. Recall from corollary~\ref{cor:setrep} that
we can choose a fixed set $\aE$ of representatives of split-exact
sequences in
$(\idemstabcat)^{\omega}$. We denote by $\Motaddun$ the
localization of $\Pre((\idemstabcat)^{\omega})_*$
\cite[5.5.4.15]{HTT} with respect to the set of maps
\begin{equation}\label{eqn:splitlocs}
\phi (\aC)/\phi(\aA) \to \phi (\aB)\,,
\end{equation}
where $\aA \to \aC \to \aB$ is a split-exact sequence in $\aE$.
Finally, let $\Umotun$ be the composite
\begin{equation}\label{un}
\xymatrix{
\stabcat \ar[r]^-{\Idem(-)}
& \idemstabcat \ar[r]^-{\phi} & \Pre(\idemstabcat)^{\omega})_* \ar[r]^-{\gamma}
& \Motaddun \,,
}
\end{equation}
where $\gamma$ is the localization functor.
\begin{theorem}\label{thm:weakadd}
The functor $\Umotun$ inverts Morita equivalences, preserves filtered colimits, and sends split-exact sequences in $\stabcat$ to cofiber sequences in $\Motaddun$.
Moreover, $\Umotun$ is universal with respect to these properties,
i.e., given any pointed presentable $\i$-category $\aD$, we have an equivalence of $\i$-categories
\[
(\Umotun)^{\ast}: \Fun^{\L}(\Motaddun, \aD) \stackrel{\sim}{\too} \Fun_{\mathrm{add}}^{\mathrm{un}}(\stabcat, \aD)\,,
\]
where the right-hand denotes the full subcategory of
$\Fun(\stabcat, \aD)$ of morphisms of $\i$-categories which satisfy
the above conditions.
\end{theorem}

\begin{proof}
The result follows from definition~\ref{defn:iMor},
lemma~\ref{lem:filtrant} and from the universal property of Bousfield
localization (see section~\ref{sec:bousloc}: The functor $\phi$
preserves filtered colimits and by proposition~\ref{prop:splitapprox}
any split-exact sequence can be approximated by a filtered colimit of
split-exact sequences in $\aE$).
\end{proof}

\subsection{Universal additive invariant}
Let $\Motadd$ be the stabilization $\Stab(\Motaddun)$
\cite[\S 1.4]{HAG} of $\Motaddun$; by construction, this is a
stable $\i$-category. Denote by $\Umot$ the following composite
$$ \stabcat \stackrel{\Umotun}{\too} \Motaddun \too \Stab(\Motaddun)\,.$$
\begin{remark}
Note that we have the equivalences
\begin{align*}
\Stab(\Pre((\idemstabcat)^{\omega})_*)
&= \Stab(\Fun(((\idemstabcat)^{\omega})^{\op}, \aT_{\i *})) \\
&\htp \Fun(((\idemstabcat)^{\omega})^{\op}, \Stab(\aT_{\i *}))) \\
&\htp \Fun(((\idemstabcat)^{\omega})^{\op}, \ispec),
\end{align*}
the last of which follows from~ \cite[1.4.4.11]{HAG}. Therefore,
defining
\[
\Pre_{\ispec}((\idemstabcat)^{\omega}) =
\Fun(((\idemstabcat)^{\omega})^{\op}, \ispec)
\] 
and writing
\[
\psi \colon \idemstabcat \too \Pre_{\ispec}((\idemstabcat)^{\omega})
\]
for the natural functor, we see that $\Motadd$ can alternately be
described as the localization of $\Pre_{\ispec}((\idemstabcat)^{\omega})$ with respect to the set of maps
\begin{equation}\label{eq:loc}
\psi(\aC)/\psi(\aA) \to \psi(\aB)\,,
\end{equation}
where $\cA \to \cC \to \cB$ is a split-exact sequence in $\cE$.
\end{remark}
\begin{theorem}\label{thm:main1}
The functor $\Umot$ is the {\em universal additive invariant}, i.e.,
given any stable presentable $\i$-category $\aD$, we have an
equivalence of $\i$-categories
\[
(\Umot)^{\ast} \colon \Fun^{\L}(\Motadd, \aD) \stackrel{\sim}{\too}
\Fun_{\mathrm{add}}(\stabcat, \aD)\,.
\]
\end{theorem}

\begin{proof}
The result follows from theorem~\ref{thm:weakadd} and from the
universal property of stabilization (i.e., \cite[1.4.5.5]{HAG}). Note
that stabilization preserves colimits and $\Umotun$ sends split-exact sequences to cofiber sequences, so the split-exact
sequence (\ref{number3}) is sent to a split cofiber sequence $\Umot(\cC)\simeq \Umot(\cA)\vee\Umot(\cB)$
in $\Motadd$.
\end{proof}

\section{Connective $K$-theory}\label{sec:kthy}

In this section, we verify that higher algebraic $K$-theory provides
an additiive invariant of small stable $\i$-categories;
theorem~\ref{thm:main1} then applies to show that ths invariant
descends to $\Motadd$.  Furthermore, following the outline
of \cite{Duke}, we prove the essential result that algebraic
$K$-theory in fact becomes co-representable in $\Motadd$ (see
theorem~\ref{thm:main2}).  The underlying point is that Waldhausen's
$\Sdot$ construction simply becomes the suspension in $\Motadd$.  This
result will allow us to understand transformations between additive
theories from algebraic $K$-theory via the Yoneda lemma; we use this
in Section~\ref{sec:cyclo} to characterize the cyclotomic trace map.
We begin by developing the necessary background on the construction of
algebraic $K$-theory for small $\i$-categories with finite colimits,
and in subsection~\ref{sub:comparison} we compare the $K$-theory of a suitable Waldhausen category with the $K$-theory of its underlying $\i$-category.

\subsection{Algebraic $K$-theory of $\i$-categories}\label{sec:back-K}

Waldhausen's algebraic $K$-theory functor takes as input a category
with cofibrations and weak equivalences.  It is now well understood
that, under mild hypotheses, the $K$-theory spectrum is determined by the Dwyer-Kan
localization $L^H\aC$ of the Waldhausen category
$\aC$ \cite{ToenVezzosi, BM4, Cisinski}.  Since
$\N((L^H \aC)^{\textrm{fib}})$ yields the $\i$-category associated to
$\aC$, these results can be interpreted as saying that the algebraic
$K$-theory of a Waldhausen category is an invariant of the underlying
$\i$-category.
Moreover, it has long been folklore that given a sufficiently good
theory of $\i$-categories one can define analogues of Waldhausen's
construction of algebraic $K$-theory (e.g., see \cite[\S
7]{ToenVezzosi} for a sketch of such a definition in the context of
Segal categories).  In this subsection we study a version of such a
direct construction of the algebraic $K$-theory of $\i$-categories in
the setting of quasicategories~\cite[1.2.2.5]{HAG}.  We prove that
Waldhausen's algebraic $K$-theory of a Waldhausen category $\aC$ is
equivalent as a spectrum to this $\i$-categorical algebraic $K$-theory
of the associated $\i$-category $\N((L^H \aC)^{\mathrm{fib}})$.

We begin by reviewing Waldhausen's $\Sdot$ construction.  Let $\aC$ be
a Waldhausen category.  Let $\Ar[n]$ denote the category of arrows in
$[n]$:  $\Ar[n]$ has objects $(i,j)$ for $0\leq i\leq j\leq n$ and a
unique map $(i,j)\to (i',j')$ for $i\leq i'$ and $j\leq j'$.  
Then $\Sdot[n]\aC$ is the full subcategory of the
category of functors $A\colon \Ar[n]\to \aC$ such that:
\begin{itemize}
\item $A_{i,i}=*$ for all $i$, 
\item The map $A_{i,j}\to A_{i,k}$ is a w-cofibration for all $i \leq
j \leq k$, and
\item The diagram
\[  \xymatrix@-1pc{%
A_{i,j}\ar[r]\ar[d]&A_{i,k}\ar[d]\\A_{j,j}\ar[r]&A_{j,k}
} \]
is a pushout square for all $i \leq j \leq k$, 
\end{itemize}

The algebraic $K$-theory space of $\aC$ is then defined to be $\Omega
|w\subdot \Sdot \aC|$, where the weak equivalences in $\Sdot \aC$ are
defined pointwise.  Furthermore, since each $\Sdot[n] \aC$ is itself a
Waldhausen category (with the Reedy cofibrations), we can iterate the
$\Sdot$ construction.  The algebraic $K$-theory spectrum of $\aC$ is
the spectrum with $n$th space $|w\subdot \Sdot^{(n)} \aC|$.

Now let $\aC$ be a small pointed $\i$-category with finite colimits.
The following definition \cite[1.2.2.2]{HAG} is the $\i$-categorical
analogue of Waldhausen's $\Sdot$ construction.

\begin{definition}\label{defn:gap}
Denote by $\Gap([n],\aC)$ the full subcategory of
$\Fun(\N(\Ar[n]),\aC)$ spanned by the functors 
$\N(\Ar[n])\to\aC$ such that, for each $i\in I$, $F(i,i)$ is a zero
object of $\aC$, and for each $i<j<k$, the square
\[
\xymatrix{
F(i,j)\ar[r]\ar[d] & F(i,k)\ar[d]\\
F(j,j)\ar[r]       & F(j,k)}
\]
is cocartesian.
\end{definition}

\begin{remark}
There is an obvious generalization of this definition to small pointed
$\i$-categories equipped with a suitable subcategory of
``cofibrations'' (satisfying the usual axioms, e.g. that cofibrations
are stable under cobase change).  However, in the presence of
factorization hypotheses, this does not yield added generality; e.g.,
see~\cite[1.3]{BM4}, which under such assumptions describes the
$K$-theory space in terms of the Dwyer-Kan localization $\aC$ regarded
as a category with weak equivalences.
\end{remark}

As with the classical $\Sdot$ construction, when $\aC$ has all
colimits, the data of the cocartesian squares (i.e., cofibers for the
maps $F(i,j) \to F(i,k)$) is necessary only for the simplicial
structure.

\begin{lemma}\label{lem:dropcof}
Let $\aC$ be an $\i$-category with finite colimits.  Then for each
$n$, the forgetful functor 
\[
\Gap([n], \aC) \to \Fun(\Delta^{1,2,\ldots,n}, \aC)
\]
is an equivalence of $\i$-categories (and observe that
$\Delta^{1,2,\ldots,n} \htp \N([n-1])$).
\end{lemma}

\begin{proof}
This follows from the fact that the space of colimits for a given
diagram in an $\i$-category is contractible \cite[1.2.12.9,
1.2.13.5]{HTT}.  Alternatively, a constructive proof along the lines
of \cite[2.9]{BM3} (using a mapping cylinder argument) can be given using
the comparison discussed in section~\ref{sub:comparison} below.
\end{proof}

\begin{remark}
Lemma~\ref{lem:dropcof} implies that $\Gap([n],\aC)$ is
stable when $\aC$ is stable.
\end{remark}

Following \cite[1.2.2.5]{HAG}, we define a simplicial $\i$-category
$\iSdot \aC$ by the rule $\iSdot[n] \C=\Gap([n],\aC)$.  Applying
passage to the largest Kan complex levelwise, we obtain a simplicial
space $\lkan{\iSdot \aC}$.  Then $\Omega |\lkan{\iSdot \aC}|$ is the
$\i$-categorical version of Waldhausen's $K$-theory space.
Furthermore, for each $n$, $\Gap([n],\aC)$ is itself a small pointed
$\i$-category with finite colimits: once again, we can can iterate
this procedure.  Since $\Gap([0],\aC)$ is contractible (with preferred
basepoint given by the point in $\aC$) and $\Gap([1],\aC)$ is
equivalent to $\aC$, there is a natural map 
\[
S^1 \sma \lkan{\aC} \to |\lkan{\iSdot \aC}|
\]
given by the inclusion into the $1$-skeleton.  Therefore, the spaces
$|\lkan{(\iSdot)^{n}(\aC)}|$ assemble to form a spectrum $K(\aC)$;
this is the $\i$-categorical version of Waldhausen's $K$-theory
spectrum.  We can see from the definition that $\Gap([n],\aC)$ is
natural in (right) exact functors, and therefore $\iSdot \aC$ and
$|\lkan{\iSdot \aC}|$ are also natural.  Since the equivalence
$\Gap([1],\aC) \to \Fun(*,\aC) \htp \aC$ induced by the restriction
map is natural in $\aC$, we deduce that the $K$-theory spectrum is
natural in exact functors.

In practice, we find it more convenient to use an ``all at once''
reformulation of the definition of the iterated $\Sdot$ construction
(e.g., see \cite[A.5.4]{BM4}, \cite[2.2]{BM5}, the appendix
to \cite{GeisserHesselholt}, and also \cite[\S 2]{Rognes}).

\begin{definition}\label{consitsdot}
Write $\Ar_{n_{1},\dotsc,n_{q}}$ for
$\Ar[n_{1}]\times \dotsb \times \Ar[n_{q}]$.  For a functor
\[
A\colon \N(\Ar_{n_{1},\dotsc,n_{q}})=\N(\Ar[n_{1}] \times \dotsb \times
\Ar[n_{q}])\to \aC,
\]
we write $A_{i_{1},j_{1};\dotsc;i_{q},j_{q}}$ for the value of $A$ on
the object $((i_{1},j_{1}),\dotsc,(i_{q},j_{q}))$.
Let $\Gap(([n_1],\dotsc,[n_q]), \aC)$ be the full subcategory of
$\Fun(\N(\Ar_{n_1,\dotsc,n_q}), \aC)$ spanned by the functors such that
\begin{itemize}
\item $A_{i_{1},j_{1};\dotsc;i_{q},j_{q}} \htp *$ whenever
$i_{k}=j_{k}$ for some $k$.
\item For every object
$(i_{1},j_{1};\dotsc;i_{q},j_{q})$ in $\Ar[n_{1}]\times \dotsb \times
\Ar[n_{q}]$, every $1\leq r\leq q$, and every $j_{r}\leq k\leq n_{r}$, the square
\[  \xymatrix@-1pc{%
A_{i_{1},j_{1};\dotsc;i_{q},j_{q}}\ar[r]\ar[d]
&A_{i_{1},j_{1};\dotsc;i_{r},k;\dotsc;i_{q},j_{q}}\ar[d]\\
A_{i_{1},j_{1};\dotsc;j_{r},j_{r};\dotsc;i_{q},i_{q}}\ar[r]
&A_{i_{1},j_{1};\dotsc;j_{r},k;\dotsc;i_{q},i_{q}}
} \]
is a cocartesian square.
\end{itemize}
\end{definition}

Now we define the multisimplicial $\i$-category
\[
\iSdotq[n_{1},\dotsc,n_{q}]\aC = \Gap(([n_1],\dotsc,[n_q]), \aC).\]
We regard $(\iSdot)^{(0)}$ as $\aC$ and it is clear that
$(\iSdot)^{(1)}_n$ is $\Gap([n],\aC)$.  Now we directly the define the
$K$-theory spectrum of an $\i$-category $\aC$ with finite colimits to
be the spectrum with $q$-th space
\[
K\aC(q)= \lkan{\iSdotq},
\]
The suspension maps $\Sigma K\aC(q)\to K(q+1)$ are induced on diagrams
by the projection map
\[
\Ar[n_{1}]\times \dotsb \times \Ar[n_{q}]\times \Ar[n_{q+1}]\to
\Ar[n_{1}]\times \dotsb \times \Ar[n_{q}].
\]

From definition~\ref{consitsdot} it is now clear that the construction
of the $K$-theory spectrum is functorial in (right) exact functors. 

\subsection{Comparison with Waldhausen's $K$-theory}\label{sub:comparison}

We now establish a comparison between Waldhausen's algebraic
$K$-theory of a Waldhausen category $\aC$ and the $\i$-categorical
version of the algebraic $K$-theory of the associated simplicial
category $L^H \aC$.  The comparison is essentially a consequence of
the theory of rigidification of homotopy coherent diagrams to strict
diagrams in a model category (originally studied by
Dwyer-Kan \cite{DKdiagrams}), which allows us to pass between
$\i$-categorical diagrams and point-set diagrams, and the
``homotopical'' $\Spdot$ construction of \cite{BM3}, which allows us
to replace the use of pushouts by homotopy pushouts for suitable
Waldhausen categories.
The version of the comparison of homotopy coherent diagrams to strict
diagrams we use is originally due to
Hirschowitz-Simpson \cite{HirschowitzSimpson} in the context of Segal
categories (see also Rezk's work in Segal spaces~\cite[8.12]{Rezk}).
Since we are using quasicategories in this paper, we work with the
version proved by Lurie in that setting \cite[4.2.4.4]{HTT}.

Let $S$ be a small simplicial set, $\aD$ a small simplicial category, and
$u \colon \mathfrak{C}[S] \to \aD$ an equivalence.  Let $\aA$ be a
combinatorial simplicial model category, and let $\aU$ be a
$\aD$-chunk of $\aA$ (see \cite[A.3.4.9]{HTT} for a discussion of
$\aD$-chunks).  Then the induced map
\[
\N((\aU^{\aD})^{\cf}) \to \Fun(S, \N(\aU^{\cf}))
\]
is a categorical equivalence of simplicial sets.  Here the notation
$(\aU^{\aD})^{\cf}$ indicates the full subcategory of $\aA^{\aD}$
consisting of cofibrant-fibrant objects (in the projective model
structure) landing in $\aU$.

Specializing to our situation, assume that $S$ is the (ordinary) nerve
$\N(J)$ of a diagram (small category) $J$; that is, $J$ is regarded as
a discrete simplicial category.  Then the counit map
$\mathfrak{C}[\N(J)] \to J$ is an equivalence and so we have that the
induced map
\[
\N((\aU^J)^{\cf}) \to \Fun(\N(J), \N(\aU^{\cf}))
\]
is a categorical equivalence of simplicial sets.

\begin{lemma}\label{lem:diagcomp}
Let $\aA$ be a combinatorial simplicial model category and $\aC
\subset \aA$ a full subcategory.  Then for each $n$ the induced map
\[
\N((\aC^{\Ar[n]})^{\cf}) \to \Fun(\N(\Ar[n]), \N(\aC^{\cf}))
\]
is a categorical equivalence of simplicial sets.
\end{lemma}

\begin{proof}
By \cite[A.3.4.15]{HTT}, we can choose a small subcategory $\aV
\subset \aA$ which contains $\aC$ and such that $\aV$ is an
$(\Ar[n])$-chunk for each $n$ and moreover $\N((\aC)^{\cf})$ is
equivalent to $\N((\aV)^{\cf})$.  Then as discussed above,
\cite[4.2.4.4]{HTT} implies that for each $n$ the natural map
\[
\N(((\aV)^{\Ar[n]})^{\cf}) \to \Fun(\N(\Ar[n]),
\N((\aV)^{\cf}))
\]
is a categorical equivalence of simplicial sets.
\end{proof}

To apply these rigidification results, we use the $\Spdot$
construction.  The $\Spdot$ construction \cite[2.7]{BM3} is a variant
of Waldhausen's $\Sdot$ construction defined by replacing the cocartesian
squares in the definition of $\Sdot$ with homotopy cocartesian
squares.  In order to define the $\Spdot$ construction, we must work
with Waldhausen categories for which there is a reasonable notion of
homotopy cocartesian squares.  We briefly recall this theory
from \cite[\S 2]{BM3}.  A map is a {\em weak cofibration} if it is
equivalent by a zig-zag to a cofibration, and a square is a homotopy
cocartesian square if it equivalent by a zig-zag to a pushout square
with one leg a cofibration.  For control on these notions, we require
the hypothesis that any map in $\aC$ can be factored (not necessarily functorially) as a
cofibration followed by a weak equivalence.

For such a Waldhausen category $\aC$, we can then define
$\Spdot[n] \aC$ to be the full subcategory of the category of functors
$\Ar[n] \to \aC$ such that:
\begin{itemize}
\item $A_{i,i} \htp *$ for all $i$, 
\item The map $A_{i,j}\to A_{i,k}$ is a weak cofibration for all $i \leq
j \leq k$, and
\item The diagram
\[  \xymatrix@-1pc{%
A_{i,j}\ar[r]\ar[d]&A_{i,k}\ar[d]\\A_{j,j}\ar[r]&A_{j,k}
} \]
is a homotopy cocartesian square for all $i \leq j \leq k$, 
\end{itemize}

By construction, the $\Spdot$ construction is
functorial in {\em weakly exact functors}, i.e., functors that
preserve weak equivalences and homotopy cocartesian squares.
Moreover, the natural inclusion 
\[
w\subdot \Sdot \aC \to w\subdot \Spdot \aC
\]
is a weak equivalence \cite[2.9]{BM3}.  Therefore, we can equivalently
define the algebraic $K$-theory space of a Waldhausen category $\aC$
as $\Omega |w\subdot \Spdot \aC|$ and similarly the algebraic
$K$-theory spectrum of $\aC$ as having $n$th space
$|w\subdot (\Spdot)^{(n)} \aC|$.

Now, let $\aC$ be a Waldhausen category that arises as a subcategory
of a model category.  Since a square is homotopy cocartesian in $\aC$
if and only if it is a pushout square in $\N((\aC)^{\cf})$, in this
setting the equivalence of Lemma~\ref{lem:diagcomp} restricts to give
an equivalence
\[
\N((\Spdot[n] \aC)^{\cf}) \to \Gap([n], \N((\aC)^{\cf})).
\]
Similar considerations for the iterated $\Spdot$
construction \cite[A.5.4]{BM4} (as modeled in
Definition~\ref{consitsdot}) yield the equivalence
\[
\N((\Spdotq[n_{1},\dotsc,n_{q}] \aC)^{\cf}) \to \Gap(([n_1], \dotsc,
[n_q]), \aC).
\]
Applying proposition~\ref{prop:nerve}, we then obtain the following
comparison of algebraic $K$-theory spaces and spectra.

\begin{corollary}\label{cor:Kcomp}
Let $\aA$ be a simplicial model category and $\aC \subset \aA$ a
small full subcategory which has all finite homotopy colimits.  Then
for each $n$ there is a weak equivalence of simplicial sets
\[
|w \subdot \Spdot[n] \aC| \htp |\lkan{\Sdot[n]^{\i}\N((\aC)^{\cf})}|.
\]
and for each $(n_1, \dotsc, n_q)$ there is a weak equivalence of
simplicial sets
\[
|w \subdot \Spdotq[n_{1},\dotsc,n_{q}] \aC| \htp
 |\lkan{\iSdotq[n_{1},\dotsc,n_{q}]\N((\aC)^{\cf})}|.
\]
\end{corollary}

In particular, this yields the following theorem:

\begin{theorem}\label{thm:kconsistency}
Let $\aA$ be a simplicial model category and $\aC \subset \aA$ a small
full subcategory of the cofibrants which admits all homotopy pushouts
and is a Waldhausen category via the model structure on $\aA$.  Then
there is an equivalence of spectra
\[
K(\aC) \htp K(\N((\aC)^{\cf}))
\]
which is natural in weakly exact functors.
\end{theorem}

Finally, specializing to our case, we find the following result.

\begin{corollary}\label{cor:kconsistency}
Let $\aC$ be a small pretriangulated spectral category and let
$\aM_{\aC}$ denote the category of perfect $\aC$-modules with its
Waldhausen structure induced by the model structure on $\aC$-modules.
Then there is an isomorphism in the stable category
 \[
K(\aM_{\aC}) \htp K(\Psi_{\perf}\aC).
\]
\end{corollary}

As a consequence, Waldhausen's additivity theorem applies to prove the
following proposition.

\begin{proposition}
The algebraic $K$-theory functor
\[
K:\idemstabcat\too\ispec
\]
is an additive invariant. 
\end{proposition}

\begin{proof} 
It suffices to show that $K$ preserves filtered colimits and split-exact sequences.
The former follows from the fact that the $\iSdot$ construction and restriction to the maximal subgroup preserve filtered colimits, as $\N(\Ar([n]))$ and $\Delta^0$ are compact $\i$-categories.  Corollary~\ref{cor:kconsistency} allows us to reduce to
consideration of split-exact sequences of spectral categories
\[
\widehat{\aA}_{\perf} \to \widehat{\aC}_{\perf} \to \widehat{\aB}_{\perf}.
\]
As in \cite{Matrix}, we observe that this sequence is Morita
equivalent to the sequence
\[
\widehat{\aA}_{\perf} \to E(\widehat{\aA}_{\perf}, \widehat{\aC}_{\perf}, \widehat{\aB}_{\perf}) \to \widehat{\aB}_{\perf}
\]
(where $E$ denotes Waldhausen's category of cofiber sequences in $\aC$
with first term in the image of $\aA$ and cofiber in the image of
$\aB$).  Now Waldhausen's additivity theorem implies the desired
splitting on $K$-theory.
\end{proof}

So far, all of our comparison results assume that the Waldhausen
category we are working with arises as a subcategory of a model
category.  In fact, we can extend our comparison and functoriality
results to Waldhausen categories $\aC$ such that all maps admit
(not necessarily functorial) factorizations as cofibrations followed
by weak equivalences and which are {\em
DKHS-saturated} (i.e., such that a map $f$ is a weak equivalence in
$\aC$ if and only if its image in the homotopy category is an
isomorphism).  We do this as follows, using a construction due to
Cisinski \cite[\S 4]{Cisinski}.

\begin{lemma}\label{lem:cis}
Let $\aC$ be a Waldhausen category with factorization and weak
equivalences that are DKHS-saturated.  Then there exists a Waldhausen
category $\aM(\aC)$ and a DK-equivalence $\aC \to \aM(\aC)$ which is
natural in weakly exact functors.
\end{lemma}

\begin{proof}
Given a Waldhausen category $\aC$, let $\aP(\aC)$ denote here the
pointed simplicial presheaves on $\aC$ with the projective model
structure (i.e., weak equivalences and fibrations are determined
pointwise).  We can successively localize $\aP(\aC)$ to produce a
category of presheaves which are pointwise Kan complexes, preserve
weak equivalences, and take homotopy cocartesian squares in $\aC$ to
homotopy pullback squares in $\aP(\aC)$; denote this category by
$\aP_{\ex}(\aC)$~\cite[4.10]{Cisinski}.  Let $\aM(\aC)$ denote the
full subcategory of the localized category consisting of the objects
which are cofibrant and weakly equivalent to representable presheaves;
this can be regarded as a Waldhausen category, inheriting structure
from the model structure on $\aP_{\ex}(\aC)$.  The Yoneda embedding
induces a DK-equivalence $\aC \to \aM(\aC)$~\cite[4.11]{Cisinski}, and
a weakly exact functor $\aC \to \aC'$ induces a left Quillen functor
$\aP_{\ex}(\aC) \to \aP_{\ex}(\aC')$ by left Kan extension, and hence
an exact functor $\aM(\aC) \to \aM(\aC')$ by restriction.
\end{proof}

As a corollary, we have the following result comparing of the
$K$-theory of Waldhausen categories that are DHKS-saturated and admit
factorization to the associated $K$-theory of $\i$-categories.

\begin{corollary}\label{cor:gencomp}
In the setting of~\ref{lem:cis}, there are equivalences
\[
K(\aC) \to K(\aM(\aC)) \to K(\N((\aM(\aC))^{\cf}))
\]
which are natural in weakly exact functors.
\end{corollary}

\begin{proof}
First, since we have a natural DK-equivalence $\aC \to \aM(\aC)$,
there is a natural equivalence $K(\aC) \to K(\aM(\aC))$ \cite{BM4,
Cisinski, ToenVezzosi}.  Next, since the category $\aM(\aC)$ satisfies
the hypothesis of Theorem~\ref{thm:kconsistency}, the second
equivalence holds.
\end{proof}

Given a (homotopically) pointed simplicial category with finite
homotopy colimits, we can use essentially the same construction to
produce a DK-equivalent Waldhausen category; see~\cite[\S
5]{ToenVezzosi} and~\cite[\S 14]{BM} for versions of such a
construction.

\subsection{Co-representability}

This subsection is entirely devoted to the proof of
theorem~\ref{thm:main2}.  The proof will follow from
propositions~\ref{prop:suspensions} and \ref{prop:local}.

\begin{theorem}\label{thm:main2}
Let $\aA$ be a small stable $\i$-category and $\aB$ be a compact
idempotent-complete small stable $\i$-category. Then there is a natural
equivalence of spectra
\[\Map(\Umot(\aB),\Umot(\aA)) \simeq K(\idemfun(\aB, \Idem(\aA)))\,.\]
When $\cB$ is the small stable $\i$-category
  $\ispec^{\omega}$ of compact spectra,
there is a natural equivalence of spectra
\[\Map(\Umot(\ispec^{\omega}),\Umot(\aA)) \simeq K(\Idem(\aA))\,.\]
In particular, we have isomorphisms of abelian groups
\[\Hom(\Umot(\ispec^{\omega})), \Sigma^{-n} \Umot(\cA)) \simeq K_n(\Idem(\cA))\]
in the triangulated category $\Ho(\Motadd)$.
\end{theorem}

\begin{notation}\label{not:prsheaves}
Given a small stable $\i$-category $\cA$, we denote by $K^{w}_{\cA}$
the object 
\[
\aB \mapsto | ( \iSdot( \idemfun(\cB, \Idem(\cA))))_{\mathrm{iso}}|
\]
in $\Pre((\idemstabcat)^{\omega})_*$ and by $K_{\aA}$
the object 
\[
\aB \mapsto K(\idemfun(\aB, \Idem(\aA)))
\]
in $\Pre_{\ispec}((\idemstabcat)^{\omega})$.  Note that the value of $K_{\cA}$
at $\ispec^{\omega}$ is precisely the $K$-theory
spectrum $K(\cA)$ of $\cA$, similarly and that $K^{w}_{\cA}$ is the delooping of the $K$-theory space.
\end{notation}

\begin{remark}
Recall that corollary~\ref{cor:rep} allow us to model the small
$\i$-category of exact functors $\idemfun(\aB, \Idem(\aA))$ as the
pretriangulated spectral category $\rep(\aB, \aA)$ of right-compact
$\spect(\aA)^{\op} \sma \spect(\aB)$-modules. Combined with
proposition~\ref{prop:nerve}, this implies that the associated mapping
space $\lkan{\idemfun(\aB,\Idem(\aA))}$ can be calculated as
$|w_\bullet\rep(\aB,\Idem(\aA))|$.  Moreover, $\rep(\aB,\Idem(\aA))$
inherits a natural Waldhausen structure as a full subcategory of the
cofibrant objects in the model structure on the category of
$\cB\text{-}\Idem(\aA)$-bimodules.  As such, we can also consider the
algebraic $K$-theory space $|w_\bullet\Sdot\rep(\aB,\Idem(\aA))|$ and
associated spectrum.
\end{remark}

In the following results, we will use the observation that Waldhausen's
$\Sdot$ construction, applied to a spectral category which is a
Waldhausen category with the cofibrations inherited from a spectral
model structure with all objects fibrant, produces a spectral category
(where the mapping spectra are given by an appropriate end)~\cite[\S
3]{BM2}.  To ensure we are in this setting, we will tacitly use the
equivalent model of spectral categories enriched in EKMM $S$-modules,
as explained in \cite[\S 15]{BM}.  Alternatively, we could stay with
spectral categories in symmetric spectra and use the ``Moore'' $\Sdot$
construction from~\cite[\S 4]{BM2}, which uses an explicit model of
the homotopy end.  We also need the following lemma which allows us to
bring the $\Sdot$ construction inside:

\begin{lemma}\label{lem:scomm}
Let $\aA$ and $\aB$ be small stable $\i$-categories.  Then we 
have an equivalence of simplicial $\i$-categories
\[
\iSdot \Fun^{\ex}(\aB, \aA) \htp \Fun^{\ex}(\aB, \iSdot \aA)
\]
and correspondingly an equivalence of spaces
\[
|\lkan{\iSdot \Fun^{\ex}(\aB, \aA)}| \htp |\lkan{\Fun^{\ex}(\aB, \iSdot \aA)}|.
\]
\end{lemma}

\begin{proof}
First, we show that for each $n$ there is an equivalence of
$\i$-categories 
\[
\Gap([n], \Fun^{\ex}(\aB,\aA)) \htp \Fun^{\ex}(\aB, \Gap([n], \aA)). 
\]
Since $\Fun(-,-)$ is defined simply as the mapping simplicial
set \cite[1.2.7.2]{HTT}, we have the equivalence
\[
\Fun(\N(\Ar[n]), \Fun(\aB, \aA)) \htp \Fun(\aB, \Fun(\N(\Ar[n]), \aA)).
\]
Since colimits in functor $\i$-categories are computed
pointwise~\cite[\S 5.1.2.3]{HTT} and the $\i$-category
$\Fun^{\ex}(\aB, \aA)$ is the full subcategory of $\Fun(\aB,\aA)$ 
spanned by the exact functors, we have a map
\[
\Gap([n], \Fun^{\ex}(\aB, \aA)) \to \Fun^{\ex}(\aB, \Gap([n], \aA)),
\]
and lemma~\ref{lem:dropcof} implies that it is an equivalence.  It
is now straightforward to check that these comparison maps assemble
into the desired simplicial equivalence. 
\end{proof}

We can now relate $\Motadd$ to the algebraic $K$-theory presheaf.

\begin{proposition}\label{prop:suspensions}
Let $\cA$ be a small stable $\i$-category. Then, we have a natural
equivalence $ \Sigma( \Umotun(\aA)) \simeq K^{w}_\aA$ in
$\Motaddun$ (see notation \ref{un}) and a natural equivalence $\Sigma \Umot(\aA) \simeq\Sigma K_{\aA}$ in $\Motadd$.
\end{proposition}

\begin{proof}
We begin by handling the unstable case.  Theorem~\ref{thm:morita}
implies that we can model $\cA$ by a small spectral category (which we
still denote by $\cA$).
Following~\cite[3.3]{McCarthy}, we consider the following
sequence of simplicial spectral categories
\[
\cA\subdot \stackrel{I}{\too} P\Sdot\cA \stackrel{Q}{\too} \Sdot\cA \,,
\]
where $\cA\subdot$ is a constant simplicial object and $P\Sdot\cA$ is
the simplicial path object of $\Sdot \cA$. By applying the functor
$\Umotun$ to this sequence, we obtain an induced morphism
\[
\Theta\colon\Umotun(P\iSdot\cA)/\Umotun(\cA\subdot) \too
\Umotun(\iSdot\cA)
\]
of simplicial objects in $\Motaddun$. We now show that each component $\Theta_n$ of $\Theta$ is an equivalence. For each $n \geq 0$, we have a split-exact sequence
\[
\xymatrix{
\cA \ar[r]_-{I_n} & PS_n\cA=S_{n+1}\cA \ar@<-1ex>[l]_-{R_n}
\ar[r]_-{Q_n} & S_n\cA \ar@<-1ex>[l]_-{S_n} \,,
}
\]
in which
\begin{align*}
I_n&(A)=(\ast \to A \stackrel{\mathrm{Id}}{\to} A \stackrel{\mathrm{Id}}{\to}
\cdots \stackrel{\mathrm{Id}}{\to} A),\\
Q_n&(\ast \to A_0 \to A_1 \to \cdots \to A_n)=(A_1/A_0 \to \cdots \to A_n/A_0),\\
S_n&(\ast \to A_0 \to A_1 \to \cdots \to A_{n-1})=(\ast \to \ast \to A_0 \to \cdots \to A_{n-1}),\\
R_n&(\ast \to A_0 \to A_1 \to \cdots \too A_{n-1})=A_0.
\end{align*}
By the construction of $\Motaddun$ (and of $\Umotun$), we conclude
that the induced morphisms 
\[
\Theta_n\colon  \Umotun(P\iSdot[n] \cA)/\Umotun(\cA) \too
\Umotun(\iSdot[n] \cA) \qquad n \geq 0\,,
\]
are equivalences in $\Motaddun$.  This allow us to obtain the
following cocartesian square
\[
\xymatrix{
\Umotun(\cA) \simeq |\Umotun(\cA)| \ar[d] \ar[r]
& |\Umotun(P\iSdot\cA)| \simeq \ast \ar[d] \\
\ast \ar[r] & |\Umotun(\iSdot \cA)|
}
\]
and so a natural equivalence
\[
\Sigma(\Umotun(\cA)) \stackrel{\sim}{\too} |\Umotun(\iSdot \cA)|
\]
in $\Motaddun$.  By combining this equivalence with the equivalences
\begin{eqnarray}
\Umotun(\iSdot \cA) & = &
|(\idemfun(-, \Idem(\iSdot \cA)))_{iso}| \nonumber \\ 
& \simeq & |(\iSdot \idemfun(-, \Idem(\cA)))_{iso}| \label{eq:new} \\
& =& K^w(\cA)\,, \nonumber
\end{eqnarray}
where \eqref{eq:new} follows from lemma~\ref{lem:scomm}, we conclude
that $ \Sigma( \Umotun(\aA)) \simeq K^{w}_\aA$ in $\Motaddun$.
The identification in the stable setting follows from the unstable considerations and the usual passage from results on the $K$-theory space to the $K$-theory spectrum.
\end{proof}

\begin{proposition}\label{prop:local}
Let $\cA$ be a small stable $\i$-category. Then, the presheaves
$K^{w}_{\aA}$ and $K_{\aA}$ (see notation~\ref{not:prsheaves}) are
local, i.e., given any split-exact sequnce $\aB \to \aC \to \aD$ in
$\aE$, the induced maps of spectra (see (\ref{eqn:splitlocs}) and
(\ref{eq:loc}))
$$ \map(\phi(\cD), K^w_{\aA}) \stackrel{\sim}{\too} \Map(\phi(\cC)/\phi(\cA), K^w_{\aA})$$
$$ \map(\psi(\cD), K_{\aA}) \stackrel{\sim}{\too} \Map(\psi(\cC)/\psi(\aA), K_{\aA})$$
are equivalences.
\end{proposition}

\begin{proof}
The argument is exactly the same in both cases. Therefore, we discuss
only the stable $K_{\aA}$.  Since $\cB$, $\cC$ and $\cD$
belong to $(\idemstabcat)^\omega$, the spectral Yoneda lemma shows us
that we need to prove that the induced sequence of spectra
\[
\xymatrix@-0.5pc{
K(\idemfun(\cD, \Idem(\cA))) \ar[r] & K(\idemfun(\aC, \Idem(\aA))) \ar[r]
& K(\idemfun(\aB, \Idem(\aA)))\\
}
\]
is a cofiber sequence.
Using corollary~\ref{cor:rep} it suffices to consider
the split-exact sequence of small spectral categories
\[
\xymatrix{
\rep(\cD, \cA) \ar@<-1ex>[r] & \rep(\aC, \cA)  \ar@<-1ex>[l]
\ar@<-1ex>[r] & \rep(\aB, \cA) \ar@<-1ex>[l] \,.
}
\]
Note that, again by corollary~\ref{cor:rep}, all of these
spectral categories carry a natural Waldhausen structure inherited
from the usual model structure on spectral modules.  We will
apply Waldhausen's fibration theorem \cite[1.6.4]{Wald}.  We have 
the Waldhausen category $v\rep(\aC,\cA)$, whose weak
equivalences are the morphisms $f$ such that $\Cone(f)$ is
contractible, as well as the Waldhausen category
$w\rep(\aC,\cA)$, with the same cofibrations as
$v\rep(\aC,\cA)$ but whose weak equivalences are those
$f$ such that $\Cone(f)$ belongs to
$\rep(\cD,\cA)$. Moreover, we have a natural inclusion
$v\rep(\aC,\cA) \subset w\rep(\aC,\cA)$ and an equivalence
$\rep(\aC,\cA)^w \simeq \rep(\cC,\cA)$; see~\cite[\S\,1.6]{Wald}. 
The conditions of~\cite[1.6.4]{Wald}
are satisfied, so we obtain a cofiber sequence of spectra
\[
K(\rep(\cD,\cA)) \too K(\rep(\aC,\cA)) \too K(\rep(\aB,\cA)).
\]
\end{proof}

Propositions~\ref{prop:suspensions} and \ref{prop:local}
allow us to prove theorem~\ref{thm:main2} as follows\,: let $\aA$ be stable $\i$-category and $\aB$ a compact small idempotent-complete stable $\i$-category. By proposition~\ref{prop:suspensions}
we have an equivalence $\Umot(\cA)\simeq K_{\cA}$ and by
proposition~\ref{prop:local} $K_{\cA}$ is local.  Therefore, we have
the following natural equivalence
$$ \Map(\Umot(\cB), \Umot(\cA)) \simeq \Map(\psi(\cB), K_{\aA})\,,$$
where the right-hand side is calculated in
$\Pre((\idemstabcat)^{\omega}; \ispec)$. Since $\cB$ belongs to
$(\idemstabcat)^{\omega}$, the presheaf $\psi(\cB)$ is representable
and so by the spectral Yoneda lemma we have $\Map(\psi(\cB),
K_{\aA}) \simeq K_{\cA}(\cB)$. Finally, since by definition of
$K_{\aA}$ we have $K_{\aA}(\cB)=K(\idemfun(\cB, \Idem(\cA)))$ the
proof is finished.

\section{Localization}\label{sec:localization}

The definition of additivity we study in this paper is given in terms
of the condition that algebraic $K$-theory takes models of split-exact
sequences of triangulated categories to (homotopy) cofiber sequences
of spectra.  This perspective is motivated in part by Neeman's
reformulation of the Thomason-Trobaugh localization
theorem \cite{NeemanLoc}.  Neeman observed that following
Thomason-Trobaugh and using the construction of Bousfield
localization, one could regard algebraic $K$-theory as in fact taking
exact sequences of triangulated categories to cofiber sequences of
spectra, provided one worked with non-connective $K$-theory (see
Theorem~\ref{thm:localization} for a version of this result).  
We will refer to such a theory as satisfying {\em localization}.  In this
section we construct the universal localizing invariant of small
stable $\infty$-categories; see theorem~\ref{thm:univ-loc}.  Our work
follows the general pattern of the analogous result for dg-categories
in \cite{CisTab}.

\begin{definition}\label{def:loc-invariant}
Let $\cD$ be a stable presentable $\i$-category. A functor
$$ E: \stabcat \too \cD$$
is called a {\em localizing invariant} of small stable $\i$-categories
if it inverts Morita equivalences (see definition~\ref{defn:iMor}),
preserves filtered colimits, and satisfies {\em localization}, \ie
sends exact sequences
\[
\cA \too \cB \too \cC
\]
of small stable $\i$-categories (see definition~\ref{def:exact}) to cofiber sequences
\[
E(\cA) \too E(\cB) \too E(\cC)
\]
in $\cD$.
We denote by $\Fun_{\mathrm{loc}}(\stabcat, \cD)$ the $\i$-category of
localizing invariants with values in $\cD$.
\end{definition}

Every localizing invariant is an additive invariant (see
definition~\ref{def:additive}), since a split-exact sequence is exact.
The converse does not hold, however: the impetus for the definition of
non-connective $K$-theory was precisely the fact that the connective
algebraic $K$-theory functor does not satisfy localization.  As we
discuss in sections~\ref{sec:non-connective} and~\ref{sec:cyclo},
non-connective algebraic $K$-theory ($\bbK$) and topological Hochschild
homology ($THH$) are localizing invariants.

Although the universal localizing invariant can be constructed by direct
localization, as in additive analogue of section~\ref{sec:additive}, we use a more involved procedure:

\begin{enumerate}
\item First, we construct a variant of the universal additive
invariant; see proposition~\ref{prop:k-additive}.  We work with a
general infinite regular cardinal $\kappa$, and we do not factor
through $\idemstabcat$; that is, Morita equivalences are not
inverted.  This produces the functor
\[
\uUmot\colon \stabcat \to \uMotadd.
\]

\item We localize $\uMotadd$ so that the exact sequences
$$ \cA \too \cB \too \cB/\cA\,,$$
with $\Ho(\cA)$ a thick triangulated subcategory of $\Ho(\cB)$, are sent
to cofiber sequences; see proposition~\ref{prop:k-Waldhausen}.  We
then obtain the functor 
\[
\uwUloc\colon \stabcat \to \uwMotloc.
\]

\item We perform a localization of $\uwMotloc$ to force Morita
equivalences to be sent to isomorphisms; see
proposition~\ref{prop:k-Morita}. We obtain then the functor
\[
\wUloc\colon \stabcat \to \wMotloc.
\]

\item Finally, we localize $\wMotloc$ so that the functor
$\wUloc$ preserves filtered colimits; see
theorem~\ref{thm:univ-loc}.  We end up with the universal localizing
invariant 
\[
\Uloc\colon \stabcat \to \Motloc.
\]
\end{enumerate}
The point of this seemingly circuitous process is that it enables
a clear, conceptual proof of the co-representability of
non-connective $K$-theory in $\Motloc$; see
section~\ref{sec:non-connective}.

\begin{notation}\label{not:important}
From now on and until the end of section~\ref{sec:non-connective} we
will work with a fixed infinite regular cardinal $\kappa$ larger than
$\omega$. We will denote by $(\stabcat)^\kappa$ the category of
$\kappa$-compact small stable $\i$-categories (see \S\ref{sec:compact}).

\end{notation}

\subsection{Additive $\kappa$-variant}\label{sec:k-add}

Let
\[
\psi\colon \idemstabcat \too \Pre((\idemstabcat)^{\kappa}; \ispec)
\]
be the functor obtained by first taking the Yoneda embedding and then
restricting the presheaves to the $\i$-category
$(\stabcat)^\kappa$.  Corollary~\ref{cor:setrep} allow us to choose a
fixed set $\underline{\cE_\mathrm{A}^\kappa}$ of representatives of split-exact
sequences in $(\idemstabcat)^\kappa$.  We denote by $\uMotadd$ the
localization of
$\Pre((\idemstabcat)^{\kappa}; \ispec)$ with respect to the set of maps
\[
\Cone(\psi(\cA) \to \psi(\cC)) \too \psi(\cB)\,,
\]
where $\cA \to \cC \to \cB$ is a split-exact sequence in
$\underline{\cE_\mathrm{A}^\kappa}$. Let $\uUmot$ be the following composite
\[
\idemstabcat \stackrel{\psi}{\too} \Pre((\idemstabcat)^{\kappa}; \ispec) \stackrel{\gamma}{\too} \uMotadd\,,
\]
where $\gamma$ is the localization functor.

\begin{proposition}\label{prop:k-additive}
The functor $\uUmot$ preserves $\kappa$-filtered colimits and sends split-exact sequences
\[
\xymatrix{
\cA \ar[r]_{i} & \cC \ar@<-1ex>[l]_f
\ar[r]_j & \cB \ar@<-1ex>[l]_{g} \,,
}
\]
in $\idemstabcat$ to (split) cofiber sequences
in $\uMotadd$.
Moreover, $\uUmot$ is universal with respect to these two properties, \ie given any stable presentable $\i$-category $\cD$, we have an equivalence of $\i$-categories
$$
(\uUmot)^{\ast}: \Fun^{\L}(\uMotadd, \cD) \stackrel{\sim}{\too} \Fun_{\underline{\mathrm{add}}}^\kappa(\idemstabcat, \cD)\,,$$
where the right-hand side denotes the full subcategory of
$\Fun(\stabcat, \cD)$ of morphisms of $\i$-categories which satisfy
the above two conditions.
\end{proposition}

\begin{proof}
The result follows from the analogue of the argument for
lemma~\ref{lem:filtrant} in the context of $\kappa$-compact objects
and $\Ind_{\kappa}$, and from the universal property of Bousfield
localization (see section~\ref{sec:bousloc} the functor $\psi$
preserves $\kappa$-filtered colimits and
proposition~\ref{prop:splitapprox} shows that any split-exact sequence
can be approximated by a $\kappa$-filtered colimit of split-exact
sequences in $\underline{\cE_\mathrm{A}^\kappa}$).
\end{proof}

Next, we localize $\uMotadd$ with respect to the set of maps

\begin{equation}\label{eq:maps1}
\Cone\left(\uUmot(\cA) \too \uUmot(\cB)\right) \too \uUmot(\cB/\cA)\,,
\end{equation}
where $\cA \to \cB \to \cB/\cA$ is a strict-exact sequence in
$\underline{\cE_\mathrm{wL}^\kappa}$ (see section \ref{sec:strict}). Let $\uwUloc$ be the following
composite
$$ \idemstabcat \stackrel{\uUmot}{\too} \uMotadd \stackrel{\gamma}{\too} \uwMotloc\,,$$
where $\gamma$ is the localization functor.

\begin{proposition}\label{prop:k-Waldhausen}
The functor $\uwUloc$ preserves $\kappa$-filtered colimits and
sends strict-exact sequences to cofiber sequences in $\uwMotloc$
\begin{eqnarray*}
\cA \too \cB \too \cB/\cA & \mapsto
& \uwUloc(\cA) \too \uwUloc(\cB) \too \uwUloc(\cB/\cA)\,.
\end{eqnarray*}
Moreover, $\uwUloc$ is universal with respect to these two
properties, \ie given any stable presentable $\i$-category $\cD$, we
have an equivalence of $\i$-categories
$$ (\uwUloc)^{\ast}\colon \Fun^{\L}(\uwMotloc, \cD) \stackrel{\sim}{\too} \Fun_{\underline{\mathrm{wloc}}}^\kappa(\idemstabcat, \cD)\,,$$
where the right-hand side denotes the full subcategory of
$\Fun(\idemstabcat, \cD)$ of morphisms of $\i$-categories which satisfy
the above two conditions.
\end{proposition}

\begin{proof}
The result follows from propositions~\ref{prop:k-additive}
and~\ref{prop:approx}, and from the universal property of localization
(see section~\ref{sec:bousloc}).
\end{proof}

\subsection{Morita equivalences}
We now localize $\uwMotloc$ with respect to the set of maps
\[
\uwUloc \left( \cA \too \Idem(\cA)\right)\,,
\]
where $\cA \to \Idem(\cA)$ belongs to $\cE^\kappa_\mathrm{L}$.
Let $\wUloc$ be the following
composition
\[
\stabcat \stackrel{\uwUloc}{\too} \uwMotloc \stackrel{\gamma}{\to} \wMotloc\,,
\]
where $\gamma$ is the localization functor.

\begin{proposition}\label{prop:k-Morita}
The functor $\wUloc$ inverts Morita equivalences, preserves $\kappa$-filtered colimits, and sends exact sequences to cofiber sequences in $\wMotloc$
\begin{eqnarray*}
\cA \too \cB \too \cC & \mapsto
& \wUloc(\cA) \to \wUloc(\cB) \to \wUloc(\cC) \,.
\end{eqnarray*}
Moreover, $\wUloc$ is universal with respect to these two properties, \ie given any stable presentable $\i$-category $\cD$, we have an equivalence of $\i$-categories
$$ (\wUloc)^{\ast}\colon \Fun^{\L}(\wMotloc, \cD) \stackrel{\sim}{\too} \Fun_{\mathrm{loc}}^\kappa(\stabcat, \cD)\,,$$
where the right-hand side denotes the full subcategory of $\Fun(\stabcat, \cD)$ of morphisms of $\i$-categories which satisfy the above three conditions.
\end{proposition}

\begin{proof}
The fact that $\wUloc$ preserves $\kappa$-filtered colimits is
clear.  Since a functor $\cA \to \cB$ is a Morita equivalence if
and only if  $\Idem(\cA) \to \Idem(\cB)$ is an
equivalence, proposition~\ref{prop:approx1} allow us to conclude that
$\wUloc$ inverts Morita equivalences.  We now show that $\wUloc$ sends
exact sequences to cofiber sequences. Let
\[
\cA \too \cB \too \cC
\]
be an exact sequence. Since we have an induced Morita equivalence
$\cB/\cA \to \cC$, it suffices to show that $\wUloc$ sends exact
sequences of shape
\[
\cA \too \cB \too \cB/\cA
\]
to cofiber sequences in $\wMotloc$. Since $\Idem:\stabcat\to\idemstabcat$ is a localization, it commutes with colimits and therefore the right-hand vertical map in the diagram
\[
\xymatrix{
\cA \ar[r] \ar[d] & \cB \ar[d] \ar[r] & \cB/\cA \ar[d] \\
\Idem(\cA) \ar[r] & \Idem(\cB) \ar[r] & \Idem(\cB)/\Idem(\cA)\,,
}
\]
is a Morita equivalence.  The bottom line is a strict-exact sequence, and so we conclude
that $\wUloc$ sends exact sequences to cofiber sequences.
Finally, the universality of $\wUloc$ follows from
propositions~\ref{prop:k-additive} and \ref{prop:k-Waldhausen}, and
from the universal property of localization (see
section~\ref{sec:bousloc}).
\end{proof}

\subsection{Universal localizing invariant}\label{sec:filtered}

We denote by $\cE_\L$ the set of maps in $\wMotloc$ of shape
\begin{equation*}
\underset{\alpha}{\mathrm{colim}}\, \wUloc(\cA_{\alpha}) \too \ \wUloc(\cA)\,,
\end{equation*}
where $\{\cA_{\alpha}\}$ is a filtered diagram of objects in $(\stabcat)^\omega$ whose
colimit is a $\kappa$-compact small stable $\i$-category $\cA$. Localize
$\wMotloc$ with respect to the set $\cE_\L$.  Let $\Uloc$ be the
following composition
\[
\stabcat \stackrel{\wUloc}{\too} \wMotloc \stackrel{\gamma}{\too} \Motloc\,,
\]
where $\gamma$ is the localization functor.

\begin{theorem}\label{thm:univ-loc}
The functor $\Uloc$ is the {\em universal localizing invariant}, \ie
given any stable presentable $\i$-category 
$\cD$, we have an equivalence of $\i$-categories 
$$(\Uloc)^{\ast}\colon \Fun^{\L}(\Motloc, \cD) \stackrel{\sim}{\too} \Fun_{\mathrm{loc}}(\stabcat, \cD)\,.$$
\end{theorem}

\begin{proof}
Let us denote by $\Motadd^\kappa$ the small stable $\infty$-category
constructed as in section \ref{sec:additive} but where we use
$(\idemstabcat)^\kappa$ instead of $(\idemstabcat)^\omega$. Similarly
to $\Uadd$ we have a well-defined functor
$\Uadd^\kappa \colon \stabcat \to \Motadd^\kappa$ and so by performing a
localization analogous to the one of subsection~\ref{sec:filtered}
(with $\cE_A$ instead of $\cE_\L$) we obtain a small stable
$\infty$-category which we denote by $\Motadd^\omega$ and a composed
functor  
$$\Uadd^\omega\colon \stabcat \stackrel{\Uadd^\kappa}{\too} \Motadd^\kappa \stackrel{\gamma}{\too} \Motadd^\omega\,.$$ 
Let us start by showing that $\Motadd^\omega$ agrees with $\Motadd$
and that $\Uadd^\omega$ agrees with $\Uadd$. For this (and because of
the universal property of $\Uadd$ and $\Uadd^\omega$) it suffices to
show that $\Uadd^\omega$ preserves filtered colimits. Consider the
composite  
$$
(\idemstabcat)^\omega \hookrightarrow
(\idemstabcat)^\kappa \subset \stabcat \stackrel{\Uadd^\kappa}{\too} \Motadd^\kappa\,.$$ 
It inverts Morita equivalences and sends
split-exact sequences to split cofiber sequences. By the construction
of $\Motadd$ we obtain then an additive invariant
$\stabcat \to \Motadd^\kappa$ and hence by the universal property of
$\Uadd$ a colimit preserving functor
$\Phi \colon \Motadd \to \Motadd^\kappa$ and a natural transformation
$\eta\colon \Phi \circ \Uadd \Rightarrow \Uadd^\kappa$. We 
now observe that the two functors 
\begin{eqnarray*}
\stabcat \stackrel{\Umot^\kappa}{\too} \Motadd^\kappa \stackrel{\gamma}{\too} \Motadd^\omega
&& \stabcat \stackrel{\Umot}{\too} \Motadd \stackrel{\Phi}{\too} \Motadd^\kappa \stackrel{\gamma}{\too} \Motadd^\omega
\end{eqnarray*}
agree. Since they preserve $\kappa$-filtered colimits and every object
in $\stabcat$ can be expressed as a $\kappa$-filtered colimit of
$\kappa$-small objects, it suffices to show that they agree for every
$\kappa$-compact small stable $\infty$-category 
$\cA$.  The stable $\infty$-category $\cA$ can be expressed as a filtered colimit
$\underset{\alpha}{\mathrm{colim}} (\cA_\alpha) \to \cA$, with 
$\cA_\alpha \in (\stabcat)^\omega$, and the evaluation of the
natural transformation $\eta$ at $\cA$ identifies with
$$\underset{\alpha}{\mathrm{colim}}\,\, \Umot^\kappa(\cA_\alpha) \too \Umot^\kappa(\cA)\,.$$ 
Since these map belongs to $\cE_A$, they become invertible in
$\Motadd^\omega$, and so we conclude that the above two functors
agree. Since the one on the right-hand side preserves filtered
colimits we conclude that $\Uadd^\omega$ also preserves filtered
colimits. This shows that $\Uadd^\omega$ agrees with $\Uadd$ (and
hence that $\Motadd^\omega$ agrees with $\Motadd$).  Now, let
$\Motloc^\omega$ be the category defined as $\Motloc^\kappa$ but 
with $\kappa$ replaced by $\omega$. Clearly, the associated functor
$\Uloc^\omega$ is the universal localizing invariant and so in order
to conclude the proof of Theorem~\ref{thm:univ-loc} it suffices to
show that $\Motloc$ agrees with $\Motloc^\omega$ and that $\Uloc$
agrees with $\Uloc^\omega$. Starting with $\Motadd^\kappa$ we can
perform the following two localizations  
\begin{eqnarray*}
\Motadd^\kappa \too \Motloc^\kappa \too \Motloc^\omega
&& \Motadd^\kappa \too \Motadd^\omega \simeq \Motadd \too \Motloc\,. 
\end{eqnarray*}
Since these localizations are independent of the order in which they
are performed, our claim follows and so the proof is finished.
\end{proof}

\section{Non-connective $K$-theory}\label{sec:non-connective}

Bass introduced the negative $K$-groups in order to measure the
failure of $K_0$ and $K_1$ to satisfy localization; this
perspective was studied in detail in Thomason-Trobaugh and led to the
definition of the Bass-Thomason non-connective $K$-theory spectrum of
rings and schemes.  In fact, any nontrivial theory which is ``like $K$-theory'' and 
satisfies localization must be non-connective; there is a nice
discussion of this in \cite{keller}.
In this section we introduce the non-connective algebraic $K$-theory of $\i$-categories
and show that it becomes co-representable in $\Motloc$; see
theorem~\ref{thm:main3}. This result depends critically on the multistage construction of Section~\ref{sec:localization}, which again follows the general pattern of the argument for dg-categories in \cite{CisTab}. For a connective ring spectrum $R$, we give a slightly different definition of the non-connective $K$-theory in terms of a ``suspension ring spectrum'' of $R$, and use this show that the negative $K$-groups of $R$ are isomorphic to those of $\pi_0 R$. 

\subsection{Non-connective $K$-theory of
$\i$-categories}\label{sub:non-conn}

In order to construct the non-connective $K$-theory spectrum associated
to a small stable $\i$-category, we use a generalization of the
axiomatic framework due to Schlichting \cite{Marco}.  For an
{\em uncountable} regular cardinal $\kappa$, we will produce functors
$\kcone$ and $\ksusp$ from $\stabcat$ to $\stabcat$ such that for any
small stable $\i$ category $\aA$:
\begin{enumerate}
\item $K(\kcone \aA)$ is contractible,
\item there are natural transformations
\[
\Id \to \kcone \to \ksusp
\]
such that $\aA \to \kcone \aA \to \ksusp \aA$ is exact,
\item the functors $\kcone$ and $\ksusp$ preserve exact sequences,
\item and $\kcone$ and $\ksusp$ preserve $\kappa$-filtered colimits in
$\stabcat$.
\end{enumerate}

The idea is that $\kcone \aA$ is a ``$K$-theoretic cone'' and so
$\ksusp \aA$ is a ``suspension'' of $\aA$.  Fix an uncountable regular
cardinal $\kappa$, and for a stable $\i$-category $\aC$ recall from
Section~\ref{sec:compact} that $\aC^{\kappa}$ denotes the
$\kappa$-compact objects in $\aC$.

\begin{definition}\label{def:setup}
Using proposition~\ref{prop:Ind}, we define $\kcone \cA
= (\Ind_\omega(\cA))^{\kappa}$ and $\ksusp \cA$ to be the
cofiber $(\Ind_\omega(\cA))^{\kappa}/\cA$.  
\end{definition}

\begin{remark}
One might wish to simply use $\Ind_{\omega} \aA$ as the cone
construction; however, this will rarely turn out to be a small
$\i$-category, whereas passing to the $\kappa$-compact objects yields
an (essentially) small $\i$-category by construction.
\end{remark}

Observe that $\kcone$ is a composite functor
\begin{equation}\label{eq:composed-functor}
\stabcat \to \stabprcat_\omega \to \kstabcat\to\stabcat.
\end{equation}
By construction and Propositions~\ref{prop:loc=cof}
and \ref{prop:quotient}, we have an exact sequence
\[
\cA \too \kcone \cA \too \ksusp \cA,
\]
which is natural in small stable $\i$-categories $\cA$.  Next, we
check that $\kcone \cA$ satisfies property (i) above.

\begin{lemma}
Let $\cA$ be a small stable $\i$-category.  Then $K(\kcone \cA)$
is trivial.
\end{lemma}

\begin{proof}
Since $\kappa$ is uncountable, $\kcone \cA$ has countable coproducts,
and so the usual Eilenberg swindle argument implies that the identity
map is null-homotopic on $K$-theory and so its $K$-theory vanishes.
Specifically, the functor $F \colon \kcone \cA \to \kcone \cA$ defined
by $X \mapsto \coprod_{\bN} X$ is exact.  Moreover, there is a natural
equivalence of exact functors $\id \coprod F \htp F$ induced by the
equivalence $X \coprod (\coprod_{\bN} X) \htp \coprod_{\bN} X$.  Applying
$K$-theory, we can split off the $F$ component of the resulting equivalence of spectra and deduce that the identity of $\kcone \cA$
is null-homotopic.
\end{proof}

We must check that $\kcone$ and $\ksusp$
preserve exact sequences of small stable $\i$-categories.  

\begin{proposition}
Let $\aA \to \aB \to \aC$ be an exact sequence of small stable
$\i$-categories.  Then the induced sequences
\[
\xymatrix{
\kcone \aA \ar[r] & \kcone \aB \ar[r] & \kcone \aC
&& \ksusp \aA \ar[r] & \ksusp \aB \ar[r]
& \ksusp \aC \\
}
\]
are exact.
\end{proposition}

\begin{proof}
It suffices to show the result for $\kcone$, as the statement
for $\ksusp$ follows because colimits commute.
Thus, we need to verify that
\[
(\Ind_\omega (\aA))^{\kappa} \to (\Ind_\omega (\aB))^{\kappa} \to (\Ind_\omega (\aC))^{\kappa}
\]
is exact.  The sequence
\[
\Ind_\omega \aA \to \Ind_\omega \aB \to \Ind_\omega \aC
\]
is exact by Definition~\ref{def:exact} and
Proposition~\ref{prop:stabquot}.  Now the result follows from
Proposition~\ref{prop:notyet}.
\end{proof}

Passing to the triangulated homotopy category by composing with the
functor $\Ho$, we get a series of functors which satisfies
Schlichting's setup of \cite[\S2.2]{Marco} and so produces negative
$K$-groups.  Furthermore, we can define the non-connective $K$-theory spectrum as
follows, following \cite[\S12]{Marco}.

\begin{definition}\label{def:knonconnspec}
Let $\cA$ be a small stable $\i$-category. Its non-connective
$K$-theory spectrum $\bbK(\cA)$ is given by
$$\bbK(\cA) := \colim_n \Omega^n K(\ksusp^{(n)}(\cA)) \,.$$
Here, $K$ stands for the $K$-theory spectrum of \S\ref{sec:back-K}, and
the structure maps are induced from the exact sequences
$$ \ksusp^{(n)}(\cA) \too \kcone \ksusp^{(n)}(\cA) \too \ksusp^{(n+1)}(\cA) \qquad
n \geq 0\,.$$
\end{definition}

Schlichting's axiomatic framework implies that this construction
agrees with his when both are defined, and therefore we deduce from
his comparison results \cite[\S 8]{Marco} that the non-connective
$K$-theory spectrum of Definition~\ref{def:knonconnspec} agrees with
the various classical constructions of non-connective $K$-theory
spectra.  

Finally, we establish the final technical condition; this will be
needed in the following sections.

\begin{lemma}\label{lem:Ffilt}
The functors $\kcone$ and $\ksusp$ preserve $\kappa$-filtered
colimits. 
\end{lemma}

\begin{proof}
Recall that $\kcone$ is the composite \eqref{eq:composed-functor}. Hence, the claim follows from the
fact that the passage to $\Ind_\kappa$ and to $\kappa$-compact objects
preserves $\kappa$-filtered colimits~\cite[5.5.7.8, 5.5.7.10,
5.5.7.11]{HTT}.  Since $\ksusp$ is the cofiber of the inclusion
$\aA \to \kcone$ and colimits commute, we deduce that $\ksusp$
preserves $\kappa$-filtered colimits if $\kcone$ does.
\end{proof}

\subsection{Co-representability}

This subsection is entirely devoted to the proof of the following
co-representability result.

\begin{theorem}\label{thm:main3}
Let $\aA$ be a small stable $\i$-category. Then there is a natural
equivalence of spectra
\[\Map(\Uloc(\ispec^{\omega}),\,\Uloc(\aA)) \simeq \bbK(\aA)\,.\]
In particular, for each integer $n$, we have isomorphisms of abelian groups
\[\Hom(\Uloc(\ispec^{\omega}),\, \Sigma^{-n}\Uloc(\cA)) \simeq \bbK_n(\cA)\]
in the triangulated category $\Ho(\Motloc)$.
\end{theorem}
The proof of theorem~\ref{thm:main3} will follow from
theorems~\ref{thm:key1} and \ref{thm:key2}, and from
propositions~\ref{prop:negK} through \ref{prop:negK1}.

\begin{theorem}\label{thm:key1}
Let $\aA$ and $\cB$ be small stable $\i$-categories such that $\cB$ is
$\kappa$-compact. Then there is a natural equivalence of spectra
\[\Map(\uUmot(\aB),\,\uUmot(\aA)) \simeq K(\Fun^{ex}(\aB, \aA))\,.\]
If $\cB=\ispec^\omega$ is the $\i$-category of
compact spectra, this reduces to an equivalence
\[\Map(\uUmot(\ispec^\omega),\,\uUmot(\aA)) \simeq K(\aA)\,.\]
\end{theorem}

\begin{proof}
The proof is analogous to the argument for theorem~\ref{thm:main2};
instead of the idempotent-complete stable $\i$-category
$\idemfun(\cB, \Idem(\cA))$ we consider the small stable 
$\i$-category $\idemfun(\cB, \cA)$.  Note that since
$\kappa > \omega$, $\ispec^\omega$ belongs to $(\stabcat)^\kappa$.
\end{proof}

\begin{theorem}\label{thm:key2}
Let $\cA$ be a small stable $\i$-category.  Then there is a natural
equivalence of spectra
\[\Map(\uwUloc(\ispec^\omega),\,\uwUloc(\aA)) \simeq K(\aA)\,.\]
\end{theorem}

\begin{proof}
By construction, the object $\uUmot(\ispec^\omega)$ is compact in
$\uMotadd$.  Let $S$ denote the set of maps in (\ref{eq:maps1}),
$\overline{S}$ the strongly saturated collection of arrows generated
by $S$ \cite[5.5.4.5]{HTT}, and let $X$ be an $S$-local
object such that the map $\uUmot(\aA)\to X$ is an $S$-local
equivalence (i.e., $\uUmot(\aA)\to X$ is in $\overline{S}$).
Then by definition,
\[
\Map(\uwUloc(\ispec^\omega),\uwUloc(\aA))\simeq\Map(\uUmot(\ispec^\omega),X),\]
so it suffices to show that the functor
\begin{equation}\label{eq:functor}
R:=\Map(\uUmot(\ispec^\omega), -): \uMotadd \too \ispec
\end{equation}
sends the maps in $\overline{S}$ to equivalences of
spectra.  Since $\uMotadd$ is a stable $\i$-category and
$\uUmot(\ispec^\omega)$ is compact, $R$ preserves small colimits,
so the two-out-of-three property allows us to reduce to checking that
$R$ sends the elements of $S$ to equivalences.

Consider the following diagram
\begin{equation}\label{eq:dia-key2}
\xymatrix{
\uUmot(\cA) \ar[r] \ar@{=}[d] & \uUmot(\cB) \ar[r] \ar@{=}[d] &  \uUmot(\cB)/\uUmot(\cA) \ar[d] \\
\uUmot(\cA) \ar[r] & \uUmot(\cB) \ar[r] & \uUmot(\cB/\cA)\,.
}
\end{equation}
By applying the functor (\ref{eq:functor}) to the above diagram
(\ref{eq:dia-key2}) we obtain by theorem~\ref{thm:key1} a diagram in
$\Spt$
\begin{equation}\label{eq:dia1-key2}
\xymatrix{
K(\cA) \ar[r] \ar@{=}[d] & K(\cB) \ar[r] \ar@{=}[d] &  K(\cB)/K(\cA) \ar[d] \\
K(\cA) \ar[r] & K(\cB) \ar[r] & K(\cB/\cA)\,,
}
\end{equation}
where the upper row is a homotopy cofiber sequence. Now, an
argument analogous to the one used in the proof of
proposition~\ref{prop:local} (where we make use of Waldhausen's
fibration theorem) allow us to conclude that the lower row in the
above diagram (\ref{eq:dia1-key2}) is also a homotopy
cofiber sequence.  This completes the argument.
\end{proof}

Let ${\bf V}$ be the partially ordered set $\{(i,j): |i-j|\leq 1,\,
i,j \geq 0 \} \subset \bbN \times \bbN$. Given a small stable
$\i$-category $\cA$, we denote by $\mathrm{Dia}(\cA)$ the $\N({\bf
V})$-diagram

\begin{equation}\label{eq:dia-spec}
\xymatrix{
 &  &  & \cdots \\
 & \kcone \ksusp^{(n)}(\cA) \ar[r] & \ksusp^{(n+1)}(\cA) \ar[u] \ar[r] &  \\
 \ar[r] & \ksusp^{(n)}(\cA) \ar[u] \ar[r] & \ksusp^{(n)}(\cA)/\ksusp^{(n)}(\cA) \ar[u] & \\
\cdots & \ar[u] & & \,,
}
\end{equation}

where $\kcone$ and $\ksusp$ are as in Definition~\ref{def:setup}.

\begin{lemma}\label{lem:trivial}
Let $\cA$ be a small stable $\i$-category. Then
$\ksusp^{(n)}(\cA)/\ksusp^{(n)}(\cA)$ and $\kcone\ksusp^{(n)}(\cA)$
become trivial after application of $\uwUloc$.
\end{lemma}

\begin{proof}
The object $\ksusp^{(n)}(\cA)/\ksusp^{(n)}(\cA)$ is already trivial in
$\stabcat$.  Since Proposition~\ref{prop:Ind} implies that
$\kcone\ksusp^{(n)}(\cA)$ admits all $\kappa$-small colimits, for any
$\cB$ in $(\stabcat)^\kappa$ the small stable $\i$-category
$\Fun(\cB, \kcone \ksusp^{(n)}(\cA))$ also admits all $\kappa$-small
colimits.  Thus, the connective $K$-theory spectrum
$K(\Fun(\cB, \kcone \ksusp^{(n)}(\cA))$ is trivial.  Finally,
theorem~\ref{thm:key1} and the fact that the objects $\uUmot(\cB)$,
with $\cB$ in $(\stabcat)^\kappa$ generate the category
$\uMotadd$ \cite[5.5.7.3]{HTT} allow us to conclude that
$\kcone \ksusp^{(n)}(\cA)$ becomes trivial after application of
$\uUmot$, and thus after application of $\uwUloc$.
\end{proof}

Let $\cA$ be a small stable $\i$-category.  We denote by $V(\cA)$ the object
\[
V(\aA) = \colim_n\Sigma^{-n}\uwUloc(\ksusp^{(n)}(\cA))
\]
in $\uwMotloc$ whose indexing maps are induced from the above diagram
(\ref{eq:dia-spec}). Note that $V(\cA)$ is functorial in $\cA$ and
that we have a natural map $\uwUloc(\cA) \to V(\cA)$. We obtain then a
well-defined functor $V$ along with a natural transformation:
\begin{eqnarray}\label{eq:Func-V}
V(-): \stabcat \too \uwMotloc && \uwUloc \Rightarrow V(-)\,.
\end{eqnarray}
\begin{proposition}\label{prop:negK}
Let $\cA$ be a small stable $\i$-category. Then, there is a natural
equivalence of spectra
$$ \Map(\uwUloc(\ispec^\omega),\, V(\cA)) \simeq \bbK(\cA)\,.$$
\end{proposition}
\begin{proof}
This follows from the following equivalences
\begin{eqnarray}
\bbK(\cA) & = & \underset{n \geq 0}{\mathrm{hocolim}}\, \Omega^n
K(\ksusp^{(n)}(\cA))  \nonumber \\ 
& \simeq & \underset{n \geq 0}{\mathrm{hocolim}}\, \Omega^n \Map(\uwUloc(\ispec^\omega),\, \uwUloc(\ksusp^{(n)}(\cA))) \label{eq:numberA}\\
& \simeq & \Map(\uwUloc(\ispec)^\kappa),\, \underset{n \geq
0}{\mathrm{colim}}\, \Sigma^{-n} \uwUloc(\ksusp^{(n)}(\cA))) \label{eq:numberB}\\
& \simeq &  \Map(\uwUloc(\ispec^\kappa),\, V(\cA))\,. \nonumber
\end{eqnarray}
Equivalence (\ref{eq:numberA}) comes from theorem~\ref{thm:key2} and equivalence (\ref{eq:numberB}) comes from the compactness of $\uwUloc(\ispec^\omega)$ in $\uwMotloc$.
\end{proof}
\begin{proposition}\label{prop:Mor-V}
The functor $V$ (\ref{eq:Func-V}) inverts Morita equivalences.
\end{proposition}
\begin{proof}
It suffices to show that $V(-)$ sends maps of shape
$\cA \to \Idem(\cA)$ to isomorphisms.  Consider the following diagram
$$
\xymatrix{
\cA \ar[d]_P \ar[r] & \kcone (\cA) \ar[d]^{\kcone(P)} \ar[r] & \ksusp(\cA) \ar[d]^{\ksusp(P)} \\
\Idem(\cA) \ar[r] & \kcone(\Idem(\cA)) \ar[r] & \ksusp(\Idem(\cA))\,.
}
$$
Proposition~\ref{prop:Ind} implies that $\kcone(P)$ is an
equivalence. Therefore, since both rows are strict-exact sequences and
$\Ho(\cA)$ and $\Ho(\Idem(\cA))$ differ by direct summands, we
conclude that $\ksusp(P)$ is an equivalence. The definition of the
functor $V(-)$ allow us to conclude the proof.
\end{proof}
\begin{proposition}\label{prop:V-loc}
The functor
\[
V(-) \colon \stabcat \too \uwMotloc
\]
inverts Morita equivalences, preserves $\kappa$-filtered colimits, and
sends exact sequences to cofiber sequences. 
\end{proposition}
\begin{proof}
Proposition~\ref{prop:Mor-V} implies that $V(-)$ inverts Morita
equivalences.  Furthermore, by Lemma \ref{lem:Ffilt}, $\kcone$ and $\ksusp$ preserve
$\kappa$-filtered colimits for $\kappa > \omega$,
and so $V(-)$ does as well. Now, let
\[
\cA \too \cB \too \cC
\]
be an exact sequence.  Proposition~\ref{prop:Mor-V} implies that we
can assume that $\Ho(\cA)$ is a thick triangulated subcategory of
$\Ho(\cB)$. Consider the following diagram
\begin{equation}\label{eq:diag:quot}
\xymatrix{
\mathrm{Dia}(\cA) \ar@{=}[d] \ar[r] & \mathrm{Dia}(\cB) \ar@{=}[d] \ar[r] & \mathrm{Dia}(\cA,\cB) := \mathrm{Dia}(\cA)/\mathrm{Dia}(\cB) \ar[d]^D \\
\mathrm{Dia}(\cA) \ar[r] & \mathrm{Dia}(\cB) \ar[r] & \mathrm{Dia}(\cC)\,,
}
\end{equation}
where $\mathrm{Dia}(\cA, \cB)$ is obtained by passage to the cofiber
objectwise. Note that since in the above diagram (\ref{eq:diag:quot})
the upper row is objectwise a strict-exact sequence, we obtain a
cofiber sequence
$$ V(\cA) \too V(\cB) \too V(\cB, \cA) \too \Sigma V(\cA)$$
in $\uwMotloc$, where
$$ V(\cB, \cA):= \colim_n\Sigma^{-n}\uwUloc (\ksusp^{(n)}(\cB)/\ksusp^{(n)}(\cA))\,.$$
We now show that the induced map
\begin{equation}\label{keymorph}
V(\aB,\aA)\too V(\aC)
\end{equation}
is an equivalence. For this,
consider the following commutative diagram
$$
\xymatrix{
\ksusp^{(n)}(\cA) \ar[d] \ar[r] & \kcone\ksusp^{(n)}(\cA) \ar[d] \ar[r] & \ksusp^{(n+1)}(\cA) \ar[d] \\
\ksusp^{(n)}(\cB) \ar[d] \ar[r] & \kcone\ksusp^{(n)}(\cB) \ar[d] \ar[r] & \ksusp^{(n+1)}(\cB) \ar[d] \\
\ksusp^{(n)}(\cB)/\ksusp^{(n)}(\cA)  \ar[d] \ar[r] & \kcone\ksusp^{(n)}(\cB)/\kcone\ksusp^{(n)}(\cA)  \ar[d]^{\theta_n} \ar[r] & \ksusp^{(n+1)}(\cB)/\ksusp^{(n+1)}(\cA) \ar[d]^{D_n} \\
\ksusp^{(n)}(\cC) \ar[r] & \kcone\ksusp^{(n)}(\cC) \ar[r] & \ksusp^{(n+1)}(\cC)\,. \\
}
$$
Since the induced triangulated functor
$$ \Ho (\kcone \ksusp^{(n)}(\cA)) \too \Ho (\kcone \ksusp^{(n)}(\cB))$$
preserves $\kappa$-small colimits, \cite[\S3.1]{Marco} implies that
the triangulated category
$$ \Ho (\kcone\ksusp^{(n)}(\cB) / \kcone\ksusp^{(n)}(\cA) )$$
is idempotent complete. Therefore, $\theta_n$ is an
equivalence, and we obtain maps
\[
\psi_n:\Sigma_\kappa^{(n)}(\aC)\to\kcone\Sigma_\kappa^{(n)}(\aB)/\kcone\Sigma_\kappa^{(n)}(\aA)
\]
which induce maps
\[
\Psi_n:\Sigma^{-n}\uwUloc(\Sigma_\kappa^{(n)}(\aC))\to\Sigma^{-n-1}\uwUloc(\Sigma_\kappa^{(n+1)}(\aB)/\Sigma^{(n+1)}_\kappa(\aA)).
\]
It follows that the natural map
\[
\colim_n\Sigma^{-n}\uwUloc(\Sigma_\kappa^{(n)}(\aB)/\Sigma_\kappa^{(n)}(\aA))\to\colim_n\Sigma^{-n}\uwUloc(\Sigma_\kappa^{(n)}(\aC))
\]
is an equivalence, which implies that the map (\ref{keymorph}) is an equivalence.
\end{proof}

\begin{corollary}[of proposition~\ref{prop:V-loc}]\label{cor:V-loc}
There is a functor $$ \mathrm{Loc}: \wMotloc \too \uwMotloc$$ such
that $\mathrm{Loc}(\wUloc(\cA)) \simeq V(\cA)$, for every small stable
$\i$-category $\cA$.
\end{corollary}

\begin{proof}
This follows from propositions~\ref{prop:V-loc} and \ref{prop:k-Morita}.
\end{proof}
\begin{proposition}\label{prop:can-isom}
The two functors
\[
\mathrm{Loc}, \gamma^{\ast} \colon \wMotloc \too \uwMotloc
\]
are canonically equivalent, where $\gamma^{\ast}$ is the right adjoint
of the localization functor. 
\end{proposition}

\begin{proof}
Let us denote by ${\bf L}$ the endofunctor $\Loc \circ \gamma$ of
$\uwMotloc$.  Note that we have a natural transformation
$\Id \Rightarrow {\bf L}$. Making use of the definition of
$V(-)$ and of the fact that colimits in $\i$-categories commute, we
observe that ${\bf L}$ is a localization functor on
$\uwMotloc$ \cite[5.2.7.4]{HTT}. Therefore, it suffices 
to show that a map in $\uwMotloc$ becomes an equivalence in $\wMotloc$
if and only if it becomes an equivalence after application of ${\bf
L}$. This follows from the fact that for every small stable
$\i$-category $\cA$, we have an equivalence
$\gamma(V(\cA)) \simeq \wUloc(\cA)$: note that we have cofiber
sequences in $\wMotloc$ 
$$ \wUloc(\ksusp^{(n)}(\cA)) \too \wUloc(\kcone\ksusp^{(n)}(\cA)) \too \wUloc(\ksusp^{(n+1)}(\cA))  \,.$$ 
\end{proof}
\begin{proposition}\label{prop:negK1}
Let $\cA$ be small stable $\i$-category.  We have a natural
isomorphism in the stable homotopy category of spectra 
$$ \Map(\wUloc(\ispec^\omega),\, \wUloc(\cA)) \simeq \bbK(\cA)\,.$$
\end{proposition}
\begin{proof}
This follows from the following equivalences
\begin{eqnarray}
 \Map(\wUloc(\ispec^\omega), \wUloc(\cA)) & \simeq &  \Map(\uwUloc(\ispec^\omega), \gamma^{\ast}(\wUloc(\cA))) \nonumber \\
& \simeq & \Map(\uwUloc(\ispec^\omega), \mathrm{Loc}(\wUloc(\cA))) \label{eq:number4} \\
& \simeq & \Map(\uwUloc(\ispec^\omega), V(\cA)) \label{eq:number5}\\
& \simeq & \bbK(\cA)\label{eq:number6}\,.
\end{eqnarray}
Equivalence (\ref{eq:number4}) comes from proposition~\ref{prop:can-isom}, equivalence (\ref{eq:number5}) comes from corollary~\ref{cor:V-loc}, and equivalence (\ref{eq:number6}) is proposition~\ref{prop:negK}.
\end{proof}
\begin{proof}[Proof of theorem~\ref{thm:main3}]
Recall from subsection~\ref{sec:filtered} that $\Motloc$ is obtained
by localizing $\wMotloc$ with respect to the set $\cE_\L$. Since
$\wUloc(\ispec^\omega)$ is compact in $\wMotloc$, it is sufficient by
proposition~\ref{prop:negK1} and the universal property of
localization (see section~\ref{sec:bousloc}) to show that the functor
\[
\Map(\wUloc(\ispec^\omega),-)\colon \wMotloc \too \ispec
\]
sends the elements of $\cE_\L$ to equivalences. This follows from the
fact that the non-connective $K$-theory construction preserves
filtered colimits (see \cite[\S7, Lemma~6]{Marco}), and so the proof
is finished. 
\end{proof}

\subsection{Non-connective $K$-theory of Waldhausen categories and localization}
In particular, theorem~\ref{thm:main3} implies that non-connective
$K$-theory satisfies localization.  This is an extremely useful fact
in practice; localization sequences provide one of the main
computation tools for understanding algebraic $K$-theory.  As such, we
state a version of this result in terms of Waldhausen categories.  We
begin by defining the non-connective $K$-theory of a Waldhausen
category.

\begin{definition}\label{def:ncKwald}
Let $\aC$ be a DHKS-saturated Waldhausen category with factorization.
Then the non-connective $K$-theory $\bbK(\aC)$ of $\aC$ is defined as
the non-connective $K$-theory $\bbK(\N(\aC)[W^{-1}][\Sigma^{-1}])$ of
the $\infty$-category 
\[
\N(\aC)[W^{-1}][\Sigma^{-1}]\simeq\colim\{\N(\aC)[W^{-1}]\overset{\Sigma}{\too}\N(\aC)[W^{-1}]\overset{\Sigma}{\too}\cdots\}
\]
obtained by inverting the suspension on the underlying $\i$-category
$\N(\aC)[W^{-1}]$ in the $\infty$-category $\Cat_\infty^\mathrm{Rex}$
of $\infty$-categories with finite colimits and right-exact functors. 
\end{definition}

This definition in terms of the stabilization is reasonable because of
the following consistency results. 

\begin{proposition}\label{prop:stab=Sigma^-1}
Let $\C$ be a presentable $\i$-category with a zero object, and let
$\C^\omega[\Sigma^{-1}]$ denote the colimit 
\[
\C^\omega[\Sigma^{-1}]\simeq\colim\{\C^\omega\overset{\Sigma}{\too}\C^\omega\overset{\Sigma}{\too}\cdots\}
\]
in $\Cat_\infty^\mathrm{Rex}$.
Then $\C^\omega$ is stable, and the induced functor
\[
\C^\omega[\Sigma^{-1}]\too\Stab(\C)
\]
identifies the idempotent-completion of $\C^\omega[\Sigma^{-1}]$ with $\Stab(\C)^\omega$.
\end{proposition}
\begin{proof}
Let $\aD$ be an idempotent-complete stable $\infty$-category.
Then
\begin{align*}
\Fun^{\ex}(\C^\omega[\Sigma^{-1}],\aD)&\simeq\lim\Fun^{\ex}(\C^\omega,\aD)\simeq\lim\Fun^{\L}_\omega(\C,\Ind(\aD))\\
&\simeq\Fun^{\L}_\omega(\Stab(\C),\Ind(\aD))\simeq\Fun^\mathrm{\ex}(\Stab(\C)^\omega,\aD).
\end{align*}
Since $\Stab(\C)^\omega$ is necessarily idempotent-complete, we conclude that it is equivalent to the idempotent-completion of $\C^\omega[\Sigma^{-1}]$.
\end{proof}

\begin{proposition}\label{prop:connconsistency}
Let $\aC$ be a DHKS-saturated Waldhausen category with factorization.
Then the natural map $\N(\aC)[W^{-1}]\to\\N(\aC)[W^{-1}][\Sigma^{-1}]$
induces a natural equivalence 
\[
K(\aC)\too K(\N(\aC)[W^{-1}][\Sigma^{-1}]).
\]
\end{proposition}

\begin{proof}
The additivity theorem implies that, for Waldhausen categories with
factorization, the suspension endomorphism $\Sigma:\aC\to\aC$ induces
$-\id:K(\aC)\to K(\aC)$.  By naturality, we conclude that
$\Sigma:\N(\aC)[W^{-1}])\to\N(\aC)[W^{-1}]$ acts invertibly on
$K$-theory. Finally, since $K$-theory (viewed as a functor of small
$\infty$-categories with finite colimits and a zero object and
right-exact functors) preserves filtered colimits, we see that 
\[
K(\N(\aC)[W^{-1}][\Sigma^{-1}])\simeq\colim K(\N(\aC)[W^{-1}])\simeq
K(\N(\aC)[W^{-1}])\simeq K(\aC), 
\]
where the last equivalence follows from Corollary~\ref{cor:gencomp}.
\end{proof}

\begin{remark}
On $0$-connective covers there is an equivalence
$K(\aC)_{>0}\simeq \bbK(\aC)_{>0}$ between this notion of
non-connective $K$-theory and the usual connective $K$-theory of
$\aC$.  In degree $0$, there an isomorphism $\pi_0
K(\aC)\cong\pi_0\bbK(\aC)$ if the underlying $\i$-category of $\aC$ is
idempotent complete. 
\end{remark}

\begin{theorem}\label{thm:localization}
Let $\aA\to\aB\to\aC$ be a sequence of DHKS-saturated Waldhausen categories with factorization such that
\[
\Ho(\N(\aA)[W^{-1}][\Sigma^{-1}])\too\Ho(\N(\aB)[W^{-1}][\Sigma^{-1}])\too\Ho(\N(\aC)[W^{-1}][\Sigma^{-1}])
\]
is a localization sequence of triangulated categories.
Then the induced map
\[
\bbK(\aA)\too\bbK(\aB)\too\bbK(\aC)
\]
is a cofiber sequence of spectra.
\end{theorem}

\begin{proof}
This follows from the natural equivalence $\bbK(-)\simeq\bbK(\N(-)[W^{-1}][\Sigma^{-1}])$ and the fact that cofiber sequence
\[
\bbK(\N(\aA)[W^{-1}][\Sigma^{-1}])\too\bbK(\N(\aB)[W^{-1}][\Sigma^{-1}])\too\bbK(\N(\aC)[W^{-1}][\Sigma^{-1}])
\]
is a cofiber sequence because $\bbK(-)$ is a localizing invariant.
\end{proof}

\subsection{Extending co-representability}

In this section, we show how to extend the co-representability of
negative $K$-theory obtained in theorem~\ref{thm:main3} to maps out of
any dualizable object, using the theory developed in
section~\ref{sec:monoidal}.  We begin with the following technical
lemma:

\begin{lemma}\label{lem:preserve}
Let $\cB$ be a small stable idempotent-complete $\i$-category.  Then
the functor given by $(-) \idemtimes \cB$ preserves equivalences,
filtered colimits, the point, and exact sequences.
\end{lemma}

\begin{proof}
It follows from the definition that $(-) \idemtimes \cB$ preserves
equivalences, filtered colimits, and the point.  The characterization
of \cite[6.3.1.16]{HAG} implies that it preserves exact sequences.
\end{proof}

We can now prove the main theorem of this section:

\begin{theorem}\label{thm:main34gen}
Let $\cB$ be a smooth and proper small stable $\infty$-category in the
sense of definitions~\ref{def:proper} and \ref{def:smooth}.  Then
$\Uloc(\cB)$ is compact in $\Motloc$ and for every small stable
$\i$-category $\cA$, we have a natural equivalence of
spectra 
\[
\Map(\Uloc(\cB),\Uloc(\cA)) \htp \bbK(\aB^{\op} \idemtimes \aA)
\]
\end{theorem}

\begin{proof}
For any small stable idempotent-complete $\i$-category $\cB$, we can
consider the functor  
\[
(-) \idemtimes \cB \colon \idemstabcat \to \idemstabcat.
\]
By lemma~\ref{lem:preserve}, the composed morphism 
\[
\idemstabcat \stackrel{(-) \idemtimes \cB}{\too} \idemstabcat \stackrel{\Uloc}{\too} \Motloc
\]
is a localizing invariant.  Thus, we obtain a commutative diagram
\[
\xymatrix{
\idemstabcat \ar[d]_{\Uloc} \ar[r]^{(-) \idemtimes \cB}
& \idemstabcat \ar[d]^{\Uloc} \\ 
\Motloc \ar[r]_{\overline{(-) \idemtimes \cB}} & \Motloc\,,
}
\]
with $\overline{(-) \idemtimes \cB}$ a colimit-preserving functor
such that 
\[
\overline{\Uloc(\cA) \idemtimes \cB} \simeq \Uloc(\cA \idemtimes \cB)\,.
\]
Now, recall from theorem~\ref{thm:dualizable} that since $\cB$ is
smooth and proper, it is also dualizable (in the symmetric monoidal
$\infty$-category $\idemstabcat$ of idempotent-complete small stable
$\infty$-categories). Therefore, we have an adjunction (on the
left) \cite[4.2.5.6]{HAG}, which induces an adjunction (on the right) 
\[
\xymatrix{
\idemstabcat \ar@<-1ex>[d]_{(-) \idemtimes \cB}
&&& \Motloc \ar@<-1ex>[d]_{\overline{(-) \idemtimes \cB}} \\ 
\idemstabcat \ar@<-1ex>[u]_{(-) \idemtimes \cB^{\op}}
&&& \Motloc \ar@<-1ex>[u]_{\overline{(-) \idemtimes \cB^{\op}}}\,, 
}
\]
with $\overline{(-) \idemtimes \cB^{\op}}$ a colimit
preserving morphism, such that 
\[
\overline{\Uloc(\cA)\idemtimes \cB^{\op}} \htp \Uloc(\cA \otimes \cB^{\op})\,.
\]
The proof now follows from the following equivalences of spectra 
\begin{eqnarray}
\Map(\Uloc(\cB), \Uloc(\cA)) & \htp
& \Map(\overline{\Uloc(\ispec^{\omega}) \idemtimes \cB}, \Uloc(\cA)) \nonumber \\
& \htp & \Map(\Uloc(\ispec^{\omega}), \overline{\Uloc(\cA) \idemtimes \cB^{\op}}) \nonumber \\
& \htp & \Map(\Uloc(\ispec^{\omega}), \Uloc(\cA \idemtimes \cB^{\op})) \nonumber \\
& \htp & \bbK(\cA \idemtimes \cB^{\op}) \nonumber \\
& \htp & \bbK(\idemfun(\aB, \Idem(\aA))) \nonumber \, .
\end{eqnarray}
Finally, since in the adjunction
$$
\xymatrix{
 \Motloc \ar@<-1ex>[d]_{\overline{-\idemtimes \cB}} \\
 \Motloc \ar@<-1ex>[u]_{\overline{-\idemtimes \cB^{\op}}}\,,
}
$$
the morphism $\overline{(-) \idemtimes \cB^{\op}}$ preserves colimits,
the object $\Uloc(\ispec^{\omega})$ is compact, and 
$\overline{\Uloc(\ispec^{\omega})\otimes \cB} \simeq \Uloc(\cB)$, we
conclude that $\Uloc(\cB)$ is compact. 
\end{proof}

\subsection{Non-connective $K$-theory of connective ring spectra}

In this section, we show that for a {\em connective} ring spectrum $R$, the
non-connective $K$-theory spectrum we associate to the category of
perfect $R$-modules has negative homotopy groups determined by the
classical non-connective $K$-theory spectrum of the ring $\pi_0 R$.

We give a proof using a model of non-connective $K$-theory
for connective ring spectra based on the construction of a
``suspension ring spectrum'' coupled with Quillen's plus construction.
We begin by recalling Wagoner's construction \cite{Wagoner} of the
non-connective $K$-theory of an ordinary ring $R$.  Given a ring $R$,
we let $\ell R$ denote the ring of locally finite (countably) infinite
matrices in $R$ --- i.e., $\bN \times \bN$ matrices such that each row
and column only has finitely many nonzero elements.  We let $mR$ denote the
finite matrices, regarded as a 2-sided ideal of $\ell R$ --- these are the
matrices with only finitely many nonzero elements.  Then we can form
the quotient ring $\mu R = \ell R / mR$, and Wagoner defines the
non-connective $K$-theory spectrum to have $n$th space
\[
K(R)_n = K_0(\mu^n R) \times BGL^+ (\mu^n R).
\]
It is known that this construction agrees with other possible
constructions of the non-connective algebraic $K$-theory spectrum of
$R$ (e.g., see \cite[\S 6]{PedersenWeibel}).

Next, we recall the generalization of this construction to connective
ring spectra.  Prior to the invention of modern notions of structured
ring spectra, May initiated the study of the algebraic $K$-theory of a
multiplicative object called an ``$A_\infty$ ring space'', which is an
$E_\infty$ space with a suitably compatible $A_\infty$ multiplication
(for a particular pair of operads) \cite{MayK, Steiner}.  The
prototype example of an $A_\infty$ ring space is $\Omega^\infty R$ for
a connective ring spectrum $R$ \cite[3.1]{MayK}.  Fiedorowicz,
Schw\"anzl, Steiner, and Vogt~\cite{noncondeloop} extended Wagoner's
constructions by defining $mR$ and $\ell R$ for $A_\infty$ ring
spaces (using the work of \cite{Steiner} to define matrices with
entries in $A_\infty$ ring spaces), and then defining $\mu R$ to be
the homotopy cofiber of the inclusion $mR \to \ell R$.  Furthermore,
they prove that there is an equivalence of spaces
\begin{eqnarray}\label{eqn:deloop}
K_0(\pi_0(R)) \times BGL^+(R) \htp \Omega(K_0(\pi_0(\mu R)) \times
BGL^+(\mu R)).
\end{eqnarray}
These constructions then allow a definition of the non-connective
algebraic $K$-theory of an $A_\infty$ ring space $R$ with spaces
\begin{eqnarray}\label{eqn:fssv}
\bbK(R)_n = K_0(\mu^n \pi_0 R) \times BGL^+ (\mu^n R).
\end{eqnarray}
This definition implies that for an $A_\infty$
ring space $R$, the natural map $R \to\pi_0 R$ induces an
isomorphism on the algebraic $K$-groups $\bbK_{-n}(R):=K_0(\mu^n\pi_0
R)$ for $n \geq 0$ \cite[1.1]{noncondeloop}.

Our approach involves constructing a variant of the ``suspension
ring'' construction that allows a construction of a non-connective
$K$-theory spectrum which agrees with the non-connective $K$-theory of
the ring space $\Omega^\infty R$ as defined in equation~\ref{eqn:fssv}
on $\pi_i$ for $i < 0$ and is equivalent to our version of the the
non-connective $K$-theory $\bbK(R):=\bbK(\widehat{R}_{\perf})$
spectrum constructed in definition~\ref{def:knonconnspec}.  Since
$\pi_0\Omega^\infty R\cong\pi_0 R$, this equivalence
implies the desired comparison.

We begin by recalling the definition of the plus construction
introduced in \cite{WaldA2}, extended to $A_\i$ ring spectra.  For
convenience, we model $A_\i$ ring spectra as EKMM $S$-algebras.  We
write $M_n R=\map_R(R^{\lor n},R^{\lor n})$ for the space of
$R$-module endomorphisms of (a cofibrant replacement of) $R^{\lor n}$,
and write $GL_n(R)\to M_n(R)$ for the full subspace of $R$-module
automorphisms of $R^{\lor n}$; that is, we have a (homotopy) pullback
of spaces
\[
\xymatrix{
GL_n(R) \ar[r] \ar[d] & M_n (R) \ar[d] \\
GL_n (\pi_0 R) \ar[r] & M_n (\pi_0 R) \cong \pi_0 M_n (R). \\
}
\]
Since $GL_n (R)$ is a topological monoid, after replacing to ensure
the inclusion of the unit is a cofibration, we can form its classifying 
space $BGL_n (R)$.  Moreover, there are natural inclusions
$GL_n(R) \to GL_{n+1}(R)$ which induce maps $BGL_n(R) \to
BGL_{n+1}(R)$.  We can form
\[
BGL(R) \cong \hocolim_n BGL_n(R).
\]
Since $\pi_1 BGL(R) \cong GL (\pi_0 R)$, we can form the plus
construction $BGL(R)^+$, and one could define the $K$-theory space to
be the infinite loop space $K_0(\pi_0 R) \times BGL(R)^+$.  The
consistency of this definition is proved in~\cite[7.1]{EKMM}, which we
restate below:

\begin{lemma}\label{lem:plusconst}
Let $R$ be a connective $A_\infty$ ring spectrum.  There is an
equivalence of infinite loop spaces
\[
\Omega^\infty K(R) \htp K_0(\pi_0 R) \times BGL(R)^+.
\]
\end{lemma}

This is consistent in the sense that a check of the definition of the
plus construction for an $A_\infty$ space \cite[\S 7]{MayK} now yields
the following proposition:

\begin{proposition}\label{prop:kcomp}
For a connective ring spectrum $R$, the connective algebraic
$K$-theory space $BGL(R)^+$ is equivalent to the algebraic
$K$-theory space $BGL^+(\Omega^\infty R)$.
\end{proposition}

We now set up analogues of the constructions of~\cite{noncondeloop}.
In order to ensure that our mapping spaces and spectra have the correct
homotopy type, we continue to work with the category of EKMM algebra
and module spectra.  Since all objects are fibrant, it then suffices
to work with cofibrant modules.  For a connective ring spectrum $R$,
in the following we let $\Map_R(x,y)$ denote the mapping spectrum
between objects $x$ and $y$ and $\map_R(x,y)$ the mapping space (which
can be computed as $\Omega^\infty \Map(x,y)$) in the category of
$R$-modules.  Moreover, when we write $R^{\lor n}$ inside a mapping
object, we will tacitly mean the wedge of a cofibrant replacement of
$R$ as an $R$-module.

\begin{definition}\label{defn:lfmat}
Let $R$ be a connective $A_\i$ ring spectrum.
We set 
\[
MR =\colim_n\Map_R(R^{\lor n},R^{\lor n}),
\]
the {\em nonunital} $A_\i$ ring spectrum of finite $R$-valued matrices.
We write $LR$ for the $A_\i$ ring spectrum of locally finite matrices,
i.e. the connective $A_\i$ ring spectrum obtained as the homotopy
pullback
\[
\xymatrix{
LR\ar[r]\ar[d] & \End_R(R^{\lor\infty})\ar[d]\\
H\ell\pi_0 R\ar[r] & H\End_{\pi_0 R}(\pi_0 R^{\lor\infty}),}
\]
where here $H$ denotes the Eilenberg-Mac Lane spectrum functor.
\end{definition}

\begin{remark}
These definitions are consistent with the those of $m$ and
$\ell$ of \cite{noncondeloop} in the case of a ringlike $A_\infty$
space of \cite{noncondeloop} --- one can construct equivalences
$\Omega^\i MR \simeq m(\Omega^\infty R)$ and $\Omega^\i
LR \simeq\ell(\Omega^\infty R)$, although we leave the details to the
interested reader, in order avoid a detailed discussion of the
technology for $A_\infty$ ring spaces.
\end{remark}

We now begin to prove the comparison theorem,
theorem~\ref{thm:nonconnKconn} below.  As explained in \cite[\S
15]{BM}, without loss of generality we can work with categories
enriched in EKMM $S$-modules as a model for spectral categories, and
we tacitly move between categories enriched in EKMM $S$-modules and
categories enriched in symmetric spectra in the following discussion.

Let $F_R$ denote the spectral category of finitely generated free
$R$-modules.  The theorem follows from
proposition~\ref{prop:suspring}, which depends on the existence of a
spectral category $\tilde{F}^\infty_R$, equipped with a homotopically
fully faithful spectral functor $F_R\to \tilde{F}^\infty_R$, whose
$\i$-category of modules has a generator $G$ such that $\pi_0$ of the
endomorphism ring spectrum of the image of G in the quotient category
$\Psi(\tilde{F}^\infty_R)/ \Psi(F_R)$ is $\mu(\pi_0\Omega^\infty R)$.
This approach to constructing analogues of $\mu (\Omega^\infty R)$ is
motivated by the explicit description of mapping spectra in the stable
quotient (see~\cite[1.3]{Drinfeld} for the dg-case and~\cite[\S 6]{BM}
for the spectral analogue) and an idea
from~\cite[6.1]{PedersenWeibel}.

We begin by giving a particular construction of such a spectral
category.  Roughly speaking, the idea is to adjoin the object
$R^{\lor\infty}$ to $F_R$ in such as way that the inclusion $F_R\to
\tilde{F}^\infty_R$ is fully faithful and $\Psi(\tilde{F}^\infty_R)$
is generated by an object $G=R^{\lor\infty}$ such that
$\pi_0\End_{\tilde{F}^\infty}(G)\cong\ell(\pi_0\Omega^\infty R)$. 

Recall that we denote by $\widehat{R}$ the category of
$R$-modules, which we can regard as a spectral category.  Let $F_R$
denote the full spectral subcategory of $\widehat{R}$ spanned by the
finite free $R$-modules $R^{\lor n}$, $n\in\mathbb{N}$, and let
$F^\infty_R$ denote the full spectral subcategory of $\widehat{R}$
spanned by the finite free $R$-modules as well as the countable wedge
$R^{\lor\infty}\simeq\colim_n R^{\lor n}$.  The inclusion gives a
fully faithful spectral functor $i \colon F_R\to F^\infty_R$.  The
spectral category $F^\infty_R$ is an intermediate construction that we
will use to construct $\tilde{F}_R^\infty$.

Write $\Psi(F_R)$ and $\Psi(F^\infty_R)$ for the presentable
stable $\i$-categories of $F_R$-modules and $F^\infty_R$-modules, and
let $i_! \colon \Psi(F_R)\to \Psi(F^\infty_R)$ denote the left
adjoint of the restriction $i^* \colon \Psi(F^\infty_R)\to\Psi(F_R)$.
Given an $R$-algebra $A$, we will also write $\Psi(A)$ for the stable
$\i$-category of $A$-modules. 

\begin{proposition}
The unit natural transformation $\Id\to i^*i_!$ is an equivalence.
\end{proposition}

\begin{proof}
This is follows from the fact that
$i_! \colon \Psi(F_R) \to \Psi(F^\infty_R)$ is fully faithful,
which in turn follows from the fact that $i$ is a fully faithful
functor of spectral categories.
\end{proof}

Since $i_!$ is fully faithful, we have an exact sequence of
presentable stable $\i$-categories 
\[
\Psi(F_R) \to \Psi(F^\infty_R) \to \aC,
\]
where $\aC$ denotes the cofiber of $i_!$.  We can regard $\aC$ as the
full subcategory of $\Psi(F^\infty_R)$ spanned by the local objects.
In mild abuse of notation, for each $0\leq n\leq\infty$, we will write
$R^{\lor n}$ for the $F^\infty_R$-module represented by $R^{\lor n}$.

\begin{proposition}
The $F_R$-module represented by any finite wedge $R^{\lor n}$ is a
compact generator of $\Psi(F_R)$ and the $F^\infty_R$-module
represented by the countably infinite wedge $R^{\lor\infty}$ is a compact
generator of $\Psi(F^\infty R)$.  In particular, we have
equivalences $\Psi(R)\simeq\Psi(\End_R(R^{\lor n}))\simeq\Psi(F_R)$ and
$\Psi(\End_R(R^{\lor\infty}))\simeq\Psi(F^\infty_R)$. 
\end{proposition}

\begin{proof}
The statement about compact generators essentially follows by
construction.  Then the $\i$-categorical version of Schwede-Shipley's
Morita theorem~\cite{SS}, \cite[\S 7.1.2]{HAG}, allows us to
characterize these categories in terms of endomorphisms of the compact
generator. 
\end{proof}

Since $i^*i_!$ is equivalent to the identity, the counit map $i_! i^*
R^{\lor\infty}\to R^{\lor\infty}$ restricts to an equivalence of
$F_R$-modules.  However, it is not an equivalence of
$F^\infty_R$-modules, since not all endomorphisms of $R^{\lor\infty}$
(e.g., the identity) factor through $i_!i^* R^{\lor\infty}$.

The following proposition is standard; we restate it for convenience.

\begin{proposition}
An $F^\infty_R$-module $M$ is in the full subcategory
$\aC\subseteq\Psi(F^\infty_R)$ spanned by the local objects if and
only if $i^*M\simeq 0$ in $\Psi(F_R)$.  Similarly, a map of
$F^\infty_R$-modules $f \colon M\to M'$ is a local equivalence if and
only if the cofiber of $f$ lies in the essential image of $i_!$. 
\end{proposition}

\begin{proof}
The first claim follows from the fact that $i^*M\simeq 0$ if and only
if for all $F_R$-modules $N$, $\Map_{F^\infty_R}(i_! N,M)\simeq 0$.
In turn, this holds if and only if for any map of $F^\infty_R$-modules
$Q\to P$ with cofiber of the form $i_! N$,
\[
\Map_{F^\infty_R}(P,M)\simeq\Map_{F^\infty_R}(Q,M).
\]
The second claim follows from the fact that, if the cofiber of $f$
lies in the essential image of $i_!$, then for any local object $L$,
$\Map_{F^\infty_R}(M',L)\simeq\Map_{F^\infty_R}(M,L)$.
\end{proof}

As a consequence, we can identify a compact generator of $\aC$.

\begin{corollary}\label{cor:compgen}
Let $G$ denote the cofiber of the counit $i_! i^* R^{\lor\infty}\to
R^{\lor\infty}$ in $\Psi(F^\infty_R)$.  Then $G$ lies in the full
subcategory $\aC\subseteq\Psi(F^\infty_R)$, i.e. $G$ is a local
object, and the map $R^{\lor\infty}\to G$ is a local equivalence.
Furthermore, $G$ is a compact generator of $\aC$.
\end{corollary}

\begin{proof}
By the previous proposition, $G$ is a local object, and the cofiber
\[
\Sigma i_!i^* R^{\lor\infty}\simeq i_!\Sigma i^* R^{\lor\infty}
\]
of $R^{\lor\infty}\to G$ is in the image of $i_!$.  $G$ is compact
because $R^{\lor \infty}$ is a compact generator of $F^\infty_R$ and
the functor $\Psi(F^\infty_R)\to\aC$ preserves compact
objects~\cite[2.9]{Marco}.
\end{proof}

This suggests that we might consider $\End_{\aC}(G)$, regarded as an
$A_\infty$ ring spectrum under composition, as an analogue of
$\hat{\mu}R$.  Note that by corollary~\ref{cor:compgen},
$\End_{\aC}(G)\simeq\Map_{F^\infty R}(R^{\lor\infty}, G)$ is 
equivalent to the cofiber (in spectra) of the map  
\[
\Map_{F^\infty_R}(R^{\lor\infty},i_!
i^*R^{\lor\infty})\to \Map_{F^\infty_R}(R^{\lor \infty}, R^{\lor \infty}).
\]
However, since $\pi_0(\Map_{F^\infty_R}(R^{\lor \infty},
R^{\lor \infty}))$ can be identified as the collection of infinite
matrices with values in $\pi_0(R)$ that have finitely many elements
per row, we need to perform a construction analogous to
definition~\ref{defn:lfmat}.

Let
\[
\pi_0\tilde{F}_R^\infty\to\pi_0 F_R^\infty
\]
denote the subcategory of $\pi_0 F_R^\infty$ consisting of those maps
\[
f\in\pi_0\Map(R^{\lor m}, R^{\lor n})\cong\Hom_{\pi_0 R}(\pi_0 R^{\lor m},\pi_0 R^{\lor n}),
\]
for $0\leq m,n\leq\infty$, which are locally finite when regarded as
elements of the group of $\pi_0 R$-valued $m \times n$-matrices.
Since the composition induces on $\pi_0$ the product of matrices and
products of locally finite matrices are locally finite, this
specification does indeed define a subcategory of $\pi_0 F_R^\infty$.
Furthermore, $\pi_0 \tilde{F}_R^{\infty}$ inherits an enrichment over
abelian groups from that of $\pi_0 F_R^\infty$.  We now perform a
categorical analogue of definition~\ref{defn:lfmat}, using the
Eilenberg-Mac Lane functor $H$ from categories enriched in abelian
groups to spectral categories~\cite[5.1.5]{SS}.

\begin{lemma}
The symmetric monoidal functor $\pi_0 \colon \Sp_{\geq
0}\to\mathrm{Ab}$ is right adjoint to the Eilenberg-MacLane spectrum
functor $H$, which is lax symmetric monoidal.  It induces a functor
\[
\pi_0 \colon \Cat_{\mathrm{Sp}_{\geq 0}}\to\Cat_{\mathrm{Ab}},
\]
from categories enriched in connective symmetric spectra to categories
enriched in abelian groups with right adjoint $H$. 
\end{lemma}

Using this we obtain a morphism of spectral categories
\[
H\pi_0\tilde{F}_R^\infty\to H\pi_0 F_R^\infty.
\]
Note that for $0\leq m,n<\infty$, the induced map of Eilenberg-Mac Lane spectra
\[
\Map_{H\pi_0\tilde{F}_R^\infty}(R^{\lor m},R^{\lor n})\to\Map_{H\pi_0 F_R^\infty}(R^{\lor m},R^{\lor n})
\]
is an equivalence, as finite matrices are locally finite.

We now define spectral categories $\tilde{F}^\infty_R$ and $\tilde{F}_R$ as the homotopy pullbacks
\[
\xymatrix{
\tilde{F}_R\ar[r]\ar[d] & \tilde{F}^\infty_R\ar[r]\ar[d] & H\pi_0\tilde{F}^\infty_R\ar[d]\\
F_R\ar[r] & F^\infty_R\ar[r] & H\pi_0 F^\infty_R}.
\]
Observe that $\tilde{F}_R$ and $\tilde{F}^\infty_R$ have the same objects
as $F_R$ and $F^\infty_R$, respectively, but
$\End_{\tilde{F}^\infty_R}(R^{\lor\infty})\simeq LR$.

\begin{proposition}
The spectral functor $\tilde{F}_R\to F_R$ is a weak equivalence of
spectral categories, and there is an equivalence of $A_\infty$ ring
spectra $\End_{\tilde{F}^\infty_R}(R^{\lor\infty})\simeq LR$.
\end{proposition}

\begin{proof}
As the functor is actually surjective on objects, it is enough to show
that it is fully faithful.  This follows from the fact that mapping
spectra in the homotopy pullback spectral category are computed as the
homotopy pullbacks of the mapping spectra.  Applying the long exact
sequence to the homotopy pullback  
\[
\xymatrix{
\Map_{\tilde{F}^\infty_R}(R^{\lor m},R^{\lor n})\ar[r]\ar[d] & \Map_{{F}^\infty_R}(R^{\lor m},R^{\lor n})\ar[d]\\
\Map_{H\pi_0\tilde{F}^\infty_R}(R^{\lor m},R^{\lor n})\ar[r] & \Map_{H\pi_0{F}^\infty_R}(R^{\lor m},R^{\lor n})}
\]
implies the desired equivalence.  A similar computation with
$m=n=\infty$ implies the second statement.
\end{proof}

The spectral functor $\tilde{F}^\infty_R\to F^\infty_R$ induces a
functor (which is not fully faithful) 
\[
\Psi(\tilde{F}_R^{\infty}) \to \Psi(F_R^{\infty}), 
\]
on $\i$-categories of modules.

Carrying out the same analysis as above, we see that the
quotient 
\[
\aC' = \Psi(\tilde{F}_R^{\infty}) / \Psi(F_R)
\]
can be described as modules over $\End_{\aC'}(G')$, where $G'$ is the
cofiber of the map
\[
i_! i^* R^{\lor \infty} \to R^{\lor \infty}
\]
(here $R^{\lor \infty}$ is regarded as an object of
$\tilde{F}_R^{\infty}$) and hence as a spectrum $\End_{\aC'}(G')$ is
equivalent to the cofiber in spectra of the map
\begin{equation}\label{eqn:cofib}
\Map_{\Psi(\tilde{F}_R^{\infty})}(R^{\lor \infty}, i_! i^*
R^{\lor \infty}) \to \End_{\Psi(\tilde{F}_R^{\infty})}(R^{\lor \infty}).
\end{equation}

\begin{lemma}\label{lem:matcomp}
There is an equivalence of rings
$\pi_0(\End_{\aC'}(G')) \simeq \pi_0(\mu\Omega^\infty R)$.
\end{lemma}

\begin{proof}
Regarding $\tilde{F}_R^{\infty}$ as a simplicial category,
$\End_{\tilde{F}_R^{\infty}}(R^{\infty})$ is (by construction) the
$A_\infty$ ring space $\ell R$.
Furthermore, we have that 
\[
\pi_0(\Map_{\Psi(\tilde{F}_R^{\infty})}(R^{\lor \infty}, i_! i^*
R^{\lor \infty})) \cong \pi_0 (i_! i^* R^{\lor \infty}) \cong m \pi_0 R
\]
and by construction
\[
\pi_0(\End_{\Psi(\tilde{F}_R^{\infty})}(R^{\lor \infty})) \cong \ell \pi_0
R.
\]
Therefore, equation~\ref{eqn:cofib} implies that as groups there is an
isomorphism  
\begin{eqnarray*}
\pi_0(\End_{\aC'}(G')) \cong \ell (\pi_0(R)) / m (\pi_0(R)) \cong \ell (\pi_0
(\Omega^\infty R)) / m (\pi_0 (\Omega^\infty R)) \\ \cong \mu \pi_0(\Omega^\infty
R) \cong \pi_0(\mu \Omega^\infty R),
\end{eqnarray*}
where the last isomorphism follows from~\cite[5.1]{noncondeloop}.
Finally, the universal property of the cofiber in spectra implies that
there is a ring structure induced on $\pi_0(\End_{\aC'}(G'))$ induced
by the ring structure on $m \pi_0(R)$ quotiented by the two-sided
ideal $\ell \pi_0 (R)$.  Inspection of $\pi_0$ shows that this
multiplication coincides with the ring structure on
$\pi_0(\End_{\aC'}(G'))$ induced by composition.
\end{proof}

Based on this, we define
\[
\hat{\mu} R \cong \End_{\aC'}(G'),
\]
using the setup described above, and we proceed to relate this suspension ring spectrum construction to an $\i$-categorical delooping.
The basic idea is that our constructions
of the suspension rings give (smaller) models of the $\i$-categorical
cone $\kcone$ from definition~\ref{def:setup} which are more closely
related to the suspension ring spectrum $\hat{\mu} R$.

\begin{proposition}\label{prop:suspring}
Let $R$ be a connective $A_\infty$ ring spectrum.  We have a natural
equivalence of spectra 
\[
\xymatrix{
K(\Psi_{\tri}(\hat{\mu} R)) \ar[r]^-{\htp} &
K((\Ind(\Psi_{\perf}(R)))^{\kappa}/ \Psi_{\perf}{R}) \\
}
\]
for any infinite cardinal $\kappa > \omega$.
\end{proposition}

\begin{proof}
For any infinite cardinal $\kappa > \omega$, there is a natural
inclusion map  
\[
\Psi_{\perf}(F_R^{\infty}) \to (\Ind(\Psi_{\perf}(R)))^{\kappa} 
\]
induced by the fact that any countable wedge of copies of $R$ is in
$(\Ind(\Psi(R)_{\perf}))^{\kappa}$, and the latter is closed under 
retracts and stable under finite colimits.  Since the inclusion
$\Psi_{\perf}(F_R) \to \Psi_{\perf}(F_R^{\infty})$ is compatible
with the (Yoneda) inclusion 
$\Psi_{\perf}(R) \to (\Ind(\Psi_{\perf}(R)))^{\kappa}$, we
have a commutative diagram
\[
\xymatrix{
\Psi_{\perf}(F_R) \ar[r] \ar[d]_-{\htp} & \Psi_{\perf}(F_R^{\infty}) \ar[d] \\
\Psi_{\perf}(R) \ar[r] & (\Ind(\Psi_{\perf}(R)))^{\kappa}. 
}
\]
Combining this with $\Psi_{\perf}(\tilde{F}_R)\to\Psi_{\perf}(\tilde{F}^\infty_R)$ we obtain the commutative diagram
\begin{equation}\label{varsusp}
\xymatrix{
\Psi_{\perf}(\tilde{F}_R) \ar[r] \ar[d]_-{\htp}
& \Psi_{\perf}(\tilde{F}_R^{\infty}) \ar[d] \\
\Psi_{\perf}(F_R) \ar[r] \ar[d]_-{\htp} & \Psi_{\perf}(F_R^{\infty}) \ar[d] \\
\Psi_{\perf}(R) \ar[r] & (\Ind(\Psi_{\perf}(R)))^{\kappa}. 
}
\end{equation}
and hence an induced composite map of quotients 
\begin{align*}
\alpha \colon \Psi_{\perf}(\tilde{F}_R^{\infty})
/ \Psi_{\perf}(F_R) \to \Psi_{\perf}(F_R^{\infty})
/ \Psi_{\perf}(F_R) \\ \to
(\Ind(\Psi_{\perf}(R)))^{\kappa} / \Psi_{\perf}(R). 
\end{align*}
By the work above, $\alpha$ can be described as a map
\[
\Psi_{\tri}(\hat{\mu} R) \to
(\Ind(\Psi_{\perf}(R)))^{\kappa} / \Psi_{\perf}(R).
\]
Finally, since $F_R^{\infty}$ has countable coproducts, the usual
Eilenberg swindle implies that $K(F_R^{\infty})$ is contractible.
We also know that $K(\Psi_{\perf}(\tilde{F}_R^{\infty})$
is contractible~\cite[6.1,6.3]{noncondeloop}.  Therefore, applying
$\Map(\uwUloc(\ispec^\omega),\,\uwUloc(-))$ to the 
commutative diagram, the fact that all of the horizontal sequences are
strict-exact allows us to apply 
theorem~\ref{thm:key2} to conclude that $\alpha$ induces an
equivalence on $K$-theory spectra.
\end{proof}

Proposition~\ref{prop:suspring} allows us finally to establish the
desired result.

\begin{theorem}\label{thm:nonconnKconn}
Let $R$ be a connective $A_\infty$ ring spectrum.  Then for $i \leq
0$, the natural map $R \to H\pi_0 R$ induces isomorphisms
$\pi_i \bbK(R) \to \pi_i \bbK(H\pi_0 R) \cong \pi_i \bbK(\pi_0R)$.
\end{theorem}

\begin{proof}
Using the proof of proposition~\ref{prop:suspring} and mimicking
definition~\ref{def:knonconnspec}, we can define a spectrum  
\[
{\bbK}^{'}(R) := \colim_n \Omega^n K((\Psi_{\tri}(\hat{\mu}^n R)).
\]
The conclusion of proposition~\ref{prop:suspring} along with
diagram~\ref{varsusp} (which implies compatibility of the structure
maps) yields an equivalence ${\bbK}^{'}(R) \htp \bbK(R)$.  By the argument
for~\cite[11.7]{Marco}, we see that we can compute the homotopy groups
of ${\bbK}^{'}(R)$ using a fibrant model that is a spectrum with $n$th
space given by the space
\[
\Omega^{\infty} K(\Psi_{\perf}(\hat{\mu}^n R)).
\]
Lastly, lemma~\ref{lem:plusconst} and lemma~\ref{lem:matcomp} implies
that there is an equivalence  
\[
\Omega^{\infty} K(\Psi_{\perf}(\hat{\mu}^n R)) \htp K_0(\mu^n \pi_0 R) \times
BGL^{+}(\hat{\mu}^n R).
\]
Therefore, for $n > 1$, $\pi_0 \Omega^{\infty}
K(\Psi_{\perf}(\hat{\mu}^n R))$ is $K_0(\mu^n \pi_0 R) =
K_{-n}(\pi_0 R)$.
\end{proof}

\section{Trace maps}\label{sec:cyclo}

In this section we apply the work of the preceding sections to give a
universal characterization of the topological Dennis trace map $K \to
THH$ \cite{Bokstedt} and the cyclotomic trace map $K \to
TC$ \cite{BokstedtHsiangMadsen}.  More generally, we identify all
natural transformations of additive functors from $K$-theory to $THH$:
they are the multiples of the topological Dennis trace map.  This
identification provides a very satisfying conceptual construction of
the cyclotomic trace map, and of course makes it clear that all known
definitions are consistent.

\subsection{$THH$ as a localization invariant}
We begin by observing that $THH$ provides a localizing invariant of
small stable $\i$-categories.  Although it is possible to do this
directly in the setting of $\i$-categories (see for instance the more
general discussion of topological chiral homology in \cite[\S
5.3]{HAG} or the constructions outlined in \cite[5.1.1]{BFN}), we use
existing constructions in the setting of spectral categories in order
to ease technical difficulties that arise in the subsequent
construction of $TR$ and $TC$.  Our basic sources for this material
are \cite{BM} and \cite{BM2}.

Recall that for a small spectral category $\aC$ we can define
$THH(\aC)$ in terms of the Hochschild-Mitchell cyclic nerve for
spectral categories \cite[\S 3]{BM}.  The cyclic nerve is defined as
the simplicial object  
\[
N^{\cyc}_{q}\aC = \bigvee \aC(c_{q-1},c_{q}) \sma \dotsb \sma
\aC(c_{0},c_{1}) \sma \aC(c_{q},c_{0}),
\]
where the sum is over the $(q+1)$-tuples $(c_{0},\dotsc,c_{q})$ of
objects of $\aC$.  This becomes a simplicial object using the usual
cyclic bar construction face and degeneracy maps: The unit
maps of $\aC$ induce the degeneracy maps, and the composition maps in
$\aC$ (along with the twist map at the end) induce the face maps.  We
denote the geometric realization as $N^{\cyc} \aC$.

The spectrum $N^{\cyc} \aC$ has the correct homotopy type only when
$\aC$ has cofibrant mapping spectra \cite[3.1]{BM}.  Since the
cofibrant objects in the Morita model structure on small spectral categories
reviewed in theorem~\ref{thm:Quillenstable} have cofibrant mapping
spectra~\cite[4.18]{Spectral}, we can define the functor
\[
THH:= N^{\cyc}\circ Q \colon \Spcat \to \Spt
\]
where $Q$ denotes the cofibrant replacement functor in the Morita 
model structure on $\Spcat$.  Since this construction preserves Morita
equivalences~\cite[5.12]{BM} (and in fact
DK-equivalences~\cite[5.9]{BM}) the functor descends to the level of
$\i$-categories.

\begin{lemma}\label{lem:thhinf}
The functor
\[
THH \colon \Spcat \to \Spt
\]
induces a functor of $\i$-categories
\[
THH \colon \stabcat \htp \N((\Spcat)^{\mathrm{c}})[W^{-1}] \to \N((\Spt)^{\mathrm{c}})[W^{-1}] \htp \ispec.
\]
\end{lemma}

This definition of $THH$ as a functor of $\i$-categories lets us
deduce the following proposition from known properties of
$THH$ in the setting of spectral categories.

\begin{proposition}
$THH$ is a localizing invariant of small stable $\i$-categories.
\end{proposition}

\begin{proof}
The cyclic bar construction commutes with filtered homotopy
colimits of spectral categories.  Furthermore, $THH(-)$ takes exact 
sequences of spectral categories to exact sequences of
spectra \cite[7.1]{BM}.  Therefore, the induced functor on
$\i$-categories is a localizing invariant. 
\end{proof}

The force of the co-representability result for algebraic $K$-theory
(theorem~\ref{thm:main2}) is that it implies, via the spectral Yoneda
lemma, the following identification of the spectrum of natural
transformations of additive functors $K \to E$.

\begin{theorem}\label{thm:nats}
Given an additive invariant $E \colon \stabcat \too \ispec$ with values in the stable $\i$-category of spectra, we have a natural equivalence
\[
\Nat(K,E) \simeq E(\ispec^{\omega}),
\]
where $\Nat(K,E)$ denotes the spectrum of natural transformations from
$K$ to $E$ as additive invariants from small stable $\i$-categories to
spectra.
\end{theorem}

\begin{proof}
By theorem~\ref{thm:main1}, we can describe the additive invariants
$K$ and $E$ as elements of $\Fun^{\L}(\Motadd, \ispec)$.
The equivalence 
\[
\Nat(\Map(\Umot(\ispec^{\omega}),-), E) \htp E(\ispec^{\omega})
\]
follows from~\ref{thm:main2} and the spectral Yoneda lemma.
\end{proof}

In particular, applying theorem~\ref{thm:nats} to $THH$ yields the
following corollary.

\begin{corollary}\label{cor:kthh}
We have an equivalence of spectra
\[
\Nat(K(-),THH(-)) \cong THH(\cS^{\omega}) \htp THH(\bbS) \htp \bbS.
\]
Passing to $\pi_0$ on both sides we obtain an isomorphism between
homotopy classes of natural transformations and
$\pi_0(\bbS) \cong \bZ$.
\end{corollary}

\subsection{The topological Dennis trace map}
Next, we want to characterize the topological Dennis trace $K \to THH$
in terms of the classification of corollary~\ref{cor:kthh}.  To do
this, we briefly recall the construction of the topological Dennis
trace map for spectral categories and verify that it descends to
provide a natural transformation of additive invariants from
$K$-theory to $THH$ on $\i$-categories.  We rely on the work
of~\cite[\S 5]{BM}.

Recall that any small spectral category $\aA$ is equivalent in the
Morita model structure to the small spectral category
$\widehat{\aA}_{\perf}$ of perfect modules.  The category
$\widehat{\aA}_{\perf}$ admits the structure of a Waldhausen category
by restriction from the spectral model structure on $\widehat{\aA}$.
As explained in \cite[\S 15]{BM}, without loss of generality we can
work with categories enriched in EKMM $S$-modules.  In this case,
because all objects are fibrant the weak equivalences of the
Waldhausen structure on $\widehat{\aA}_{\perf}$ are compatible with
the spectral enrichment in the sense of~\cite[\S 1]{BM2}.

Now, following~\cite[\S 5]{BM2}, we construct a trace using the
perspective of \cite{McCarthy2} by ``mixing'' a cyclic bar
construction and Waldhausen's $\Sdot$ construction; this is
definition~\cite[5.12]{BM2}.  Upon passage to underlying
$\i$-categories, we end up with a natural transformation of localizing
invariants.

\begin{lemma}\label{lem:dendesc}
The topological Dennis trace above induces a natural transformation of 
localizing invariants
\[
K \to THH.
\]
\end{lemma}

\begin{proof}
It is clear from the construction of the trace described above that it
descends to a natural transformation of functors of $\i$-categories $K \to THH$.
One checks on each side that the trace commutes with filtered homotopy
colimits of spectral categories, and so the result follows.
\end{proof}

As a corollary, we know that there exists an element $x \in \bZ$ such
that $x$ corresponds to the homotopy class of the topological Dennis
trace under the identification of corollary~\ref{cor:kthh}.  The
following theorem identifies this element as the unit.

\begin{theorem}\label{thm:tunit}
The topological Dennis trace is (up to homotopy) the natural
transformation given by the identity element $1 \in \pi_0(THH(\bbS))\cong\pi_0(\bbS)
\cong\bZ$. 
\end{theorem}

\begin{proof}
Given a point $\phi$ in $\Nat(K(-),THH(-))$ (a specific natural
transformation, that is), we can describe the corresponding element in
$\pi_0(\bbS)$ as the homotopy class represented by the composite
\[
\xymatrix{
\bbS \ar[r] & \Map(\Umot(\ispec^{\omega}),\Umot(\ispec^{\omega})) \htp K(\bbS)
\ar[r]^-{\phi} & THH(\bbS) \htp \bbS,
}
\]
where the first map picks out the identity map in
$\Map(\Umot(\ispec^{\omega}),\Umot(\ispec^{\omega}))$.  There is also
a classical map $i \colon \bbS \to K(\bbS)$ constructed (for instance)
as the canonical inclusion of the finite sets into finite spaces.
Waldhausen's calculations~\cite[\S 5]{WaldA2} imply that the homotopy
class of $i$ is represented by $1 \in \pi_0(K(\bbS))$.  On the other
hand, since the identity map is the unit for the multiplication on
$\pi_0(\Map(\Umot(\ispec^{\omega}),\Umot(\ispec^{\omega})))$ induced
by the composition, it must also be represented by
$1 \in \pi_0(K(\bbS))$. 
Finally, specializing to the case when $\phi$ is the topological
Dennis trace, Waldhausen \cite[5.2]{WaldA2} proves that the composite
\[
\xymatrix{
\bbS \ar[r]^i & K(\bbS) \ar[r] & THH(\bbS) \htp \bbS
}
\]
is homotopic to the identity.
\end{proof}

\subsection{$TC$ and the cyclotomic trace map}
Much of the interest in the topological Dennis trace comes from the
fact that the $THH$ spectrum comes with an additional equivariant
structure (generalizing the classical connection between the cyclic
bar construction of a space and the free loop space) which allows a
refinement into a theory called $TC$, the {\em topological cyclic homology}.
The topological Dennis trace lifts to a map
\[
K \to TC
\]
called the cyclotomic trace map \cite{BokstedtHsiangMadsen}.

We will once again apply co-representability and point-set models to
characterize the cyclotomic trace.  However, when attempting to apply
our co-representability results to $TC$, we run into certain
obstacles.  Although $TC$ is Morita invariant and satisfies
localization \cite{BM}, it does not preserve filtered colimits and is
therefore not a localizing (or additive) invariant.  Nonetheless, we
can adapt our results to characterize the cyclotomic trace in this
setting.

We begin by recalling the definition of $TC$ in the context of
spectral categories.  Our review is brief; once again, we refer the
interested reader to~\cite[\S 5]{BM2} and \cite[\S 4]{BM} for
authoritative treatment.  Fix a prime $p$.  For a spectral category 
$\aC$, we can realize $THH(\aC)$ as a cyclotomic $S^1$-spectrum.  It is convenient to use the Bokstedt model of $THH$, which is revised in
detail in \cite[\S 4]{BM}.  Since the Bokstedt model is naturally
weakly equivalent to $N^{\cyc} Q\aC$ for any spectral category
$\aC$~\cite[3.1]{BM}, we can leverage the work above.  A cyclotomic
structure is additional structure on an equivariant spectrum arising
from the self-equivalence $S^1 / H \cong S^1$ (for finite $H \subset
S^1$) that models the structure of the free loop space.

Roughly speaking, what we have is a set of compatible maps
\[
\phi^H T \to T
\]
for finite $H \subset S^1$.  The equivariant structure allows us to
consider the associated non-equivariant spectra
\[
TR^n(\aC) = THH(\aC)^{C_{p^{n-1}}},
\]
the fixed points with respect to the induced $C_{p^{n-1}}$ action.
The inclusion of fixed points and the cyclotomic structure give rise
to maps $F$ and $R$ respectively
\[
F,R\colon TR^n \to TR^{n-1}.
\]
We define $TC^{n}(\aC)$ to be the homotopy equalizer 
\[
\holim_{F,R} TR^n(\aC) \to TR^{n-1}(\aC).\]
We then have that
\[
TC(\aC) = \holim_n TC^n (\aC),
\]
where we form the homotopy limit over the maps induced by the
restriction $R$; this definition is equivalent to the one originally
given in \cite{BokstedtHsiangMadsen}.

The work of \cite[\S 5]{BM} produces a description of $THH(\aC)$ for
a spectral category $\aC$ as a cyclotomic spectrum, and hence
constructions of $TR^n$, $TC^n$, and $TC$.  Moreover, a cyclotomic
trace map $K \to TC$ is constructed which arises from compatible maps
$K \to TC^n$.  As above, we import these constructions into the
setting of $\i$-categories.  First, we have the following lemma:

\begin{lemma}
The functors
\[
TR^n, TC^n, TC \colon \Spcat \to \Spt
\]
induces functors of $\i$-categories
\[
TR^n, TC^n, TC \colon \stabcat \htp \N((\Spcat)^{\mathrm{c}})[W^{-1}] \to \N((\Spt)^{\mathrm{c}})[W^{-1}] \htp \ispec.
\]
\end{lemma}

\begin{proof}
By~\cite[3.9]{BM}, maps which induce equivalences on $THH$ induce
equivalences on $TR^n$, $TC^n$, and $TC$.  As a consequence, the
result follows from lemma~\ref{lem:thhinf}. 
\end{proof}

Next, we observe that each of the objects $TC^{n}$ provides a
localizing invariant.

\begin{proposition}
The functor $TC^{n}$ is a localizing invariant of stable
$\i$-categories with values in the stable $\i$-category of spectra.
\end{proposition}

\begin{proof}
The localization theorem of \cite[7.1]{BM} implies $\{TC^{n}(\aC)\}$
takes strict-exact sequences of spectral categories to strict exact
sequences of small stable $\i$-categories.  Thus, we need to show that
$TR^n(\aC)$ preserve filtered colimits.  We know this for $THH$, and
the result now follows inductively from consideration of the
fundamental cofibration sequence (e.g., \cite[2.1.4]{Hesselholt})
\[
THH(\aC)_{C_{p^{n-1}}} \to TR^{n}(\aC) \to TR^{n-1}(\aC)
\]
(where the left-hand term denotes the homotopy orbit space) and the
fact that homotopy orbits commute with filtered colimits.
\end{proof}

The topological Dennis trace lifts through the constructions of
$TC^{n}$, essentially by construction.  Roughly speaking
(see~\cite[5.12]{BM2} for a detailed construction), the trace is
induced by an ``inclusion of objects'' map  $\Sdot \aC \to
THH(\Sdot \aC)$, given by taking an object to its identity map in the
0-skeleton.  Since the trace lands in the fixed set with respect
to the spacewise $S^1$-action on the cyclotomic $THH$ spectrum, it is
compatible with the maps $R$ and $F$ (see
e.g., \cite[1.2]{HesselholtMadsen} for a more detailed discussion of
this).  Moreover, we can check that the trace descends to a natural
transformation of localizing invariant using lemma~\ref{lem:dendesc}. 

\begin{lemma}
The topological Dennis trace above induces a natural transformation of 
localizing invariants
\[
K \to TC^{n}.
\]
\end{lemma}

Furthermore, the cyclotomic trace $K \to TC$ provides a natural
transformation in this setting which is assembled from natural
transformations of localizing invariants.

\begin{lemma}
The natural transformations of localizing invariants
\[
K \to TC^{n}
\]
induce a natural transformation of spectrum-valued functors
\[
K \to TC^{n}.
\]
\end{lemma}

Although $TC = \holim TC^{n}$ is not itself a localizing invariant (it
does not preserve filtered colimits in general), any natural
transformation of  functors $K \to TC$
is equivalent to the data of compatible maps to each $TC^n$.
Therefore, if we consider the spectrum of natural transformations of
functors to spectra from $K \to TC$ which restrict to localizing
invariants on each component, the spectrum is given by the limit (in the $\i$-category of spectra)
\[
\lim_n \Nat(K(-),TC^n(-)).
\]
Finally, this yields the following characterization of the cyclotomic
trace.

\begin{theorem}
After $p$-completion, the set of homotopy classes of compatible
localizing invariants $\{K \to TC^{n}\}$ is isomorphic to $\bbZ_p$.
The cyclotomic trace is represented by $1 \in \bbZ_p$.
\end{theorem}

\begin{proof}
Using theorem~\ref{thm:nats} as in the proof of corollary~\ref{cor:kthh},
we see that the spectrum of natural transformations $K \to TC$ which
restrict to localizing invariants on each component can be computed as the limits (in the $\i$-category of spectra)
\[
\lim_n \Nat(K(-), TC^n(-)) \htp \lim_n TC^{n}(\mathbb{S}) =
TC(\mathbb{S}).
\]
Completing at the prime $p$, recall that
$TC(\mathbb{S}) \htp \mathbb{S} \vee \Sigma
CP^{\infty}_{-1}$~\cite[\S 1]{Rognes2}.  Since $\pi_0(\Sigma
CP^{\infty}_{-1}) = 0$, we deduce that the set of homotopy classes of
compatible invariants is $\bbZ_p$.  Furthermore, using the argument for
theorem~\ref{thm:tunit} and passing to the limit, we can identify the
class of the cyclotomic trace by understanding the homotopy class of
the composite
\[
\mathbb{S} \to K(\mathbb{S}) \to TC(\mathbb{S}) \to
THH(\mathbb{S}) \htp \mathbb{S}
\]
(after $p$-completion).  An elaboration of Waldhausen's
results~\cite[\S 5]{WaldA2} (see~\cite[\S5]{BokstedtHsiangMadsen}
or~\cite[\S 1]{Rognes2}) implies that this homotopy class is the
identity (i.e., the unit splits the trace $TC(\mathbb{S}) \to
THH(\mathbb{S})\htp \mathbb{S}$, which gives the identification of
$TC(\mathbb{S})$ above), and so using the work of
theorem~\ref{thm:tunit} we again deduce that the cyclotomic trace is
represented by the unit.
\end{proof}

\end{document}